\documentclass{tufte-handout-tai}

\usepackage{graphicx} % allow embedded images
  \setkeys{Gin}{width=\linewidth,totalheight=\textheight,keepaspectratio}
  \graphicspath{{graphics/}} % set of paths to search for images
\usepackage{amsmath,amssymb,amsthm,amsfonts,ulem,tikz}  % extended mathematics
\usepackage{booktabs} % book-quality tables
\usepackage{units}    % non-stacked fractions and better unit spacing
\usepackage{multicol} % multiple column layout facilities
\usepackage{lipsum,cancel}   % filler text
\usepackage{fancyvrb,xcolor} % extended verbatim environments
  \fvset{fontsize=\normalsize}% default font size for fancy-verbatim environments
\usepackage{tikz-cd,bbm,mathbbol}
\DeclareSymbolFontAlphabet{\mathbbl}{bbold} %let's you use \mathbbl{k} for a field k
\hypersetup{colorlinks} %puts color to hyperlinks
\usepackage{enumerate}
\usepackage{mathabx}
\usepackage{mathtools}
\usepackage{nccmath}
\usepackage{hyperref}
\hypersetup{
    colorlinks=true,    
    urlcolor=Cerulean,
    linkcolor = ForestGreen,
}
\newcounter{dummy} %so that \pageref works properly
\usepackage[absolute]{textpos}
\usetikzlibrary{decorations.pathreplacing,} %for braces with itemize
\newcommand{\tikzmark}[1]{\tikz[baseline={(#1.base)},overlay,remember picture] \node[outer sep=0pt, inner sep=0pt] (#1) {\phantom{A}};}

\newcommand{\cat}[1]{{\normalfont\textsf{#1}}}

\DeclareMathOperator{\id}{id}
\newcommand{\adj}[4]{\begin{tikzcd}[ampersand replacement=\&, column sep=4ex]
					  	   #1 \colon #2	\ar[yshift=+.6ex]{r}
					  	\& #3 \colon #4	\ar[yshift=-.4ex]{l}
					 \end{tikzcd}}

\theoremstyle{plain}
\newtheorem{thm}{Theorem}[section]

\theoremstyle{definition}
\newtheorem{defn}[thm]{Definition}

\theoremstyle{remark}
\newtheorem{ex}[thm]{Example}

\title{What is Applied \\ \noindent Category Theory?}
\author{Tai-Danae Bradley}
\publisher{Department of Mathematics\\ \noindent
CUNY Graduate Center \\ \noindent 
New York, New York\\ \noindent 
\MakeLowercase{\texttt{tbradley@gradcenter.cuny.edu}}

}

\begin{document}
\maketitlepage

\section*{What is applied category theory?}
\begin{fullwidth}
Upon hearing the phrase ``applied category theory,'' you \textit{might} be thinking either one of two thoughts:

	\begin{enumerate}[\#1]
		\item \textit{Applied category theory? Isn't that an oxymoron?}
		\item \textit{Applied category theory? What's the hoopla? Hasn't category theory always been applied?}
	\end{enumerate}

For those thinking thought \#1, I hope to convince you that the answer is \textit{No way}! It's true that category theory sometimes goes by the name of \textit{general abstract nonsense}, which might incline you to think that category theory is too pie-in-the-sky to have any impact on the ``real world.'' My hope is that these notes will convince you that that's \textit{far} from the truth!

For those thinking thought \#2, \textit{yes} it's true that ideas and results from category theory have found applications in computer science and quantum physics (not to mention pure mathematics itself), but these are not the only applications to which the word \textit{applied} in \textit{applied category theory} is being applied. So what \textit{is} applied category theory?

\vspace{.4cm}

\noindent \textit{Read on.} 
\end{fullwidth}

%%%%%%%%%%%%% SECTION 1 %%%%%%%%%%%%%%
\section*{A quick note to the reader}
\begin{fullwidth}
Before we get started, I'll mention that this document is a collection of notes I amassed while participating in the \href{http://www.appliedcategorytheory.org/school/}{2018 Applied Category Theory Adjoint School}---a wonderful online seminar that ran from January - April 2018 and culminated in a two-week workshop at the \href{http://www.lorentzcenter.nl/lc/web/2018/969/info.php3?wsid=969&venue=Oort}{Lorentz Center} in May 2018. I had a blast learning from the folks there, and I want to share some of the things I learned with anyone who's interested. (So thanks for being interested!) Later, I'll describe a couple of the research projects discussed during the workshop. Much of the information in this PDF can be found in various journal articles, blog posts, and videos of conference talks, most of which are freely available online. I've provided citations to these throughout. Here are a few other things to know:
\end{fullwidth}
	\begin{itemize}
		\item I'll assume the reader is comfortable with the basics of category theory: \textit{categories, functors} and \textit{natural transformations}. For a friendly introduction to these topics, feel free to browse through the articles from my blog \href{www.math3ma.com}{Math3ma} listed in the margin.
		\marginnote[-1cm]{For a gentle introduction to (pure) category theory, here are a few places to start:
			\begin{itemize}
				\item \href{http://www.math3ma.com/mathema/2017/1/17/what-is-category-theory-anyway}{What is Category Theory, Anyway?}
				\item \href{https://www.math3ma.com/blog/what-is-a-category-definition-and-examples}{What is a Category?}
				\item \href{https://www.math3ma.com/blog/what-is-a-functor-part-1}{What is a Functor?}
				\item \href{https://www.math3ma.com/blog/what-is-a-natural-transformation}{What is a Natural Transformation?}
			\end{itemize}
			At the first link, you'll find a list of other recommended resources for learning about category theory.
		}

		\item I'll make heavy use of hyperlinks, as I have already, and I'll also incorporate the occasional use of \textcolor{Green}{c}\textcolor{OrangeRed}{o}\textcolor{Melon}{l}\textcolor{RoyalBlue}{o}\textcolor{Magenta}{r} throughout the text. For these reasons, it's probably best to read this PDF on a computer rather than in print form. 

		\item Finally, a fair warning: I use italics \textit{a lot} (along with frequent parenthetical remarks). I also like exclamation points! And many of my sentences begin with a conjunction.

		%\item (For those who are interested...) This PDF was made using \LaTeX's beautiful tufte-handout document class!
	\end{itemize}
\begin{textblock}{0.8}[0.5,0.5](0.52,.92)
\rule{3cm}{0.4pt}\\[10pt]
\noindent \small{\textit{My gratitude goes to the participants and mentors of the 2018 ACT Workshop from whom I learned a great deal. I also thank John Baez, Joseph Hirsh, Maximilien P{\'e}roux, and Todd Trimble for providing valuable feedback on a first draft of these notes.}}
\end{textblock}

\newpage
\section*{Introduction}
\begin{fullwidth}
One of the great features of category theory, birthed in the 1940s, is that its organizing principles have been used to reshape and reformulate problems within pure mathematics, including topology, homotopy theory and algebraic geometry. Category theory has light on those problems, making them easier to solve and opening doors for new avenues of research. Historically, then, category theory has found immense application within mathematics. As John Baez \href{https://forum.azimuthproject.org/discussion/1808/lecture-2-what-is-applied-category-theory}{recently noted}, ``[category theory] \textit{was meant to be applied.}'' 

More recently, however, category theory has found applications in a wide range of disciplines outside of pure mathematics---even beyond the closely related fields of computer science and quantum physics. These disciplines include chemistry, neuroscience, systems biology, natural language processing, causality, network theory, dynamical systems, and database theory to name a few. And what do they all have in common? That's much of what current-day applied category theory is seeking to discover. In other words, the techniques, tools, and ideas of category theory are being used to identify recurring themes across these various disciplines with the purpose of making them a little more formal. And that's what the phrase \textbf{applied category theory} (ACT) refers to in these notes. As explained on the \href{http://www.appliedcategorytheory.org/workshops/}{ACT 2018 workshop webpage},
	\begin{quote}
	...we should treat the use of categorical concepts as a natural part of transferring and integrating knowledge across disciplines. The restructuring employed in applied category theory cuts through jargon, helping to elucidate common themes across disciplines. Indeed, the drive for a common language and comparison of similar structures in algebra and topology is what led to the development category theory in the first place, and recent hints show that this approach is not only useful between mathematical disciplines, but between scientific ones as well.
	\end{quote}

Of course, one of the challenges of using category theory to transfer and integrate knowledge across disciplines is making category theory itself accessible to the broader scientific audience. John Baez and Brendan Fong address this very point in their 2016 paper on electrical circuit diagrams\footnote[][2cm]{\textit{A Compositional Framework for Passive Linear Networks}, \href{https://arxiv.org/pdf/1504.05625.pdf}{https://arxiv.org/pdf/1504.05625.pdf}}: 
\end{fullwidth}

	\begin{quote}While diagrams of networks have been independently introduced in many disciplines, we do not expect formalizing these diagrams to immediately help the practitioners of these disciplines. At first the flow of information will mainly go in the other direction: by translating ideas from these disciplines into the language of modern mathematics, we can provide mathematicians with food for thought and interesting new problems to solve. We hope that in the long run mathematicians can return the favor by bringing new insights to the table.
	\end{quote} 

\begin{fullwidth}
Although their comments refer to a particular project, they can apply to the field at large, too.

% And with that in mind, this document was written primarily for two reasons:
% \begin{itemize}
% 	\item I like to tell people about things I'm excited about. And I'm excited about applied category theory!
% 		\marginnote[2cm]{Fair warning: I use italics \textit{a lot}.\\\noindent (I also use lots of parentheses.) And exclamation points!} 
% 	\item I'd love to help make the language and ideas of (applied) category theory more accessible to a broader audience. In particular, if you \textit{think} you might be interested in joining the applied category theory community, I hope this collection of notes might be helpful in some way.
% \end{itemize}

\newthought{The goal of this document} is to give a taste of applied category from a graduate student's perspective. In doing so, I'll share \textbf{two themes} and \textbf{two constructions} that appeared frequently during the ACT 2018 workshop. The math underlying these themes and constructions is not new. The newness, rather, is in how they are being \textit{applied}. To illustrate the themes and constructions, I'll also share \textbf{two examples}---two research projects in the field of ACT. The first project relates to chemistry and the second to natural language processing, though the expositions are weighted unevenly. I'll devote considerably more time on the second example since that's where my own research interests lie. And that's what's on the carte du jour! \textbf{Two themes} and \textbf{two constructions} and \textbf{two examples}, along with a few crumbs (i.e. digressions) in between. Here's the menu in more detail:
\end{fullwidth}

\newpage
\tableofcontents

\marginnote[-5cm]{Although the items are listed \textit{linearly,} they are very much intertwined. The themes motivate the constructions; the constructions embody the themes, and both the themes and the constructions come to life in the examples.}

	\marginnote{
	\begin{center}
	\includegraphics[width=!,totalheight=!,scale=0.3]{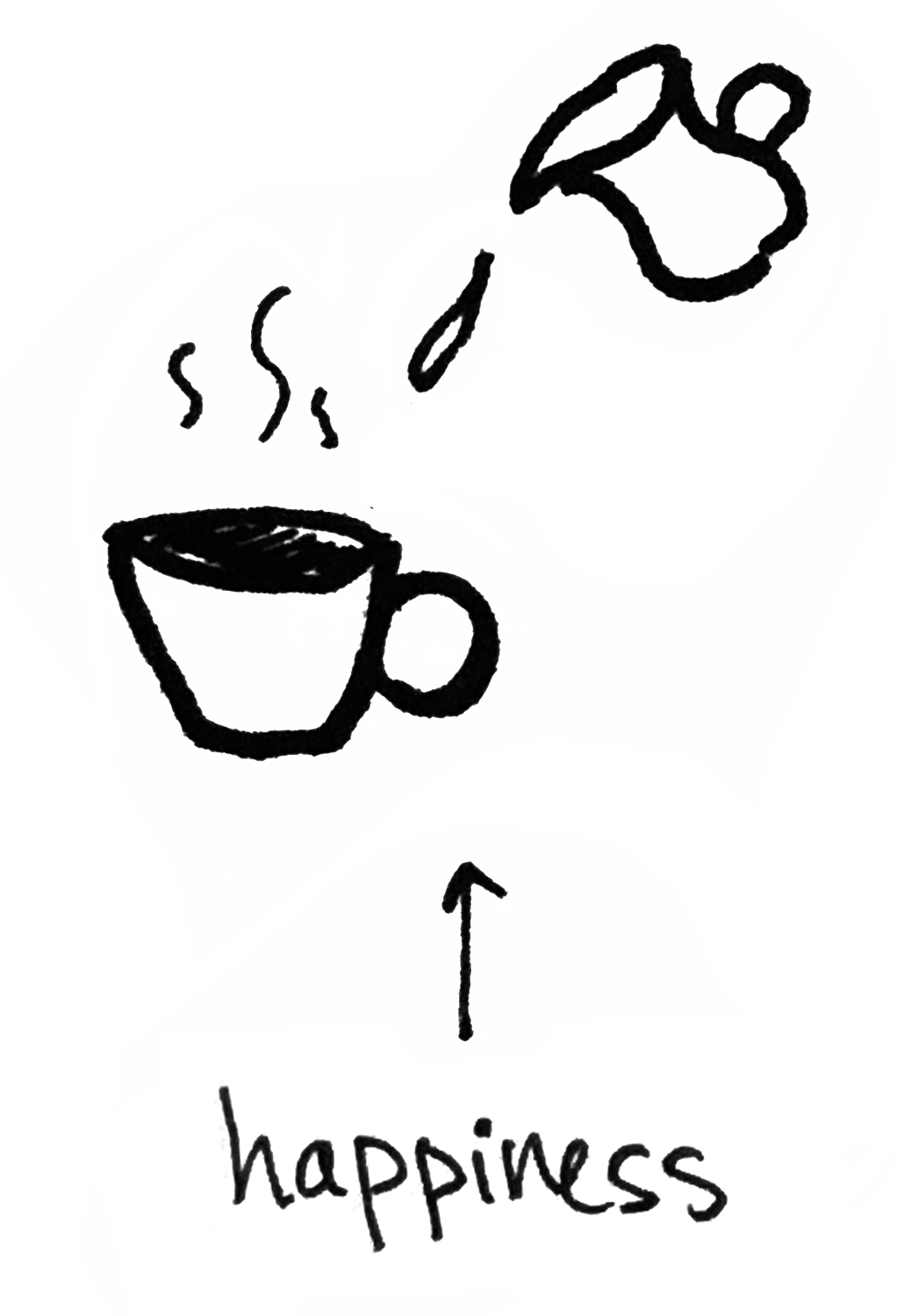}
	\end{center}}
\newthought{I like to imagine} that category theory is like a cup of black coffee, while fields outside of pure mathematics are like fresh cream. Both are lovely on their own, but blending them makes for a beverage \textit{par excellence}.
\vspace{.5cm}

\noindent I hope you'll enjoy it as much as I do!
\vspace{3cm}
% \begin{fullwidth}
% \begin{center}
% \includegraphics[width=!,totalheight=!,scale=0.13]{ACToutline.jpg}
% \end{center}
% \end{fullwidth}

\newpage
%%%%%%%%%%%%% SECTION 2 %%%%%%%%%%%%%%
\section{Two Themes}\label{sec:Themes}
Two themes that appear over and over (and \textit{over and over and over}) in applied category theory are \textbf{functorial semantics} and \textbf{compositionality.} Let's talk about the first one first.
\subsection{Functorial Semantics}\label{sec:Func}
Functorial semantics relates to the idea that a structure-preserving functor between categories
	\marginnote[-0.5cm]{The phrase ``functorial semantics'' was coined by William Lawvere.\\[0.8cm]
	\hspace{-.2cm}$\leftarrow$ This is how Lawvere \textit{defines} the word ``functor'' in his book with  Stephen  Schanuel, \href{https://www.amazon.com/Conceptual-Mathematics-First-Introduction-Categories/dp/052171916X}{Conceptual Mathematics}! It's a nice introductory text, by the way.}
\[\cat{C}\to\cat{D}\]
can be viewed as an \textit{interpretation} of $\cat{C}$ within $\cat{D}.$ 
	%\marginnote{This idea is common in computer science. Also! According to the \href{https://ncatlab.org/nlab/show/type+theory\#CategoricalSemantics}{nLab,} ``One way to look at type theory, from the point of view of a category theorist, is as a syntax for describing the construction of objects and morphisms in a category.'' a.k.a \textbf{type theory is syntax!}} 
It's often helpful to think of $\cat{C}$ as somehow encoding for \textit{syntax} while $\cat{D}$ provides \textit{semantics}. Syntax refers to rules for putting things together and semantics refers to the meaning of those things. A functor
\[\cat{syntax} \to \cat{semantics}\]
provides a way to bring the syntax to \textit{life}.

\newthought{To get a better idea} of syntax vs. semantics, think of the English language where two important features of communication are 1) grammar, which provides rules for combining words to form sentences, and 2) the actual \textit{meaning} conveyed by those words and sentences.
	\marginnote[-2cm]{I'm using English language as an analogy to illustrate syntax vs. semantics, but it's \textbf{more} than an analogy! As we'll see in Section \ref{sec:Second}, the pairings \begin{center}``grammar $\rightsquigarrow$  \textsf{syntax}'' \\ ``meanings of words $\rightsquigarrow$ \textsf{semantics}''\end{center} become quite literal in applied category theory!} 
Grammar is the syntax, and the meaning is the semantics.
	\begin{center}
		grammar $\rightsquigarrow$ \cat{syntax} 
		\qquad \qquad meaning  $\rightsquigarrow$ \cat{semantics}.
	\end{center}
Of course, neither is useful on their own. For instance, it's easy to come up with a sentence that is grammatically correct and yet has no meaning. That's the whole point behind MadLibs!
	\marginnote[-0.5cm]{
	\includegraphics{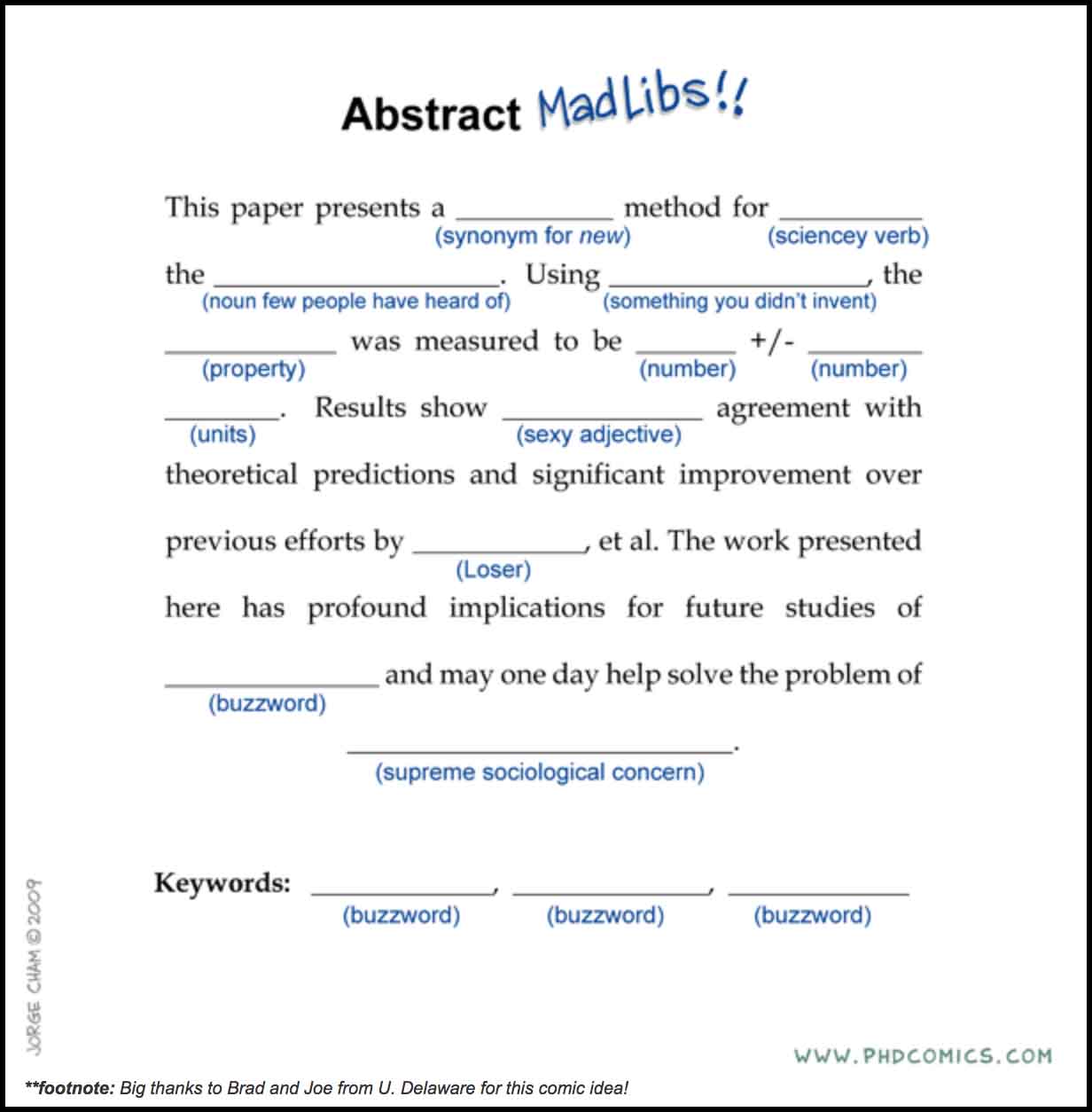}}
As another example, here's a sentence attributed to linguist \href{https://en.wikipedia.org/wiki/Colorless_green_ideas_sleep_furiously}{Noam Chomsky}:
\begin{quote}
\textit{Colorless green ideas sleep furiously.}
\end{quote}
It is grammatically correct, yet it has no meaning. The point here is that \textit{syntax vs. semantics} is nothing new. So when an applied category theorist wants to model some phenomena in the ``real world,'' don't be surprised if their model is ultimately a functor from a syntax category to a semantics category! 

\newpage

% Let me give you a very silly, non-math example. Suppose that you, like Dr. Frankenstein, wanted to make a full-grown human in a laboratory. You'd probably need 1) a skeleton and 2) some flesh. The skeleton provides the underlying structure and tells you how the little pieces fits together to make a working whole. A skeleton, however, does \textit{not} have life! You only get a human once you stick some flesh on those bones.
% 	\begin{center}
% 		\textbf{syntax} $\rightsquigarrow$ skeleton
% 		\qquad \qquad \textbf{semantics} $\rightsquigarrow$ flesh.
% 	\end{center}
% \textit{But} you need a way to do that. That is, you need a process for putting the flesh onto the skeleton in a way that's consistent with skeleton's structure\footnote{This process might also require some electricity and a maniacal laugh.}---i.e. so that it fits it like a glove. That process is precisely the morphism skeleton $\to$ flesh.

% The functor then assigns \textit{meaning} to the syntax. It provides a way to \textit{model} the syntax by sticking meaning to it.
% \begin{center}
% 		$F\colon \mathsf{syntax}  \to \mathsf{semantics}.$
% \end{center}

\subsection*{A small-ish digression...}
Even though the idea goes by the fancy name of \textit{functorial semantics}, it is \textit{not} just a ``category theory thing.'' Mind if I digress for a while to elaborate on this?\marginnote[-0.25cm]{I'll take your silence as a No.} 

\newthought{If you know a little bit} about \textit{groups}, then you've seen functorial semantics in action before! How so? A group is a set endowed with some extra structure, though that \href{https://twitter.com/evelynjlamb/status/968618077003825152}{tells us nothing} about why groups are useful. It's better to think of a group as encoding for some kind of action or transformation.
	\marginnote{Group elements are like verbs. They DO stuff! For more on this notion from a categorical perspective, check out the article \href{http://www.math3ma.com/mathema/2017/2/16/group-elements-categorically}{Group Elements, Categorically} on Math3ma.}
And this is why group representations are so great! A \textit{group representation} provides a way to view your \textit{abstract} group elements as \textit{concrete} linear transformations of some vector space. Explicitly, given a vector space $V$, a \textbf{group representation} is a group homomorphism from $G$ to $\text{Aut}(V)$, the group of all automorphisms of $V$
	\marginnote{If we replace $\text{Aut}(V)$ by $\text{Aut}(X)$ for some set $X$ (i.e. the group of automorphisms, i.e. bijections, on $X$), then a group homomorphism $G\to \text{Aut}(X)$ is precisely a group action on $X.$}
\[G\to \text{Aut}(V)\]
It assigns to each group element a linear isomorphism $V\to V.$

As a quick example, suppose our group is $D_3$, the dihedral group of order $6$, which is the group of symmetries of an equilateral triangle. If we were to look at a presentation of the group,
\[D_3=\langle r,s \:|\: r^3=s^2=rsrs=1 \rangle\]
it might not seem to have anything to do with triangles. Fortunately, a representation of $D_3$ makes the connection clearer by assigning to each group element $r$ and $s$ a linear transformation of the real plane. Specifically, the standard representation of $D_3\to \text{Aut}(\mathbb{R}^2)$ assigns to each of $r$ and $s$ an invertible $2\times 2$ matrix with real entries:
	\marginnote{
	\includegraphics{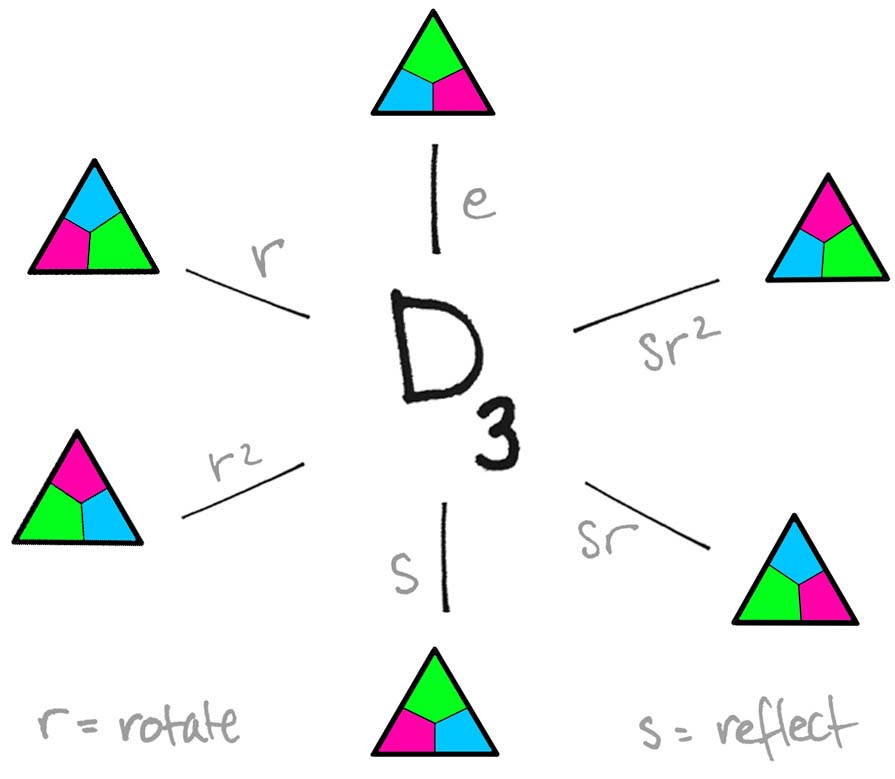}
	}
\[r\mapsto R=\begin{bmatrix}\cos(2\pi/3) & -\sin(2\pi/3)\\ \sin(2\pi/3) & \cos(2\pi/3)\end{bmatrix}
\qquad\qquad
s\mapsto S=\begin{bmatrix}1 & 0\\ 0 & -1\end{bmatrix}
\]
Here $R$ is a rotation by $60^\circ$ while $S$ is reflection across the $x$-axis. Moreover $R^3$ and $S^2$ and $RSRS$ each are equal to the $2\times 2$ identity matrix, which is exactly what we would expect: rotating an equilateral triangle by one full revolution leaves it unchanged, as does reflecting it twice in a row, and so on.

More generally then, we can think of a group $G$ as providing the \textit{syntax} while automorphisms $\text{Aut}(V)$ provide the \textit{semantics}
\[G\rightsquigarrow \cat{syntax} \qquad \qquad \text{Aut}(V)\rightsquigarrow \cat{semantics}\]
So a group representation is like a (structure-preserving) morphism
\[\cat{syntax}\to\cat{semantics}\]
In fact... it's not \textit{like} that. It \textbf{IS} that. If we view both the groups $G$ and $\text{Aut}(V)$ as one-object\footnote{Every group $G$ gives rise to a category having a single object $\bullet$ (the group itself) and a morphism $\bullet\overset{g}{\to}\bullet$ for each group element $g\in G$. Composition is given by the group operation.} categories,
	\marginnote[1cm]{Here's another example I can't resist sharing: \textbf{operads!} If you're not familiar with \href{http://www.math3ma.com/mathema/2017/10/23/what-is-an-operad-part-1}{operads}, just know that this is a souped-up version of the group theory example. If you \textit{are} familiar with operads, then you know this is the souped-up version of the group theory example. 

	An operad is an example of syntax, while an algebra over that operad provides the semantics.  For example, given a vector space $V,$ an operad homomorphism from the [commutative, associative, Lie, Poisson,...] operad  to the endomorphism operad on $V$ \textbf{IS} a [commutative, associative, Lie, Poisson,...]-algebra! That is, the structure-preserving homomorphism provides an \textit{interpretation} of each abstract $n$-ary operation as a actual, concrete operation $V^{\otimes n}\to V$.}
then a group representation \[G\to \text{Aut}(V)\] \textbf{IS} a functor from syntax to semantics. That's because every group homomorphism \textit{is} a functor when the groups are viewed as one-object categories! So although \textit{\textbf{functor}ial semantics} has the word ``functor'' in it, don't think that the idea behind it is unique to category theory. Indeed, representation theory capitalizes on the relationship between syntax and semantics: a \textit{representation} assigns to an abstract algebraic gadget (the syntax) some concrete meaning (the semantics).

\newthought{I could end our digression} here, but I'd like to share one more instance of functorial semantics at work in pure mathematics. The next few examples involve monoids and monoidal categories, so I'll assume you are familiar with those words. If you are \textit{not} familiar with those words, don't fret---you're in luck! Section \ref{sec:Monoidal} is all about monoids and monoidal categories, so feel free to read that section first then come back here. In either case, let's proceed with another neat example of functorial semantics in action:
\marginnote[1.5cm]{This next comment is really digressing from the digression, but: I also like to think of simplicial sets as an instance of functorial semantics. A simplicial set $X$ is a bit like syntax, while a topological space is like semantics. \href{https://ncatlab.org/nlab/show/geometric+realization}{Geometric realization} $X\mapsto |X|$ provides a map from one to the other.} 
% \begin{itemize}
		% \item (an example from category theory) Consider $\mathbb{N}$ viewed as a monoid under multiplication with unit the number 1. This gives rise to a category with a single object $\bullet$ (which we think of as $\mathbb{N}$ itself) and a morphism $\bullet\overset{n}{\to}\bullet$ for each natural number $n.$ Composition is given by multiplication, that is, the composite of $\bullet\overset{n}{\to}\bullet\overset{m}{\to}$ is $\bullet\overset{nm}{\to}\bullet$. Now define a functor $\bullet\to \cat{Set}$ that sends $\bullet$ to the \textit{set} $\mathbb{N}$ and sends each morphism $\bullet\overset{n}{\to}\bullet$ to the function $f\colon \mathbb{N}\to\mathbb{N}$ given by multiplication by $n$. In other words, $f(x)=nx$. This functor is \textit{one interpretation} of the monoid $\bullet$ within the category $\cat{Set}$. \tai{This actually doesn't make sense. What?}
		% \item 

\vspace{0.5cm}
\noindent \underline{Example}: \textit{a \textbf{monoid} is the image of a functor from a certain syntax category to a certain semantics category.}
\vspace{0.5cm}

\noindent More specifically,\footnote[][1cm]{\label{fn:Fc} A functor $F\colon \cat{C}\to\cat{D}$ between monoidal categories is called \textbf{lax monoidal} if for every pair of objects $c,c'$ in $\cat{C}$ there is a morphism \[Fc\otimes Fc'\to F(c\otimes c')\] (which assembles into a natural transformation.) It's called \textbf{strong monoidal} if $Fc\otimes Fc'\cong F(c\otimes c')$, and it's called \textbf{strict monoidal} if $Fc\otimes Fc' = F(c\otimes c')$.}
\begin{center}
\refstepcounter{dummy}
\label{fig:monfun1}
\includegraphics[width=!,totalheight=!,scale=0.1]{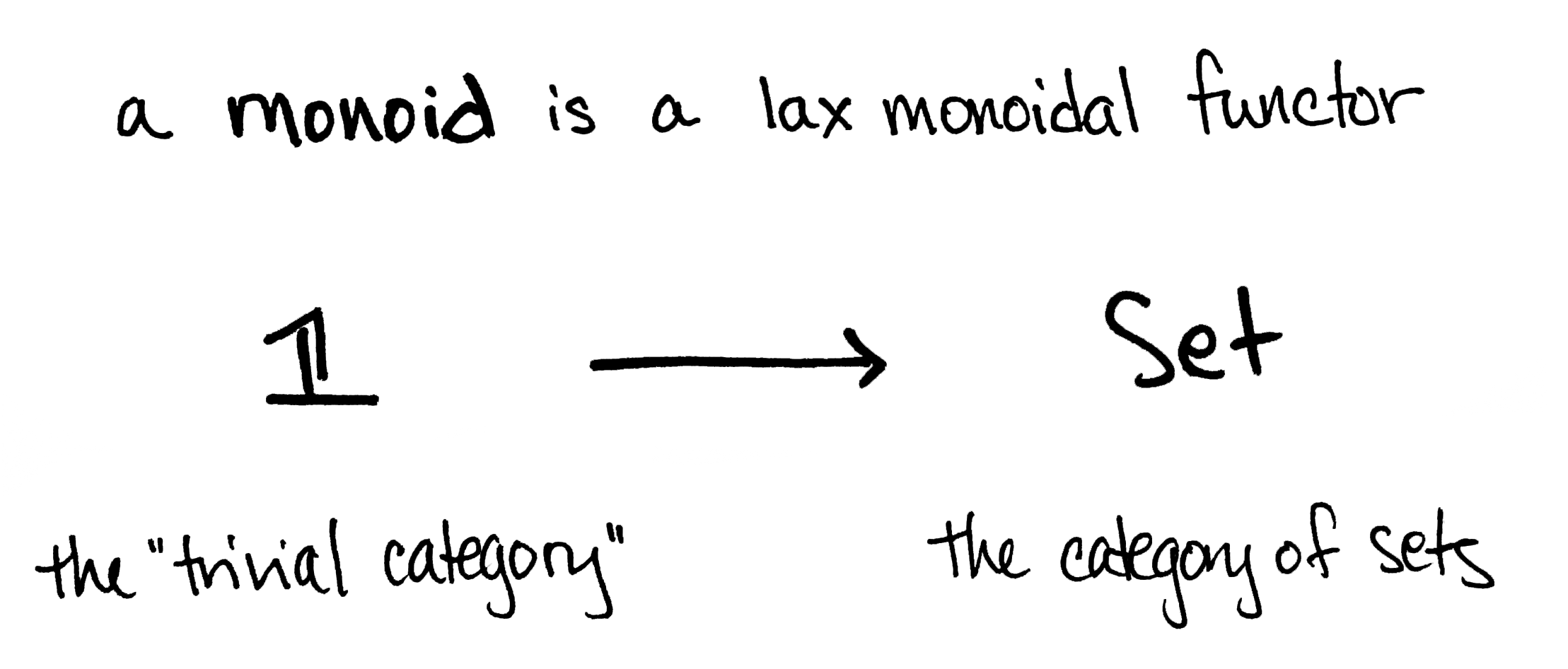}
\end{center}
Here I'm viewing both $\mathbb{1}$ and $\cat{Set}$ as monoidal categories. The symbol $\mathbb{1}$ is meant to represent the category with one object and only one morphism (the identity), which we can view as a monoidal category $(\mathbb{1},\otimes,\mathbb{1})$ in exactly one way. The category of sets has a monoidal structure given by the Cartesian product with the set containing one element, denoted $\{\ast\}$, as the monoidal unit. \textit{Technically} then,
\begin{center}
\includegraphics[width=!,totalheight=!,scale=0.1]{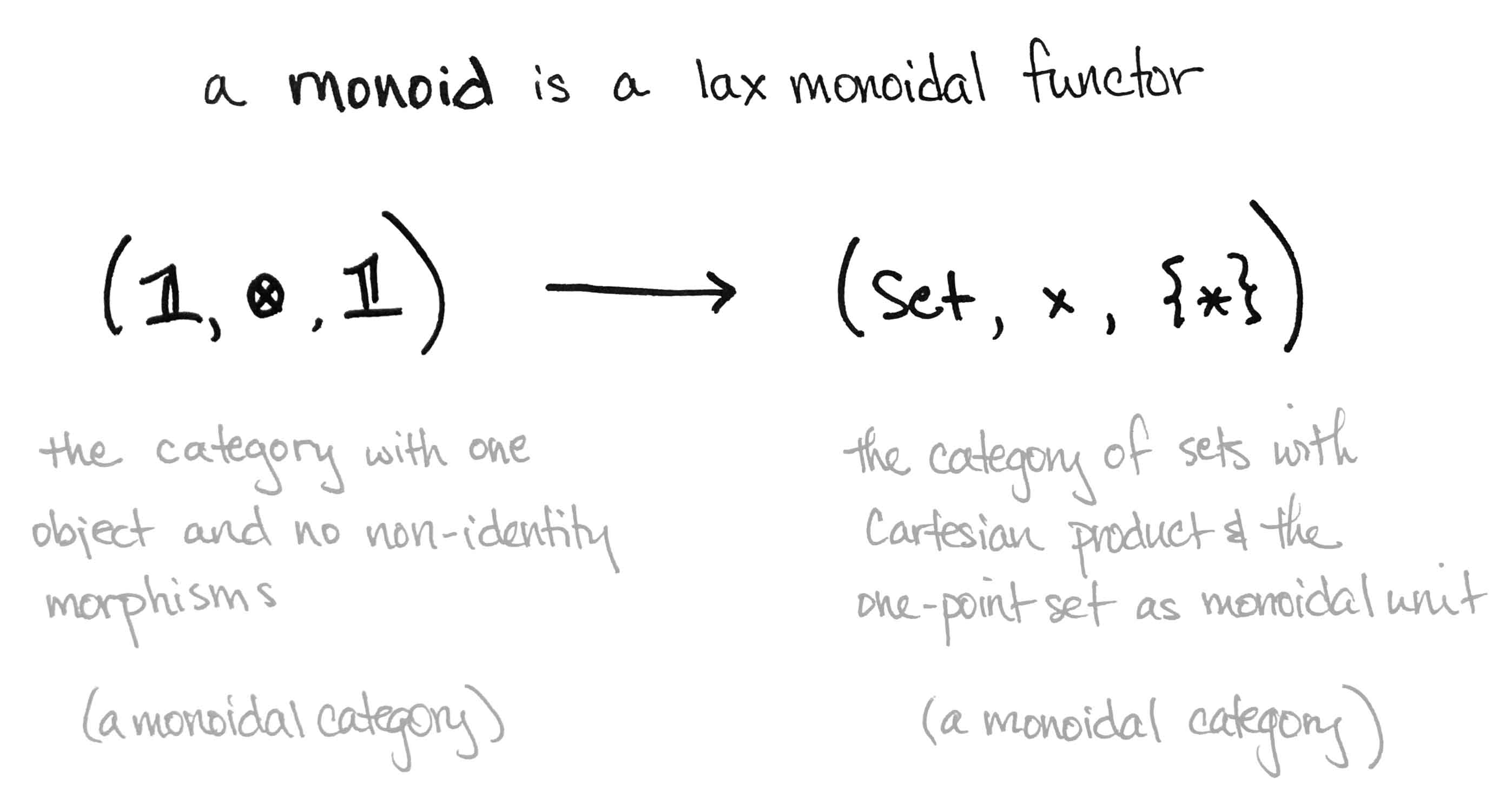}
\end{center}
	\marginnote[-4cm]{In footnote \ref{fn:Fc} you'll notice I dropped parentheses and wrote $Fc$ rather than $F(c)$. The reason for my preference is \textit{categorical!} Let me explain by saying a few words about sets: Did you know that an element $x$ in a set $X$ is the same thing as a function $\{\ast\}\to X$, where $\{\ast\}$ denotes the one-element set? It's true. The function $\{\ast\}\to X$ is uniquely determined by where it sends that one point $\ast$. So since an element $x$ is the same thing as an arrow $\{\ast\}\to X$, we might as well label \textit{that} arrow by $x,$
	\[x\colon \{\ast\}\to X \]
	Now if we have another function $f\colon X\to Y,$ then an element $f(x)\in Y$ is \textit{precisely} the composition
	\[\{\ast\}\overset{x}{\to} X\overset{f}{\to}Y\]
	That is,
	\[f(x)=f\circ x =fx\]
	where on the right hand side, I've omitted the composition symbol $\circ$ because it's cleaner. So there you have it! An element $f(x)\in Y$ is the same as a function $fx$. And since categorically-minded folks (such as you and I) prefer arrows over elements (Because of \href{http://www.math3ma.com/mathema/2017/8/30/the-yoneda-perspective}{Yoneda}. Also, we might be \href{https://arxiv.org/pdf/1212.6543.pdf}{rethinking set theory}.), the notation $fx$---and more generally, $Fc$ as above---is preferred. By the way, this is all related to Lawvere's \textbf{philosophy of generalized elements}, which is the idea that a morphism $A\to B$ is really an ``$A$-shaped element in $B$.'' For some examples, check out the articles \href{http://www.math3ma.com/mathema/2018/1/10/a-diagram-is-a-functor}{``A Diagram is a Functor''} as well as \href{http://www.math3ma.com/mathema/2017/9/6/the-yoneda-embedding}{``The Yoneda Embedding''} on Math3ma. Generalized elements are closely related to functorial semantics, so both links are worth a read!} 
Why is this true? First observe that a functor $F\colon \mathbb{1}\to\cat{Set}$ picks out a set, $F\mathbb{1}:=M$. And the data of a \textit{lax monoidal} functor $F\colon (\mathbb{1},\otimes,\mathbb{1}) \to (\cat{Set},\times,\{\ast\})$ consists of a morphism
\[\bullet\colon F\mathbb{1}\times F\mathbb{1}\to F(\mathbb{1}\otimes \mathbb{1}) \qquad \text{i.e.}\qquad \bullet\colon M\times M\to M\]
along with a morphism
\[1\colon \{\ast\}\to F\mathbb{1}\qquad \text{i.e.}\qquad 1\colon \{\ast\}\to M\]
both of which are required to fit into some commuting diagrams. I won't write them here, but one diagram says ``$\bullet$ is associative'' and the other diagram says, ``$1$ serves as an identity for $\bullet.$'' In summary, the data of a lax monoidal functor $F\colon (\mathbb{1},\otimes,\mathbb{1}) \to (\cat{Set},\times,\{\ast\})$ are
	\begin{enumerate}[i)]
		\item a set $M$
		\item an associative binary operation $\bullet\colon M\times M\to M$
		\item a special element $1:=1(\{\ast\})\in M$ that serves as a ``multiplicative identity'' for $\bullet$.
	\end{enumerate} 
This triple $(M,\bullet,1)$ is \textit{precisely} a monoid! Or to borrow from Lawvere's terminology, the functor (equivalently, the monoid)  is one \textit{interpretation} of the category $\mathbb{1}$ in the $\cat{Set}$. Interestingly, $\mathbb{1}$ may be interpreted in other categories as well. This leads to other familiar monoidal structures. Indeed, if we replace $(\cat{Set},\times,\{\ast\})$ by any monoidal category $(\cat{C},\otimes,1)$, then a lax monoidal functor 
\[(\mathbb{1},\otimes,\mathbb{1})\to(\cat{C},\otimes,1) \]
is a \textbf{monoid in the category $\cat{C}$}. Sometimes this monoid goes by a familiar name. Here are some examples.\newpage
\begin{itemize}
	\item [\textsc{\textbf{1. Topological Monoid.}}] Let $(\cat{Top},\times,\ast)$ denote the category of topological spaces and continuous functions, 
		\marginnote[3cm]{The idea that ``a monoid in $\cat{C}$ is a lax monoidal functor $\mathbb{1}\to \cat{C}$'' is \textbf{completely analogous} to claim that ``an element in $X$ is a function $\{\ast\}\to X$'' made in the margin on the previous page. In both cases, with have two objects $A$ and $B$ of the same kind (monoidal categories on this page; sets on the previous page) together with a structure-preserving map $A\to B$. 

		(Caveat: a lax monoidal functor is only \textit{somewhat} structure-preserving. That's why it's called \textit{lax}. And a function is vacuously structure-preserving since sets don't have any structure! But I digress...) 

		In both cases the object $A$ is trivial (technically, \href{https://en.wikipedia.org/wiki/Initial_and_terminal_objects}{terminal})---it's just a \textit{point,} so to speak. And in both cases the arrow provides an interpretation of that point within the context of $B$. 

		To phrase it another way, we are \textit{probing} $B$ with a point-shaped object. In the case when $B$ is a set, probing it with a point will pick out an element. In the case when $B$ is a monoidal category, probing it with a point will pick out a monoid!
		\begin{center}\refstepcounter{dummy} \label{fig:probe}
		\includegraphics[width=!,totalheight=!,scale=0.4]{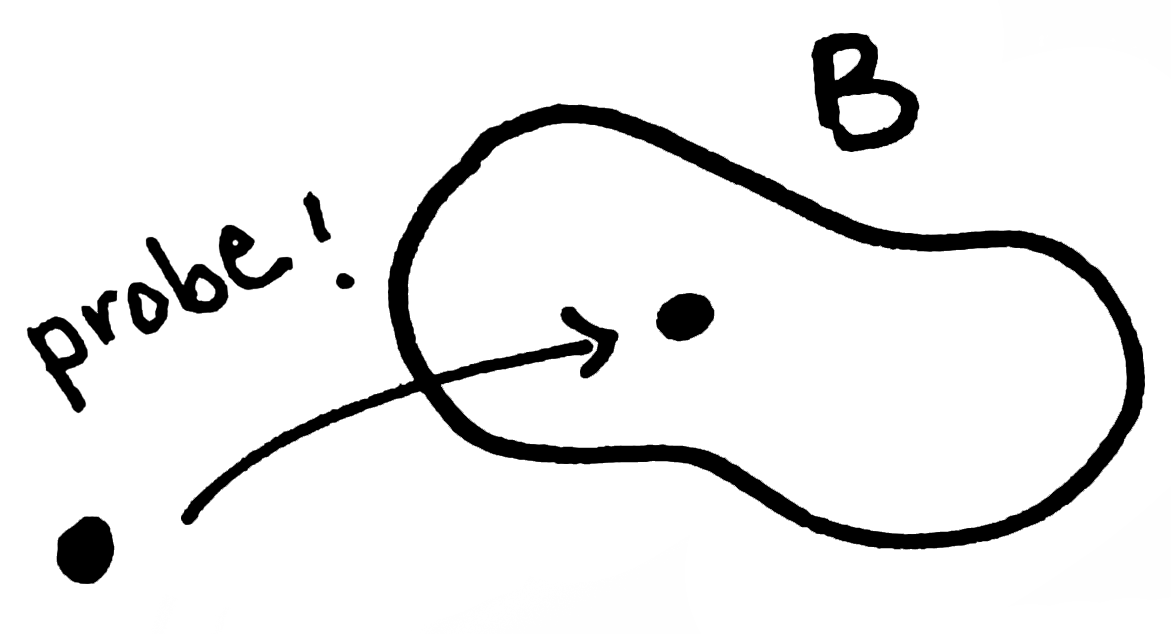}
		\end{center}}
	viewed as a monoidal category with the Cartesian product $\times$, with the one-point space $\ast$ as monoidal unit. A lax monoidal functor
		\[(\mathbb{1},\otimes,\mathbb{1}) \to (\cat{Top},\times,\ast)\] 
	is a topological monoid. That is, \textbf{a topological monoid is a monoid in the category of topological spaces.}

	\item [\textsc{\textbf{2. Ring.}}] Let $(\cat{AbGroup},\otimes,\mathbb{Z})$ denote the category of abelian groups (viewed as $\mathbb{Z}$-modules) and abelian group homomorphisms, viewed as a monoidal category with the tensor product $\otimes$, with the integers $\mathbb{Z}$ as monoidal unit. A lax monoidal functor
		\[(\mathbb{1},\otimes,\mathbb{1}) \to (\cat{AbGroup},\otimes,\mathbb{Z})\]
	is a ring (with unit). That is, \textbf{a ring is a monoid in the category of abelian groups.}
				%YES. IT'S TRUE. See Proposition 3.3 on the nLab: https://ncatlab.org/nlab/show/tensor+product+of+abelian+groups 
				%\tai{This is clear bc a ring $R$ is an abelian group under $+$ AND a monoid-with-unit under $\cdot.$ But $+$ and $\cdot$ must obey distributive property. Is \textit{that} part of the data? I think so. Probably comes from fact that morphisms aren't just functions, they are homomorphisms. Same for algebra below. Also, Emily doesn't say anything about this, so...}
	\item [\textsc{\textbf{3. Algebra.}}]  Let $(\cat{FVect},\otimes,\mathbbl{k})$ denote the category of finite-dimensional vector spaces over a field $\mathbbl{k}$ and linear maps, viewed as a monoidal category with the tensor product $\otimes$, with $\mathbbl{k}$ as monoidal unit. A lax monoidal functor 
		\[(\mathbb{1},\otimes,\mathbb{1}) \to (\cat{FVect},\otimes,\mathbbl{k})\]
	is an algebra (with unit). That is, \textbf{an algebra is a monoid in the category of vector spaces.}

	\item [\textsc{\textbf{4. Monad.}}]	Let $\cat{C}$ be a category and let $\cat{End}_\cat{C}$ denote the category whose objects are functors $\cat{C}\to\cat{C}$ and whose morphisms are natural transformations. (So $\cat{End}_\cat{C}$ is the category of \textit{endofunctors} on $\cat{C}$.) Note that $\cat{End}_\cat{C}$ can be given the structure of a monoidal category: the monoidal product is composition of functors (i.e. if $F,G$ are objects in $\cat{End}_\cat{C}$, then the monoidal product of $F$ and $G$ is $F\circ G$), and the monoidal unit is the identity functor $1_{\cat{C}}$ on $\cat{C}$ (i.e. $1_{\cat{C}}$ assigns each object and morphism in $\cat{C}$ to itself). Then a lax monoidal functor
		\[(\mathbb{1},\otimes,\mathbb{1}) \to (\cat{End}_\cat{C},\circ,1_{\cat{C}})\]
	is a monad. That is, \textbf{a monad is a monoid in the category of endofunctors on $\cat{C}$.}
\end{itemize}
You'll notice that in each of these examples, a change in the semantics category $\cat{C}$ gives rise to a different interpretation of $\mathbb{1}$, which served as our syntax category. Pretty neat, right? For more details on the examples, see Emily Riehl's \textit{\href{http://www.math.jhu.edu/~eriehl/context.pdf}{Category Theory in Context}} Definitions 1.6.3 and 5.1.1.

\newthought{That's the idea behind} functorial semantics. Now, how is it used in applied category theory? We'll see the answer when we look at the two examples---two research projects from the field---one from chemistry and one from natural language processing. In both examples, the \textit{key} is the existence of a (structure-preserving) functor from a syntax category to a semantics category. Here's a sneak preview: 
\vspace{0.3cm}
\begin{center}
\refstepcounter{dummy}
\label{fig:4cats}
\includegraphics[width=!,totalheight=!,scale=0.45]{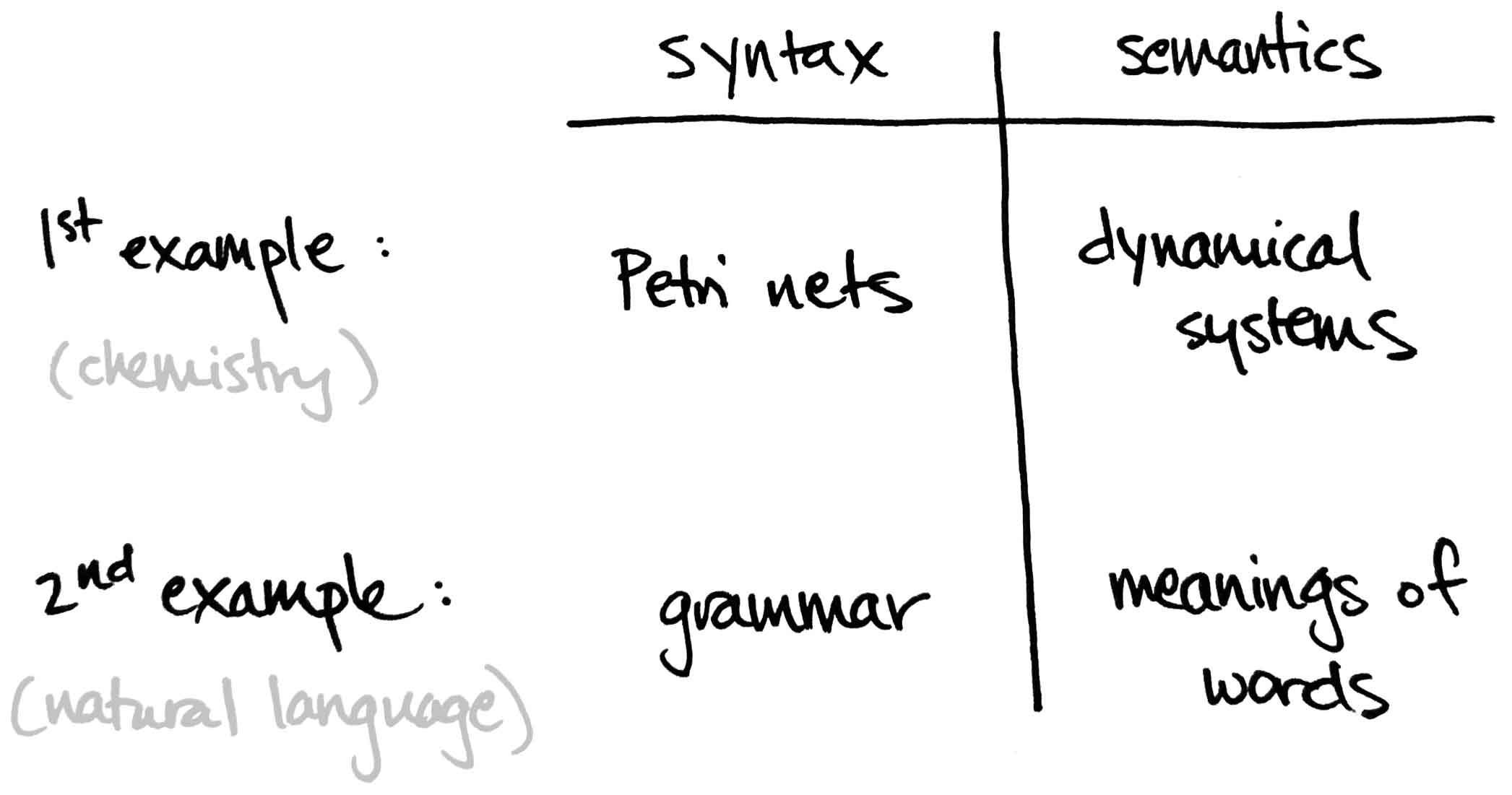}
\end{center}
\vspace{0.1cm}
% \begin{table}[h!]
% \centering
% \begin{tabular}{c|c|c|}
% \cline{2-3}
%  & \textsf{syntax} & \textsf{symantics} \\ \cline{2-3} 
% \begin{tabular}[c]{@{}c@{}}first example\\ \textcolor{gray}{(chemistry)} \end{tabular} & Petri nets & dynamical systems \\ \cline{2-3} 
% \multicolumn{1}{l|}{\begin{tabular}[c]{@{}l@{}}second example\\ \textcolor{gray}{(natural language)}\end{tabular}} & \multicolumn{1}{l|}{grammar} & \multicolumn{1}{l|}{meanings of words} \\ \cline{2-3}
% \end{tabular}
% \end{table} 
% \vspace{0.4cm}
In Section \ref{sec:First} we'll see how the behavior of a chemical reaction network is modeled by a functor
\[\mathsf{Petri\;nets}  \to  \mathsf{dynamical\;systems}\]
as shown in \href{https://arxiv.org/abs/1704.02051}{``A Compositional Framework for Reaction Networks''} by John Baez and Blake Pollard. In Section \ref{sec:Second}, we'll see how a model for natural language can be described by a functor
\[\mathsf{grammar}  \to \mathsf{meanings\;of\;words}\]
via the work of Bob Coecke, Mehrnoosh Sadrzadeh, and Stephen Clark in \href{https://arxiv.org/abs/1003.4394}{``Mathematical Foundations for a Compositional Distributional Model of Meaning''}. And perhaps you're wondering, ``How do $\mathsf{Petri\;nets}$, $\mathsf{dynamical\;systems}$, $\mathsf{grammar}$, and $\mathsf{meanings\;of\;words}$ form categories? And what's a Petri net, anyway?'' We'll answer these questions in the pages to come, but first I'd like to introduce another important theme in applied category theory: compositionality. 

% \vspace{0.6cm}
% \begin{center}
% \textit{The next section won't be as long-winded as this one. I promise!}
% \end{center}

\newpage
\subsection{Compositionality}\label{sec:Comp}
\textbf{Compositionality}, also known as the principal of compositionality, also known as Frege's principle\marginnote{Frege as in \textit{Gottlob Frege.}}, is the idea that the meaning of a complex \textcolor{RubineRed}{expression} is determined by
	\begin{enumerate}
		\item the meanings of its constituent \textcolor{BurntOrange}{parts}, and 
		\item the rules for how those parts are combined.
	\end{enumerate}
Or, as succinctly stated on the homepage of the \href{http://www.compositionality-journal.org/}{brand new journal}
of applied-category-theory-and-related-fields, 
	\begin{quote}
	\textbf{compositionality} describes and quantifies how complex things can be assembled out of simpler parts.
	\end{quote}
As it turns out, the name of that journal is itself \textit{Compositionality},\footnote{A \href{https://golem.ph.utexas.edu/category/2018/05/compositionality.html\#c053932}{contending title} was \textit{Applied Category Theory}, but in the end \textit{Compositionality} had the most votes.} which hints at the importance of this concept within the field.
	% \begin{center}
	% \begin{figure}[h]
	% \includegraphics{compositionality.png}
	% {\caption*{\textcolor{gray}{A snapshot of the banner on the new journal's \href{http://www.compositionality-journal.org/}{homepage.}}}}
	% \end{figure}
	% \end{center}

In Section \ref{sec:First}, which contains our example from chemistry, the complex \textcolor{RubineRed}{expression} will be a network---a big complicated directed multigraph, if you like. Its constituent \textcolor{BurntOrange}{parts} are simply smaller chunks of the network.

\begin{fullwidth}
\begin{center}
\includegraphics[width=!,totalheight=!,scale=0.14]{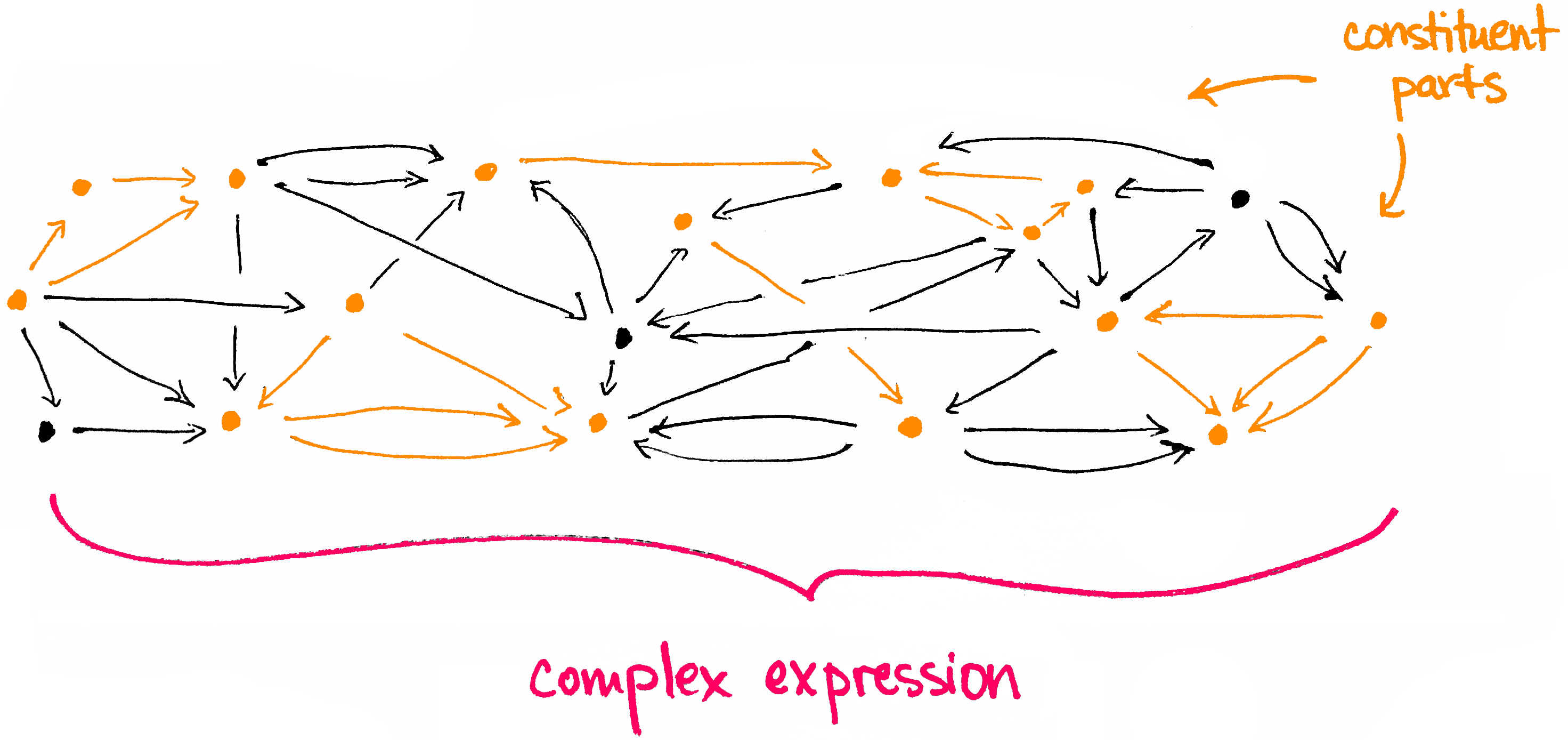}
\end{center}
\end{fullwidth}

In Section \ref{sec:Second}, 
	\marginnote{Matrix factorization provides another illustration of compositionality in mathematics. As an example, every $n\times m$ matrix $M$ has a \textit{singular value decomposition}, which means it can be written as a product of three matrices $M=UDV^\dagger$ where $U$ and $V$ are unitary square matrices (here $V^\dagger$ denotes the conjugate transpose of $V$) and $D$ is a rectangular diagonal matrix. Intuitively then, the linear transformation $M$ can be broken down into a rotation followed by a shear followed by another rotation. So you can analyze your transformation (or your data set, if that's what $M$ is encoding) by understanding its constituent pieces---the factors---and how they compose together. More generally, I like to think that tensor networks are a good example of compositionality, but such a discussion might take us too far off course. Perhaps another day!}
our example from natural language, the complex \textcolor{RubineRed}{expression} will be a sentence; its constituent \textcolor{BurntOrange}{parts} are the words that comprise the sentence.

%image source from Baez: http://math.ucr.edu/home/baez/networks/networks_16.html

In both examples, functorial semantics and the principle of compositionality will go hand-in-hand. The former prompts us to model behavior using a functor between syntax and semantics categories. The latter encourages us to take things one at a time: \textit{To model a huge system,} compositionality tells us, \textit{it's enough to model smaller pieces of it and then stick those pieces together.} Simple enough. But what does it mean to ``stick pieces together'' \textit{mathematically}? The answer is provided by the structure of a \textbf{monoidal category}. And \textit{that} is the first of our two main constructions in ACT.

\newpage
\subsection{Further Reading}
\noindent\textbf{For more on functorial semantics and compositionality:}
\begin{itemize}
	\item Take a look at (this small notice on) William Lawvere's 1963 PhD thesis \href{http://www.pnas.org/content/pnas/50/5/869.full.pdf}{``Functorial Semantics of Algebraic Theories''} for the formal foundations for functorial semantics.

	\item You might also enjoy this discussion on \textit{doctrines} over at the $n$-Category Caf{\'e}: \href{https://golem.ph.utexas.edu/category/2006/09/doctrines.html}{https://golem.ph.utexas.edu/category/2006/09/doctrines.html}

	\item The preface to Brendan Fong's PhD thesis, \href{https://arxiv.org/abs/1609.05382}{``The Algebra of Open and Interconnected Systems''}, has a nice discussion on the principal of compositionality and includes various references. And while you're at it, take a look at the entire thesis, which is wonderfully written and provides the backbone of much of John Baez's current research in \href{http://math.ucr.edu/home/baez/networks/}{network theory}, which we'll talk a little bit about in Section \ref{sec:First}.

	\item Applied category theorist Jules Hedges has also written a nice exposition on compositionality, appropriately entitled \href{https://julesh.com/2017/04/22/on-compositionality/}{``On Compositionality.''} In the article, you'll find a link to the Stanford Encyclopedia on Philosophy's \href{https://plato.stanford.edu/entries/compositionality/}{entry on compositionality}, which gives a thorough overview of the topic.
\end{itemize}

\newpage
%%%%%%%%%%%%% SECTION 3 %%%%%%%%%%%%%%
\section{Two Constructions}
Two constructions that appear over and over (and \textit{over and over and over}) in (some projects in) applied category theory are \textbf{monoidal categories} and \textbf{decorated cospans}. Let's talk about the first one first.

\subsection{Monoidal Categories}\label{sec:Monoidal}

Actually, before we talk about monoidal categories, let's talk about \textit{monoids}.
	\marginnote[-1cm]{Monoids and monoidal categories were the main focus in the digression on page \pageref{fig:monfun1}, but now I'll proceed as if they are new to the reader.}
Here are three examples of monoids: the \textit{integers} $\mathbb{Z}$, the \textit{rational numbers} $\mathbb{Q}$, and the set of all $n\times n$ \textit{matrices} with real-number entries $ M_n(\mathbb{R})$. Well, technically \textit{these} are monoids:
\begin{fullwidth}
\[(\mathbb{Z},+,0) \hspace{4.4cm} (\mathbb{Q},\cdot,1) \hspace{4.6cm} (M_n(\mathbb{R}),\cdot,\mathbb{1})\]
\includegraphics{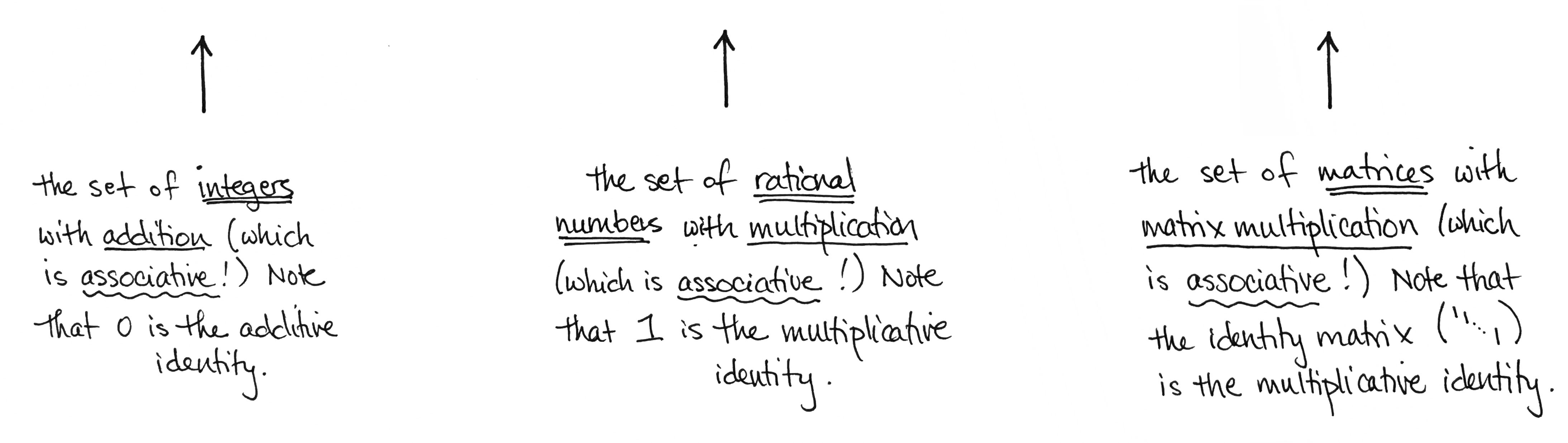}
\end{fullwidth}

Each example consists of a \textbf{set} $X$ equipped with an \textbf{associative binary operation}, which I'll denote by $\bullet$.
	\marginnote{A ``\textit{binary operation} on $X$'' is just the fancy name for a function $X\times X\to X.$ So $\bullet$ is a function $\bullet\colon X\times X\to X.$ It assigns to a pair $(x,y)\in X\times X$ an element $x\bullet y\in X$. To say that an element $1\in X$ serves as an \textit{identity} for $\bullet$ means it satisfies $1\bullet x=x\bullet 1=x$ for all $x\in X.$ } 
Moreover, there is a special element in the set, let's call it $1$, that serves as an \textbf{identity} for the  operation. \textit{Those three things}---a set, an associative binary operation, an identity---comprise a \textbf{monoid}. Usually, we write this triple as \[(X,\bullet,1)\]

\noindent Not too bad, right?

\vspace{0.4cm}

Great. Now imagine replacing the set $X$ by a \textit{category} $\cat{C}$, and replacing the binary operation $\bullet\colon X\times X\to X$ by a \textit{functor} $\bullet\colon  \cat{C}\times \cat{C}\to \cat{C}$, and replacing the identity element $1\in X$ by an \textit{object} $1$ in $\cat{C}$. The resulting triad
\[(\cat{C},\bullet,1)\] 
is called a \textbf{monoidal category}. The object $1$ is often called the \textit{monoidal unit}, and people usually prefer to write $\otimes$ (and call it the \textit{monoidal product}) instead of $\bullet$ so let's do that too:
\begin{center}
\includegraphics{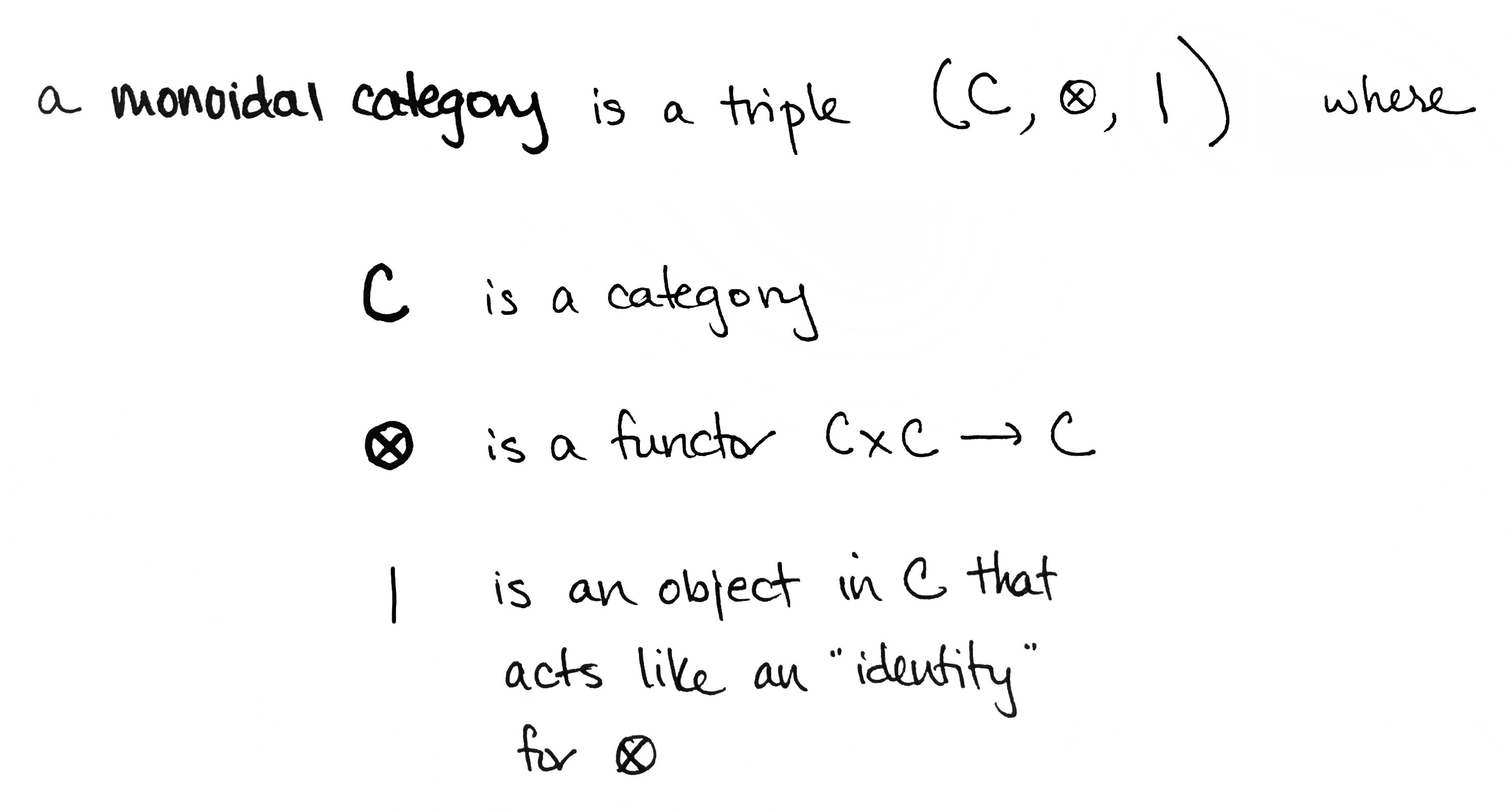}
\end{center}
	\marginnote[-5cm]{There's a little more to the story here since we want $\otimes$ to be associative: for any objects $A,B,C$ in the category $\cat{C}$, we want $A\otimes(B\otimes C)=(A\otimes B)\otimes C$. Alas, things are rarely equal on the nose. To compensate for this, we ask instead that there exist an \textit{isomorphism} $A\otimes(B\otimes C)\cong (A\otimes B)\otimes C$, which should behave nicely. I won't go into the details here, but of course you can find more on the Wikipedia page on \href{https://en.wikipedia.org/wiki/Monoidal_category}{monoidal categories}. For a delightful exposition on the richness of monoidal categories, I strongly recommend \href{http://math.ucr.edu/home/baez/rosetta.pdf}{``A Rosetta Stone''} by John Baez and Mike Stay. It is a \textit{gem}.}

In short, a monoidal category is a category in which it makes sense to ``combine'' objects and morphisms.\footnote[][10.6cm]{Allow me to explain the ``and morphisms'' part. First remember, $\otimes$ is a \textit{functor}! That means it's an assignment on objects \textit{and} on morphisms. Consider $(\cat{Set},\times,\{\ast\})$, for example, where $\times$ assigns to a pair of sets $(A,B)$ their Cartesian product $A\times B$. And given two functions $f\colon A\to B$ and $g\colon A'\to B'$, it assigns to the pair $(f,g)$ the function $f\times g\colon A\times B\to A'\times B'$, which is defined by: $(f\times g)(a,b):=(fa,gb)$. This is one example of the action of the monoidal product on morphisms. More generally, what $f\otimes g$ \textit{is} depends on the explicit definition of $\otimes$.} As we'll see in Sections \ref{sec:First} and \ref{sec:Second}, each of the four categories mentioned on page~\pageref{fig:4cats}---$\mathsf{Petri\;nets},$ $\mathsf{dynamical}$ $\mathsf{systems},$ $\mathsf{grammar},$ and $\mathsf{meanings\;of\;words},$---are monoidal categories! Here are some more examples.
\begin{fullwidth}
\[(\cat{Set},\times,\{\ast\}) \hspace{4.4cm} (\cat{Top},\sqcup,\varnothing) \hspace{4.6cm} (\cat{FVect},\otimes,\mathbbl{k})\]
\includegraphics{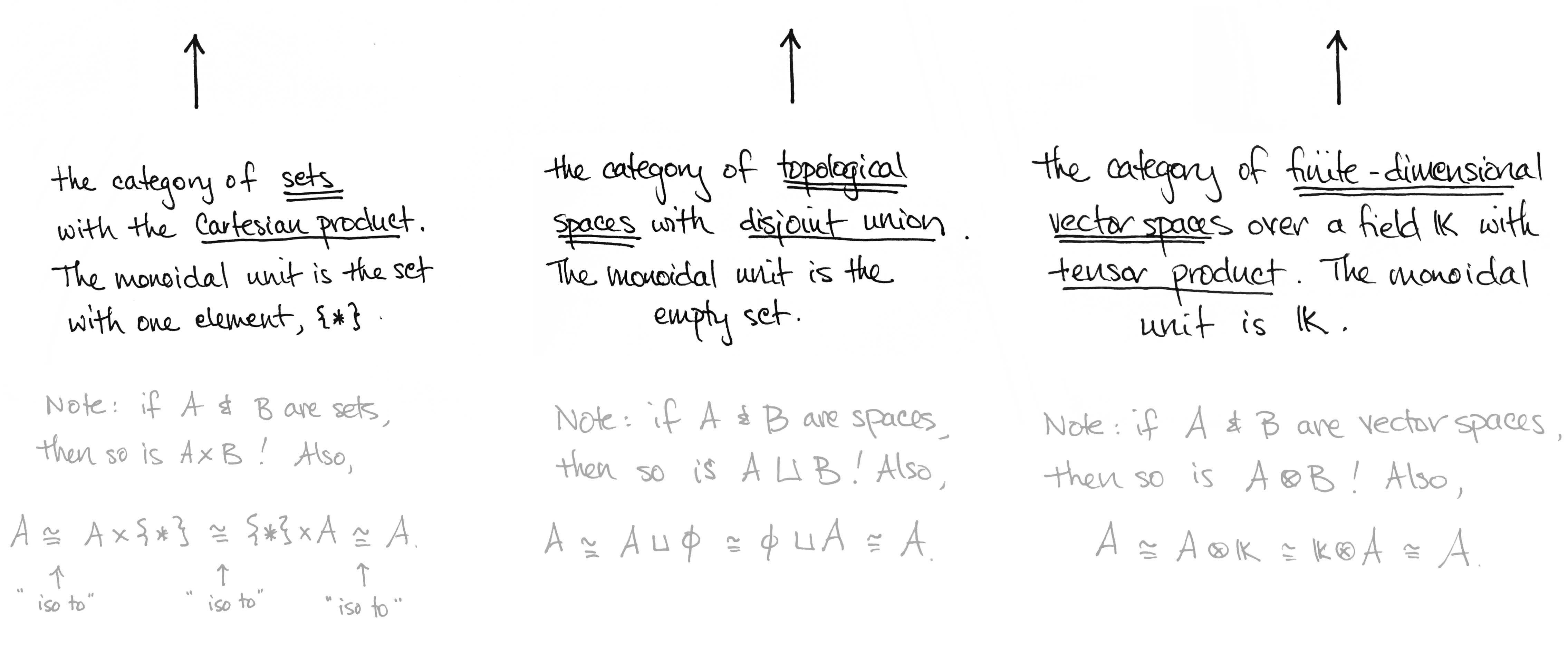}
\end{fullwidth}

By the way, if there is an isomorphism $A\otimes B\cong B\otimes A$ for all objects $A$ and $B$ that behaves nicely in a sense that can be made precise, then we say that $(\cat{C},\otimes,1)$ is a \textbf{symmetric monoidal category}. Each of the three examples above are symmetric monoidal. Monoidal categories come in other flavors too (braided, Cartesian, closed, Cartesian closed, closed braided,...), depending on which properties are satisfied.

\newthought{The main takeaway} here is that monoidal categories are the bread and butter of many applied category theorists. One reason for this is that monoidal categories provide a good setting in which to view morphisms ${\color{red}\to}$ as \textit{physical processes} and objects $A,B,\ldots$ as \textit{states}. As a non-technical example, let's suppose $A$ is a bunch of lemon meringue pie ingredients while $B$ is a fully-assembled-yet-unbaked lemon meringue pie. We might view a morphism $A\:{\color{red}\to}\:B$ as the process of mixing the raw ingredients together and then pouring the resulting concoction into a pre-baked crust.
\begin{center}
\refstepcounter{dummy}
\label{fig:pie}
\includegraphics[width=!,totalheight=!,scale=0.4]{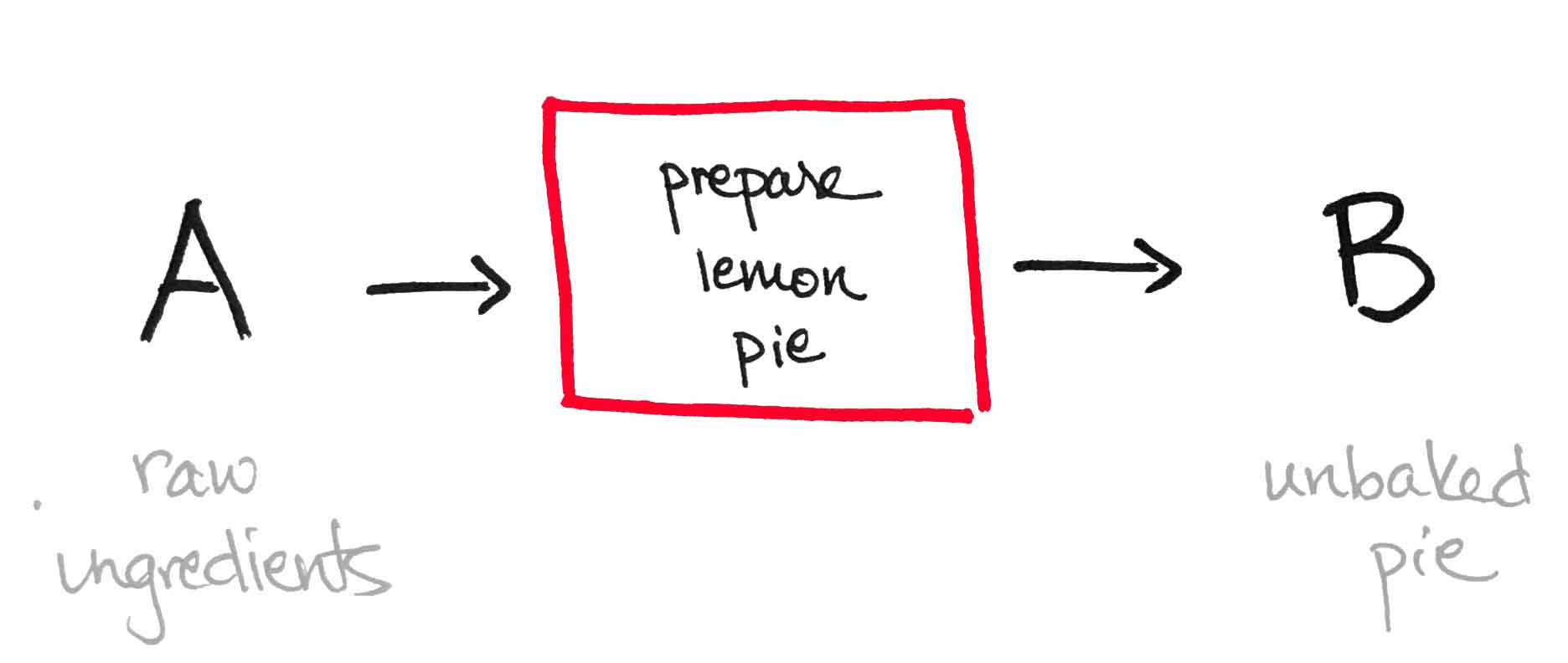}
\end{center}
	\marginnote[-4cm]{\includegraphics[width=!,totalheight=!,scale=0.6]{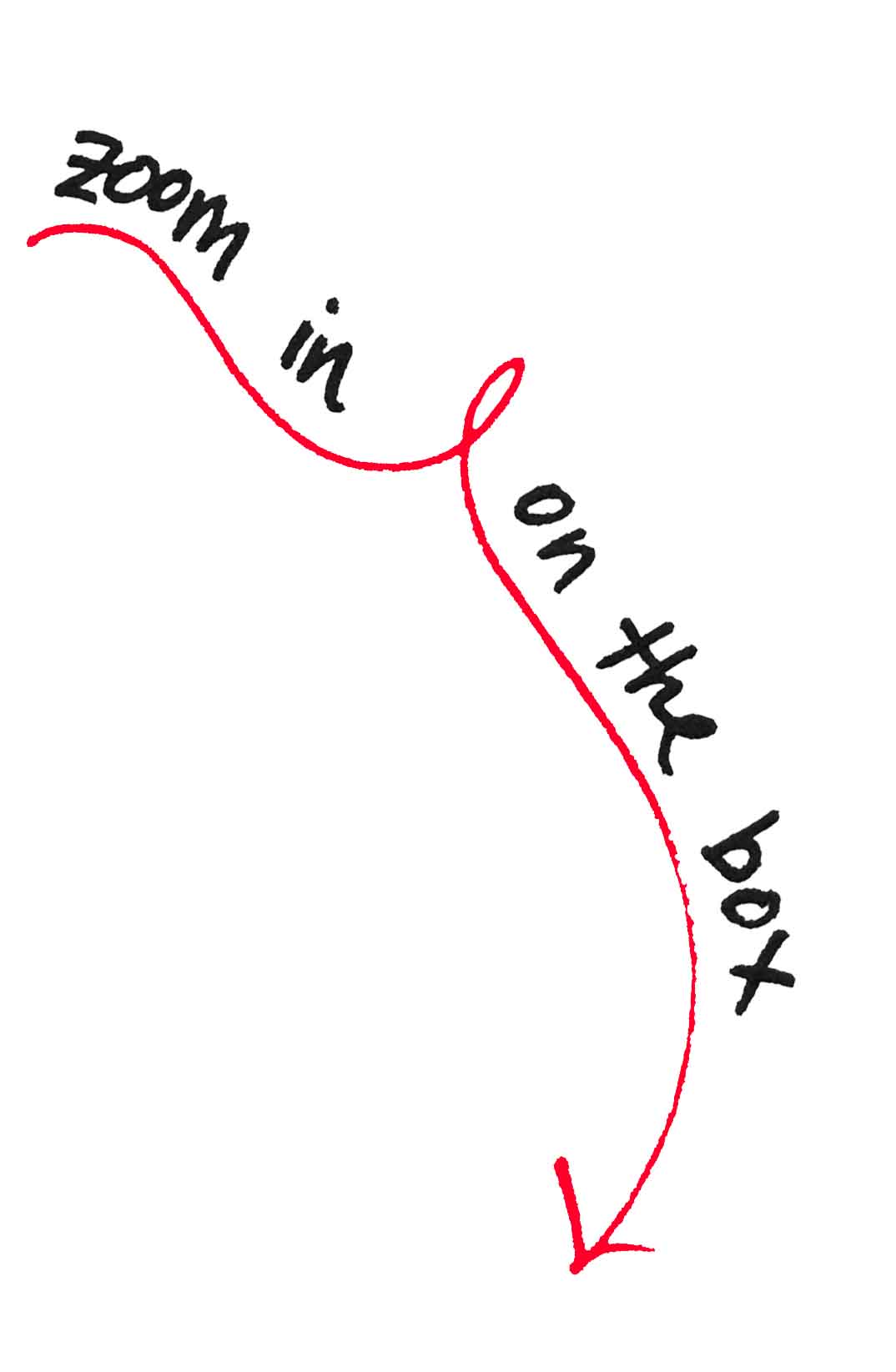}}
As it turns out, this pie example isn't so silly after all. It's one of the motivating examples that Brendan Fong and David Spivak use in their \textbf{excellent} book \href{https://arxiv.org/pdf/1803.05316.pdf}{\textit{Seven Sketches in Compositionality: An Invitation to Applied Category Theory}} to illustrate both the ubiquity and the simplicity of monoidal categories. (If you haven't read \textit{Seven Sketches} yet, you really \textbf{must}.) Below is a copy of their lemon meringue pie diagram, where I've drawn our $A$ and $B$ on the left as input and right as output.
\begin{fullwidth}
\begin{center}
\includegraphics[width=!,totalheight=!,scale=0.5]{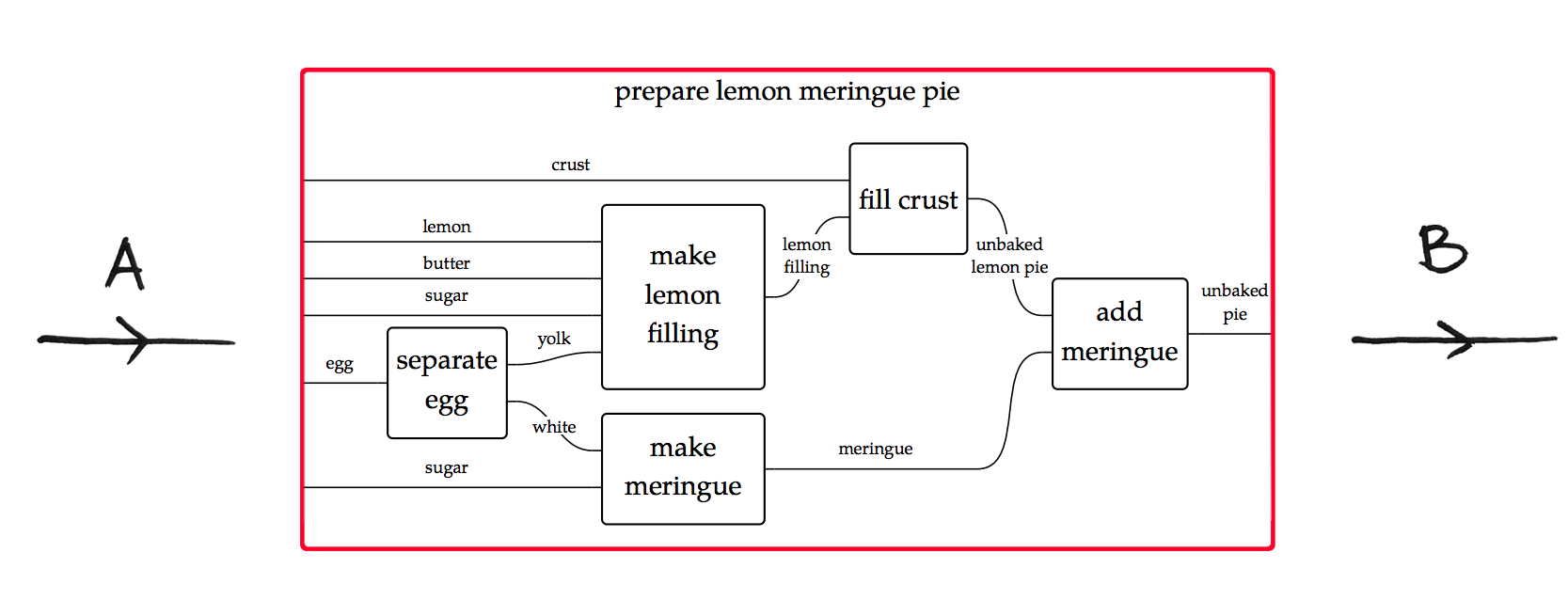}
\end{center}
\end{fullwidth}
Now that we've zoomed in, we can see that our process {\includegraphics[width=!,totalheight=!,scale=0.11]{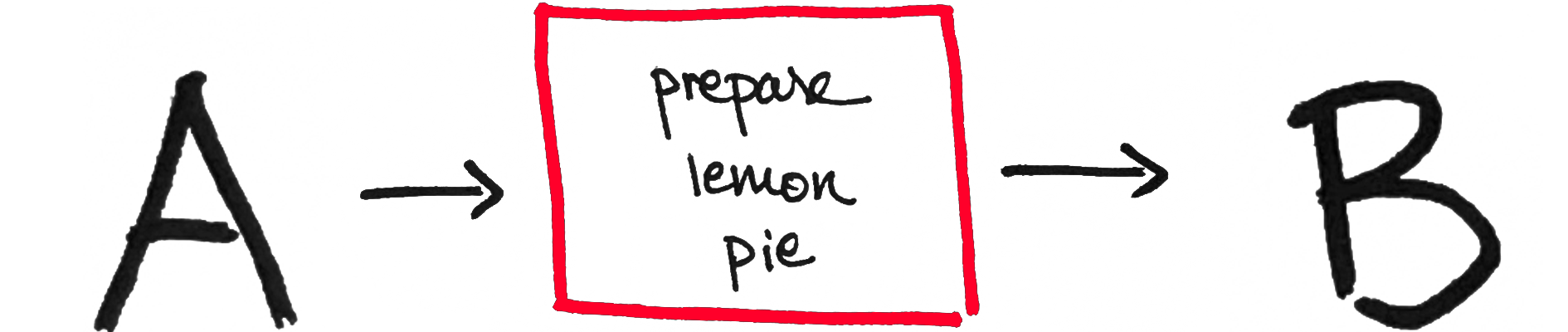}} is actually made up of a bunch of other processes! This isn't too surprising as there are several steps that go into preparing a lemon pie: \textit{separating the eggs, making the lemon filling, filling the crust}, and so on. Fong and Spivak's diagram illustrates just how those those individual steps combine to form the single process \textit{prepare lemon meringue pie}. What's neat is that we can describe these steps using the language of monoidal categories! We'll go into more detail later in this section, but here's a quick preview:

	\begin{itemize}
	\item The category's composition $\circ$ corresponds to using one box's output wire as another box's input wire. For example,
	\begin{center}
	\includegraphics[width=!,totalheight=!,scale=0.5]{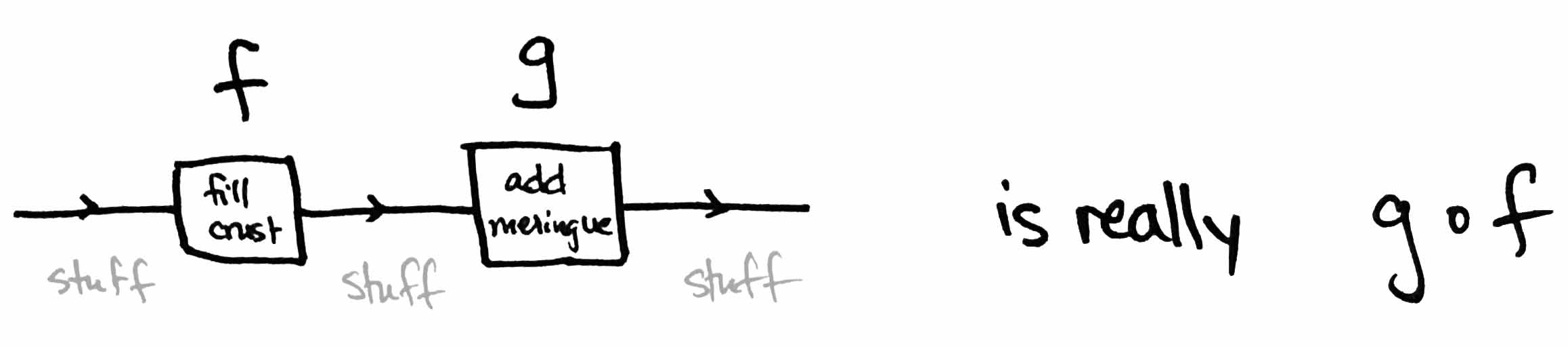}
	\end{center}
	\item The monoidal product $\otimes$ corresponds to stacking boxes on top of each other. For example,
	\begin{center}
	\includegraphics[width=!,totalheight=!,scale=0.5]{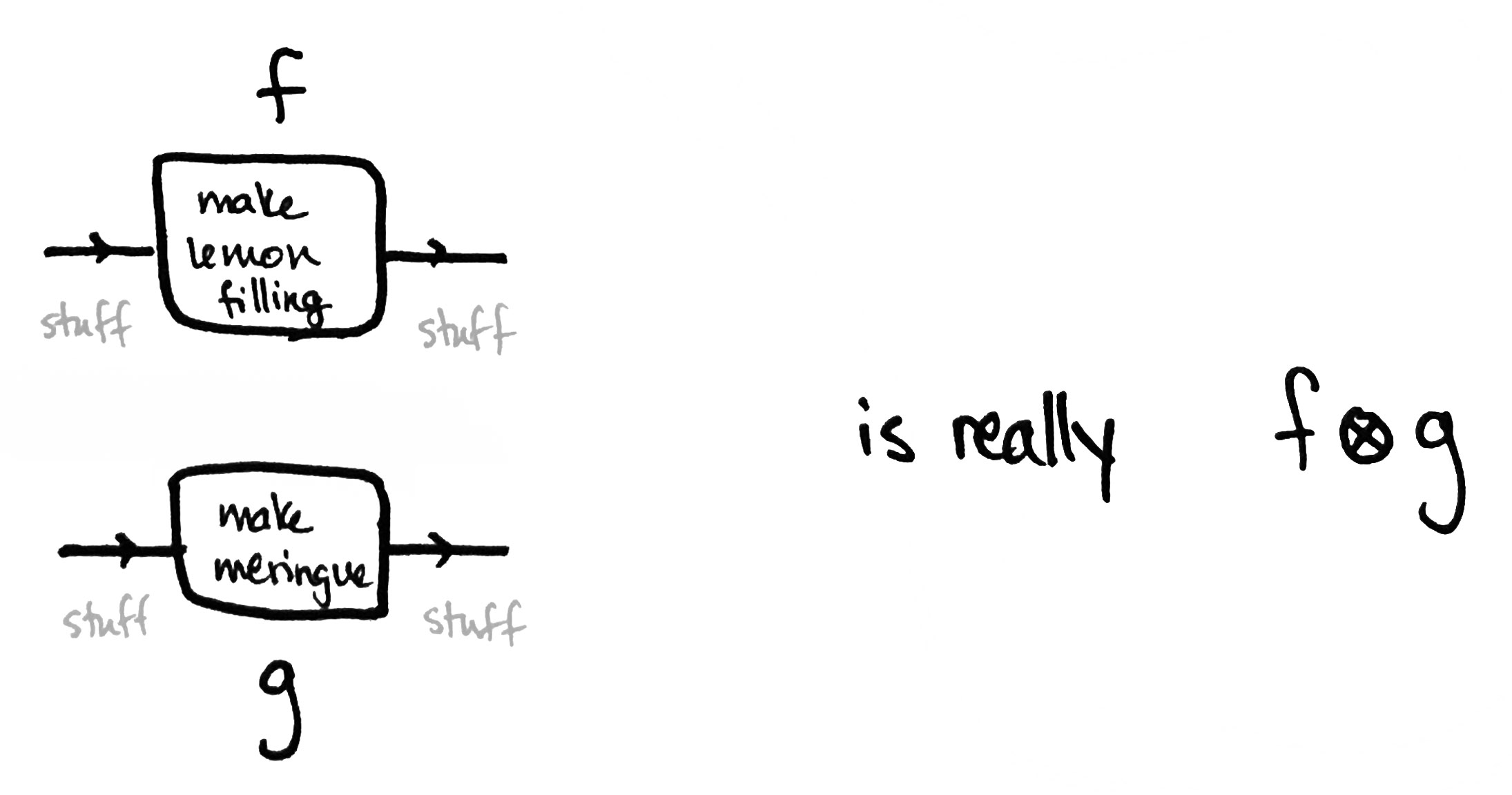}
	\end{center}
	\end{itemize}
In other words, $\circ$ means ``do the processes in series'' while the monoidal product $\otimes$ means ``do the processes in parallel.'' The resulting picture is called a \textbf{string diagram}---a graphical representation of a process (or equation of processes) in a monoidal category. I'll give more detail on how string diagrams work in a second.  But first, I'm reminded of something else about monoidal categories that I want to tell you! So let me tell you this new bit of information, then we'll come back to string diagrams. This small digression will, in fact, tie things together quite nicely. Bear with me.

\newthought{Earlier, I mentioned} that the word ``symmetric'' can be used as an adjective for ``monoidal categories'':
\begin{center}
\includegraphics[width=!,totalheight=!,scale=0.5]{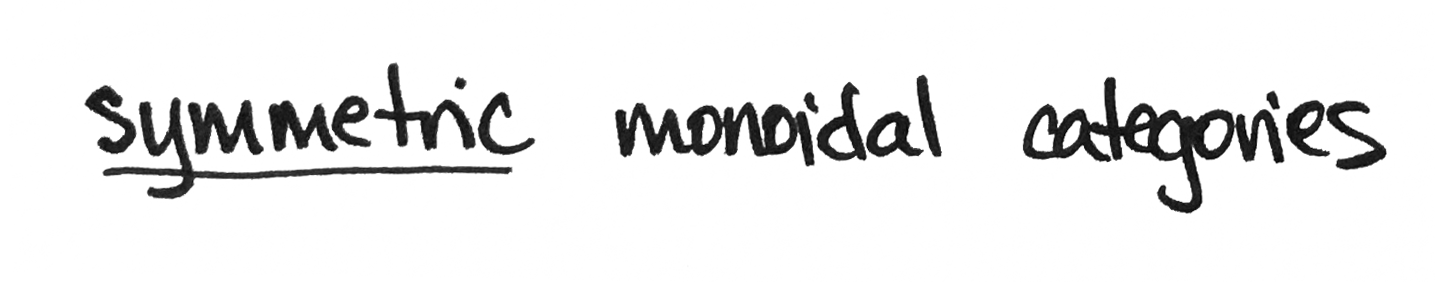}
\end{center}
Similarly, there is another flavor of monoidal categories that we should know about. This one will provide the main setting for our example in Section \ref{sec:Second}:
\begin{center}
\includegraphics[width=!,totalheight=!,scale=0.5]{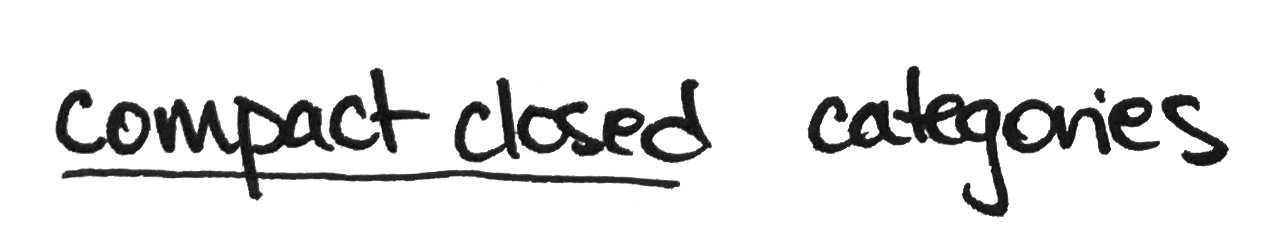}
\end{center}
I'll explain.

\subsection*{Interlude: Compact Closed Categories and String Diagrams}
Do you know what a finite dimensional vector space over $\mathbb{R}$ is? Then you know what a compact closed category is! Or rather, you know an example of one: 
\[\cat{FVect} \textit{ is a compact closed category}.\]
What makes that sentence true? Answer: every finite dimensional vector space $V$ has a \textit{dual} space $V^*=\hom(V,\mathbb{R})$.

\vspace{0.5cm}
\noindent \textit{That's it.} 
\vspace{0.5cm}

\noindent A \textbf{compact closed category} just the name for a monoidal category in which every object has a dual!\footnote{Technically, every object must have a \textit{left} dual and a \textit{right} dual. We need the distinction because not all monoidal categories are \textit{symmetric} monoidal.} But what does ``has a dual'' mean? In other words, what makes a dual \textit{dual}? Before I \sout{tell you the  answer} cite the definition, let's think back to the category $\cat{FVect}$, by way of motivation, and let's assume each vector space comes with an inner product. In this case, there are two very important linear maps between the ground field---let's say it's $\mathbb{R}$ for now---and a vector space $V$ tensored with its dual:
\[
\eta_V\colon \mathbb{R}\to V\otimes V^* \qquad \qquad \epsilon_V\colon V^*\otimes V\to\mathbb{R}
\]
In fact, there's a nice fact from linear algebra, namely that once we fix a basis $\{\mathbf{e}_1,\ldots,\mathbf{e}_n\}$ for $V$ then there is an isomorphism $V\cong V^*$. So let's fix that basis (the standard one) and write $V$ instead of $V^*$. Also, the subscript is a little cumbersome, so let's drop it for now.
		\marginnote[-1cm]{We'll need $\eta$ and $\epsilon$ for a computation in Section \ref{sec:Second}, so it's good to see what they look like explicitly.}
So we have two maps
\begin{align}\label{lis:ucounit}
\eta\colon \mathbb{R}\to V\otimes V \qquad \qquad \epsilon\colon V\otimes V\to\mathbb{R}
\end{align}
The map $\eta$ is called the \textbf{unit} \footnote{Note: this unit is not to be confused with ``monoidal unit''!}, and it assigns to every real number a vector in $V\otimes V$, namely:
\[\eta(1)=\sum_{i=1}^n \mathbf{e}_i\otimes \mathbf{e}_i \qquad \qquad\text{(and extend linearly)}\]
The map $\epsilon$ is called the \textbf{counit}\footnote{Note: \textit{counit} is \href{https://twitter.com/math3ma/status/965271641847685120}{pronounced} ``coh-yew-nit'' not ``cow-nit.'' This is \href{https://twitter.com/math3ma/status/965255364286574594}{important}.}, and it assigns to every vector in $V\otimes V$ a real number, namely:
\[\epsilon\left(\sum_{i,j}c_{ij}\mathbf{v}_i\otimes \mathbf{w}_j\right)=\sum_{i,j}c_{ij}(\mathbf{v}_i\cdot\mathbf{w}_j)\qquad\text{where $\cdot$ is the inner product}\]
	Intuitively, we can think of $\epsilon$ as an \textit{evaluation} map. That's because there is always a map $V^*\otimes V\to\mathbb{R}$ given by evaluation. Indeed, if $\mathbf{v}\in V$ and $f\in V^*$, then we can pair the two together to obtain $f\mathbf{v}\in \mathbb{R}.$ And if we view $f$ as a $1\times n$ matrix and $\mathbf{v}$ as an $n\times 1$ matrix, then $f\mathbf{v}$ is their inner product:
	\[
	\begin{bmatrix}\cdots f \cdots
	\end{bmatrix}
	\begin{bmatrix}\vdots\\\mathbf{v}\\\vdots
	\end{bmatrix} = \text{ a number}
	\]
	The $\epsilon$ map just extends this linearly. That is, if we now have \textit{any} vector $\sum_{i,j}c_{ij}\mathbf{v}_i\otimes \mathbf{w}_j$ in $V\otimes V,$ then  $\epsilon\colon V\otimes V\to\mathbb{R}$ is given by
	\[\epsilon\left( \sum_{i,j}c_{ij} \mathbf{v}_i\otimes \mathbf{w}_i\right)=\sum_{i,j}c_{ij}(\mathbf{v}_i\cdot \mathbf{w}_j)\]
	as above.
	% \marginnote[-2.5cm]{\textcolor{gray}{By the way, this is exactly what you get if you wanted to just define $\epsilon$ on basis vectors $\{e_i\otimes e_j\}$, assuming the $\{e_i\}$ form an orthonormal basis. For in that case, $\epsilon$ assigns to $\sum_{ij}c_{ij}e_i\otimes e_j$ the value $c_{11}+c_{22}+\cdots+c_{nn}.$}}

Finally, the unit $\eta$ and counit $\epsilon$ interact nicely with each other because they satisfy some equations called the \textit{yanking equations,} which I'll explain shortly. The bottom line is that all the above---the maps $\eta$ and $\epsilon$ and the equations they satisfy---makes $V^*$ into a bona fide \textit{dual} for $V$. The upshot is that compact closed categories generalize these notions.

\begin{defn}A \textbf{compact closed category} is a monoidal category $(\cat{C},\otimes,1)$ where for every object $c$ in $\cat{C}$ there exists objects $c^l$ and $c^r$ and morphisms	
	\marginnote{The $\eta$ maps are called the left and right \textbf{units}, and the $\epsilon$ maps are called the left and right \textbf{counits}.}
\begin{align*}
\eta_c^l\colon 1&\to c\otimes c^l\qquad\qquad \epsilon_c^l\colon c^l\otimes c\to 1\\
\eta_c^r\colon 1&\to c^r\otimes c\qquad\qquad \epsilon_c^r\colon c\otimes c^r\to 1
\end{align*}
that satisfy the ``yanking (or snake) equations''
	\marginnote{Here $\id_c$ denotes the identity morphism $\id_c\colon c\to c.$}
\begin{equation}\label{eq:yank}
\begin{aligned}
&(\id_c \otimes \epsilon^l )\circ (\eta^l \otimes \id_c) = \id_c &\:& (\epsilon^r \otimes \id_c) \circ (\id_c \otimes \eta^r ) = \id_c
 \\
&(\epsilon^l \otimes \id_{c^l})\circ (\id_{c^l} \otimes \eta^l)= \id_{c^l} &\:& (\id_{c^r} \otimes \epsilon^r ) \circ (\eta^r \otimes \id_{c^r} ) = \id_{c^r}
\end{aligned}
\end{equation}
\end{defn}

\vspace{0.5cm}
\noindent \textit{Yikes. What do these equations MEAN?}
\vspace{0.5cm}
 
\noindent I'm glad you asked. 
\vspace{0.5cm}

To answer, it's time to revisit our previous discussion on \textbf{string diagrams!} String diagrams are loved by applied category theorists far and wide because \textit{they make life SO much easier.}
	\marginnote[-0.5cm]{$\leftarrow$ In that sentence, I meant to convey that \textit{string diagrams} make life so much easier, but one may argue that \textit{applied category theorists} also (are working to) make life so much easier.}
As we saw earlier, a string diagram is  a picture that represents morphisms in a monoidal category $\cat{C}$. For now let's take that category to be $\cat{FVect}$ so that our objects are vector spaces $V,W,\ldots$. In this case, the left and right dual of space $V$ is its vector space dual
\[V^*= V^r=V^l\]
We'll get to the yanking equations shortly, but first: If this document has been your first introduction into string diagrams, then here is THE KEY thing to know:
	\begin{quote}
	In category theory, we often draw an object as a dot $\bullet$ and a morphism as an arrow $\bullet\to\circ$. To draw a string diagram, just \textbf{do the opposite!} (This goes back to \href{https://en.wikipedia.org/wiki/Poincar%C3%A9_duality}{Poincar{\'e} duality} in topology.) To draw a string diagram, draw an object as an arrow and a morphism as a dot or, even better, a box. 
		\marginnote[-1cm]{The lemon pie diagram that we saw on page \pageref{fig:pie} is an example of a string diagram! \begin{center}\includegraphics[width=!,totalheight=!,scale=0.17]{pie.png}\end{center}}
	\begin{center}
	\includegraphics[width=!,totalheight=!,scale=0.6]{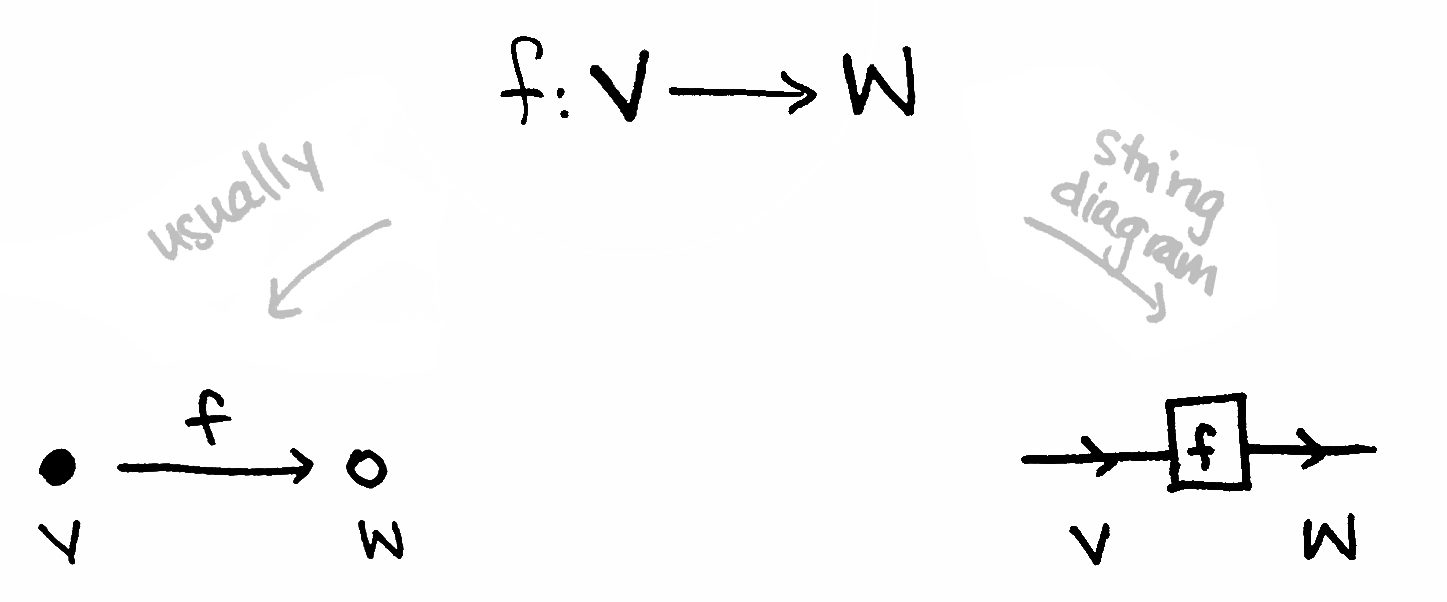}
	\end{center}
	\end{quote}
With this small artistic adjustment, we can represent the monoidal product $\otimes$ pictorially as well. The product of two spaces $V\otimes W$ is drawn as two lines, side-by-side. A similar picture holds for the product of two morphisms. Composition $\circ$ is represented by gluing strings together. 
\begin{center}
\includegraphics[width=!,totalheight=!,scale=0.5]{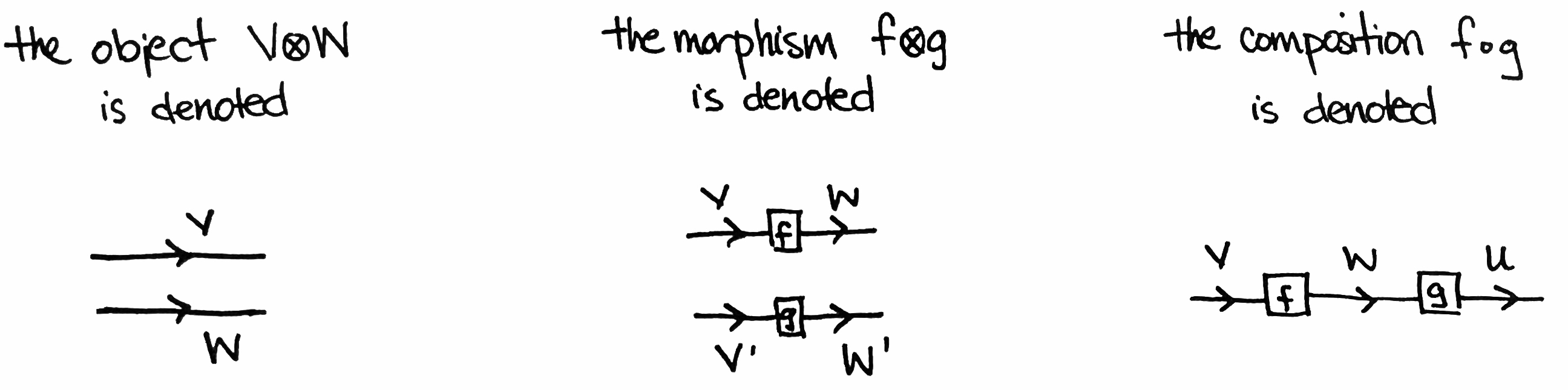}
\end{center}
And as we saw above, every object in a compact closed category such as $\cat{FVect}$ has a dual. Its picture is given by an arrow pointing in the \textit{opposite} direction. 
\begin{center}
\includegraphics[width=!,totalheight=!,scale=0.5]{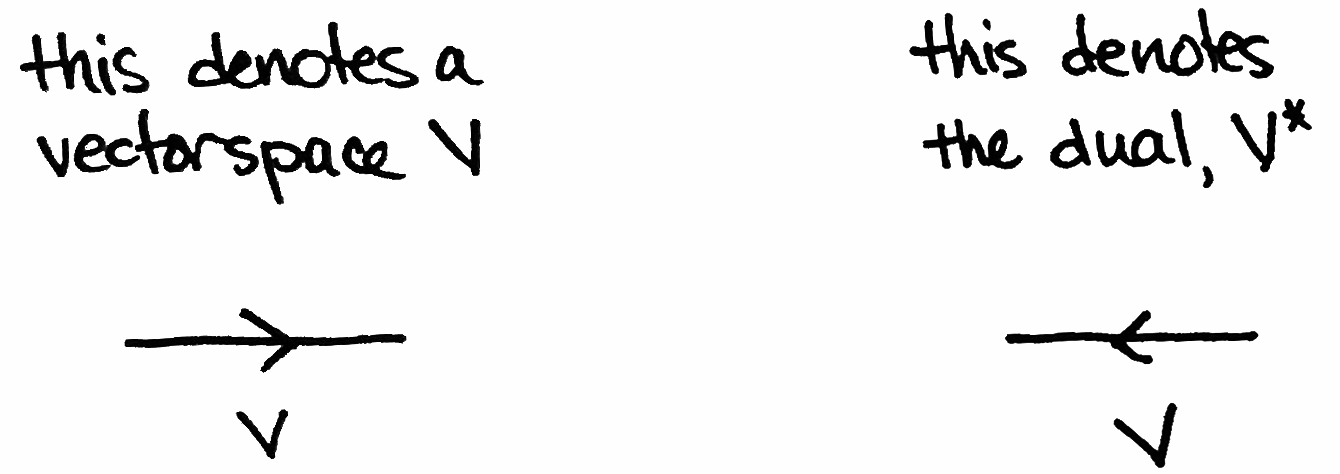}
\end{center}
% Any map $f\colon V\to W$ \textit{also} has a dual (or \textit{adjoint}) $f^*\colon W^*\to V^*.$ In the land of string diagrams $f^*$ is drawn by simply reflecting $f$.
% \begin{center}
% \includegraphics[width=!,totalheight=!,scale=0.5]{dualmaps.jpg}
% \end{center}
% Notice that even though the reflection causes the arrows to point right-to-left, we still read it from left-to-right. That is, {\includegraphics[width=!,totalheight=!,scale=0.35]{dualmaps3.png}} is a map \textit{from} $W^*$ \textit{to} $V^*$. In fact, the folks in \href{https://en.wikipedia.org/wiki/Categorical_quantum_mechanics}{categorical quantum mechanics} like to denote $f$ as a trapezoid like this, to emphasize the relationship.
% \begin{center}
% \includegraphics[width=!,totalheight=!,scale=0.5]{dualmaps2.jpg}
% \end{center}
Another special object in a compact closed category is the monoidal unit, for instance $\mathbb{R}$ in $\cat{FVect}$. Because the unit is an \textit{object}, it's depicted as an arrow, too. People like to draw this arrow in a special way, namely as the ``empty'' arrow. In other words, people don't like to draw an arrow. That's because the monoidal unit $\mathbb{R}$ satisfies\footnote{More generally, the monoidal unit $1$ in a monoidal category $(\mathsf{C},\otimes,1)$ satisfies $1\otimes c\cong c \cong c\otimes 1$ for all objects $c$ in $\mathsf{C}$.}
\[V\otimes \mathbb{R}\cong V\cong \mathbb{R}\otimes V \qquad\text{for all $V$},\]
which suggests that the unit is ``invisible.'' But I like to draw it anyway, shaded:
\begin{center}
\includegraphics[width=!,totalheight=!,scale=0.5]{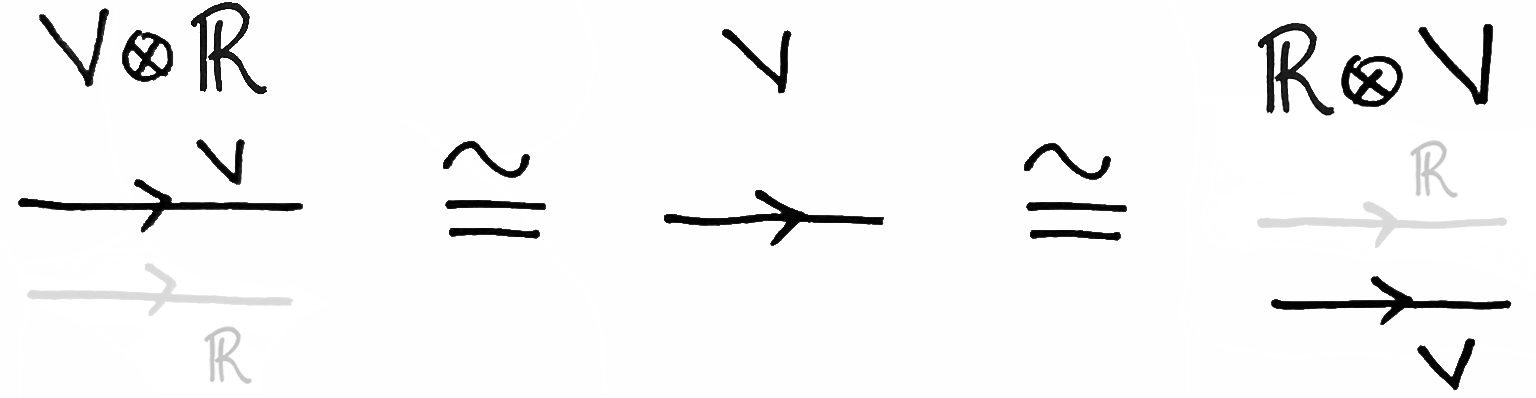}
\end{center}

\newthought{Now we are ready to} get back to the yanking equations. Remember, part of the data of a compact closed category is that each object $V$ has left and right duals $V^r, V^l$ \textit{together with} morphisms
\begin{align*}
\eta_V^l\colon \mathbb{R}&\to V\otimes V^l\qquad\qquad \epsilon_V^l\colon V^l\otimes V\to \mathbb{R}\\
\eta_V^r\colon \mathbb{R}&\to V^r\otimes V\qquad\qquad \epsilon_V^r\colon V\otimes V^r\to \mathbb{R}
\end{align*}
Again, to simplify the notation we'll use the fact that that for vector spaces, $V^*=V^r=V^l.$ I'll also drop the subscripts to keep things clean.

Graphically, the $\eta$s and $\epsilon$s are drawn as below.
	\marginnote[3cm]{Alternatively, some folks will rotate the $\epsilon$ and $\eta$ diagrams by $90^\circ$ clockwise and counterclockwise, respectively, which is the reason for their common nickname of ``cups and caps.''}
The reason we have \textit{two} versions of each map is because the ``information flow'' can either flow up or it can flow down.\marginnote[4cm]{Note: the direction of the (invisible) arrow for the unit $\mathbb{R}$ can go either way. The monoidal unit is always self dual!}
\begin{center}
\includegraphics{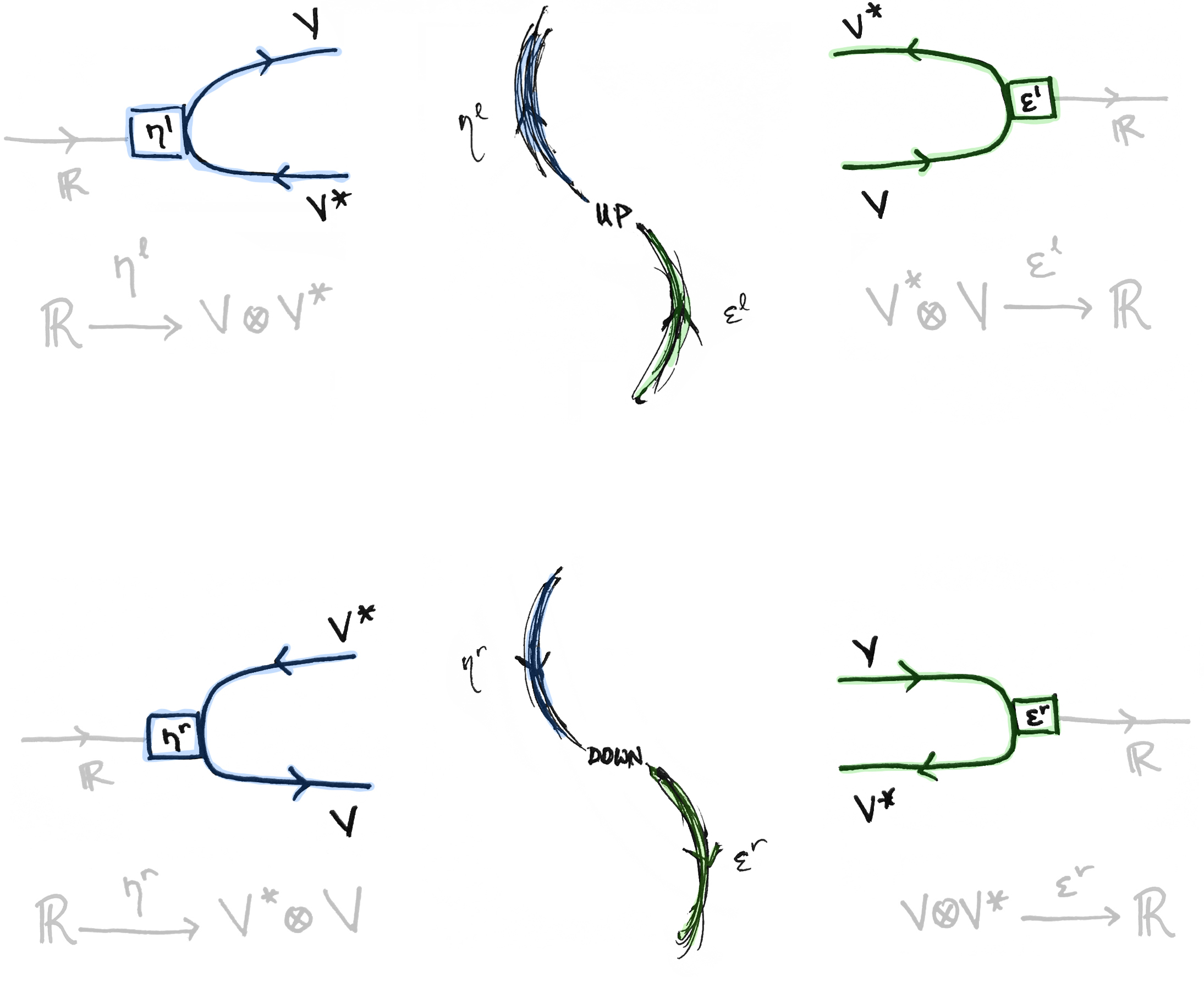}
\end{center}
And since $\mathsf{FVect}$ is a \textit{symmetric} monoidal category, and since the left and right duals are both $V^*$, there is really only \textit{one} unit and \textit{one} counit for vector spaces.\footnote[][-2cm]{Remember, the sentence ``$\mathsf{FVect}$ is symmetric monoidal'' means there is an isomorphism $V\otimes W\cong W\otimes V$ for every pair of vector spaces $V$ and $W.$ In string diagram calculus, this means that the order in which we draw our arrows doesn't matter:
	\begin{center}
	\includegraphics{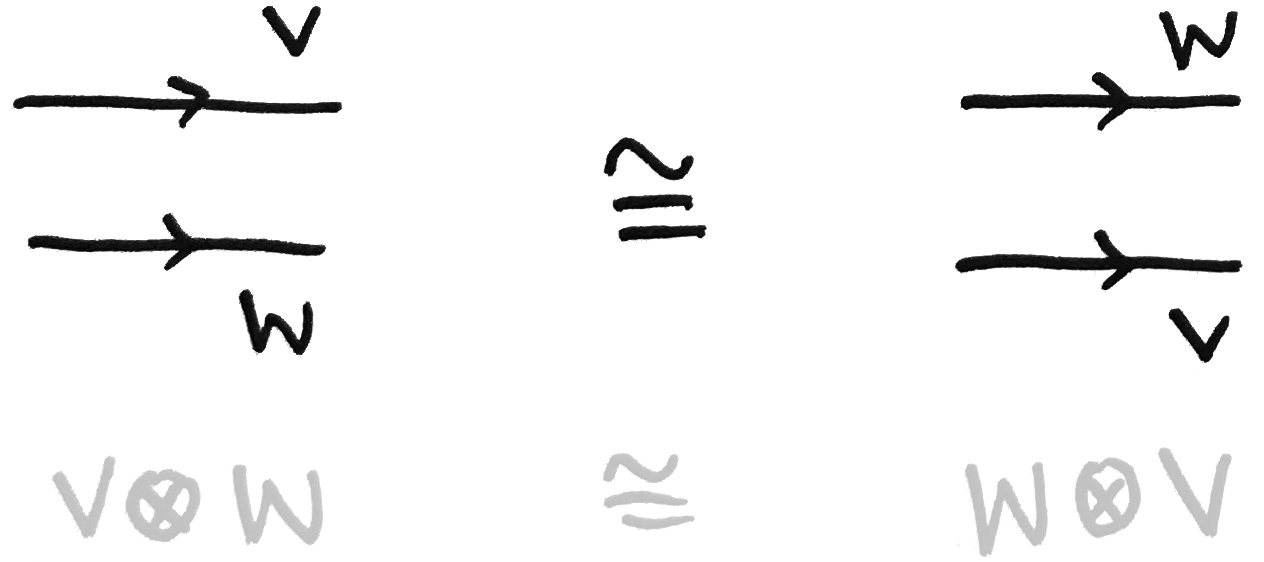}
	\end{center}}
\begin{center}
\includegraphics[width=!,totalheight=!,scale=0.35]{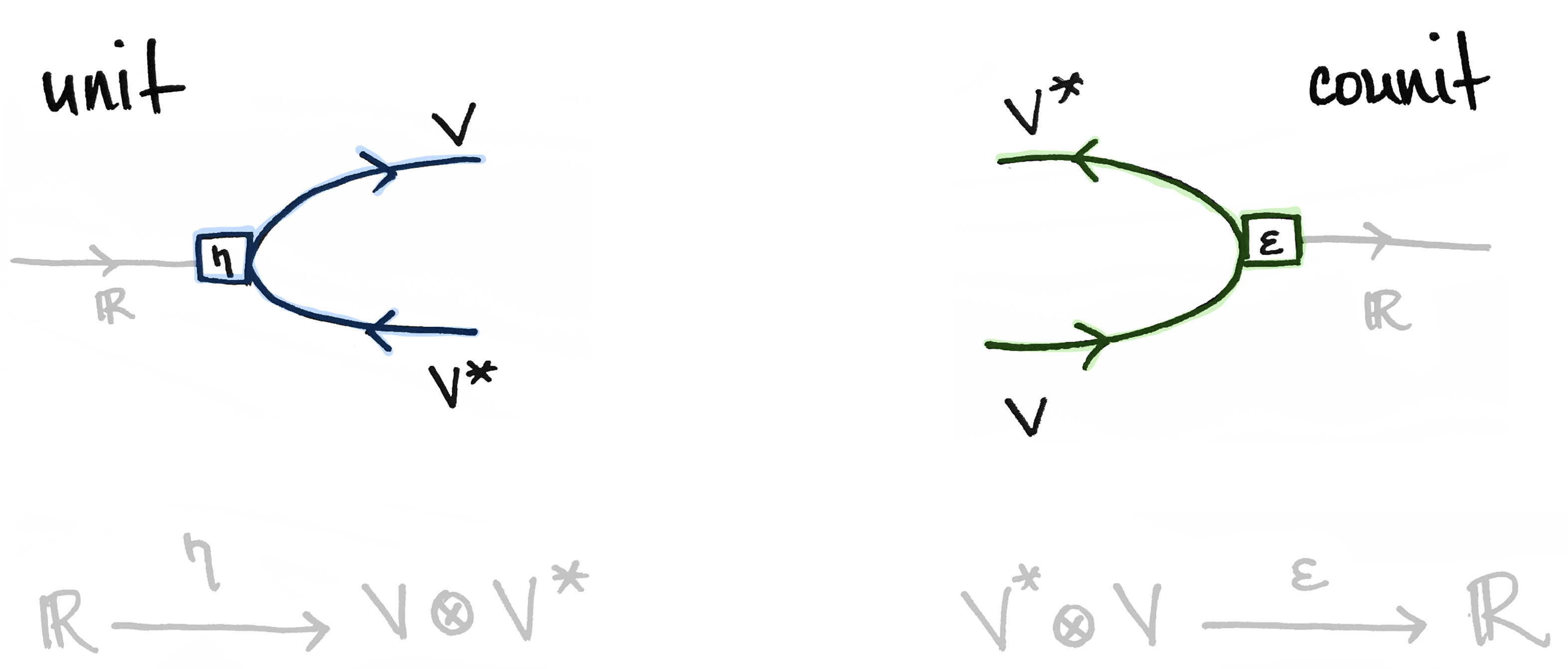}
\end{center}
That is, in $\mathsf{FVect}$ \[\eta=\eta^r=\eta^l \qquad\text{and}\qquad\epsilon=\epsilon^r=\epsilon^l,\]
and these are \textit{precisely} the $\eta$ and $\epsilon$ defined on page \pageref{lis:ucounit}! So in this example, the four yanking equations of (\ref{eq:yank}) reduce down to just \textit{two}:
\begin{equation}\label{eq:yankvect}
\begin{aligned}
(\epsilon \otimes \id_V)\circ (\id_V\otimes \eta)&=\id_V\\[7pt]
(\id_{V^*}\otimes\:\epsilon)\circ (\eta\otimes \id_{V^*})&=\id_{V^*}
\end{aligned}
\end{equation}
Graphically, these equations can be represented as follows:
	\marginnote[-3cm]{For fun, verify that the unit and counit maps on page \pageref{lis:ucounit} do indeed satisfy these two equations.}
\begin{fullwidth}
\begin{center}
\refstepcounter{dummy}
\label{fig:yank}
\includegraphics{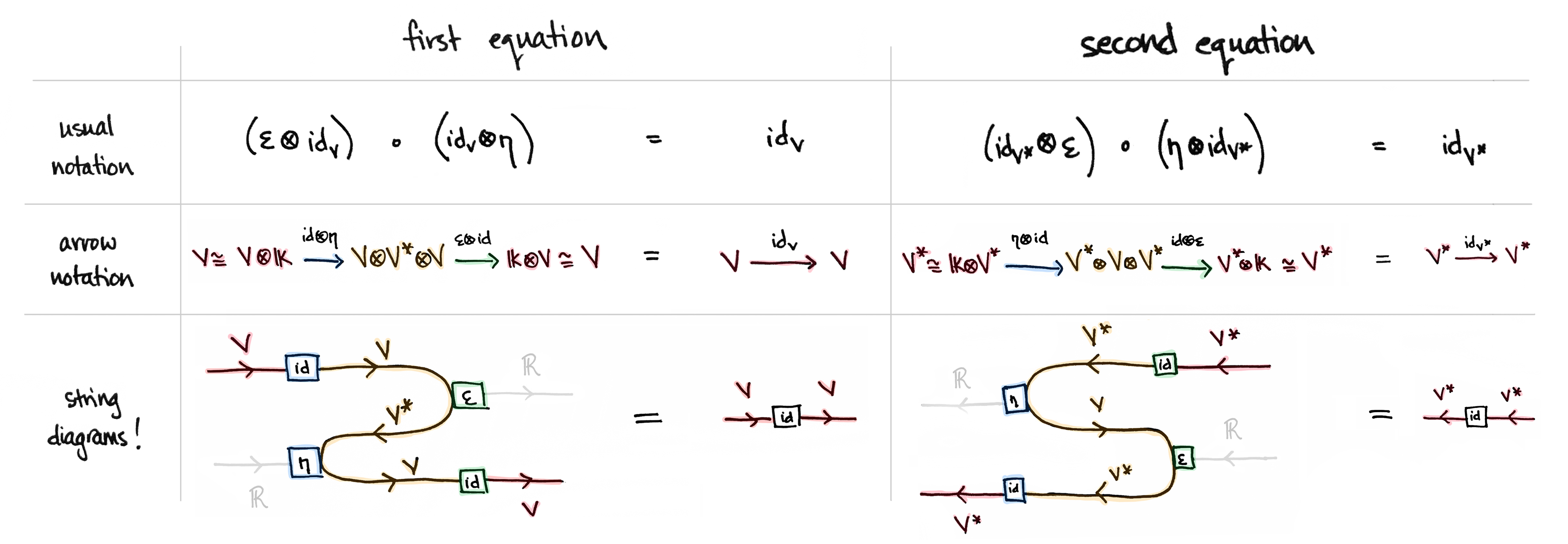}
\end{center}
\end{fullwidth}
\noindent After yanking the strings taut, you'll notice that information flows \textit{rightwards} in the first equation, while it flows \textit{leftwards} in the second equation.
%Notice that in both string diagrams, the overall information flow is downwards.
\begin{center}
\includegraphics[width=!,totalheight=!,scale=0.4]{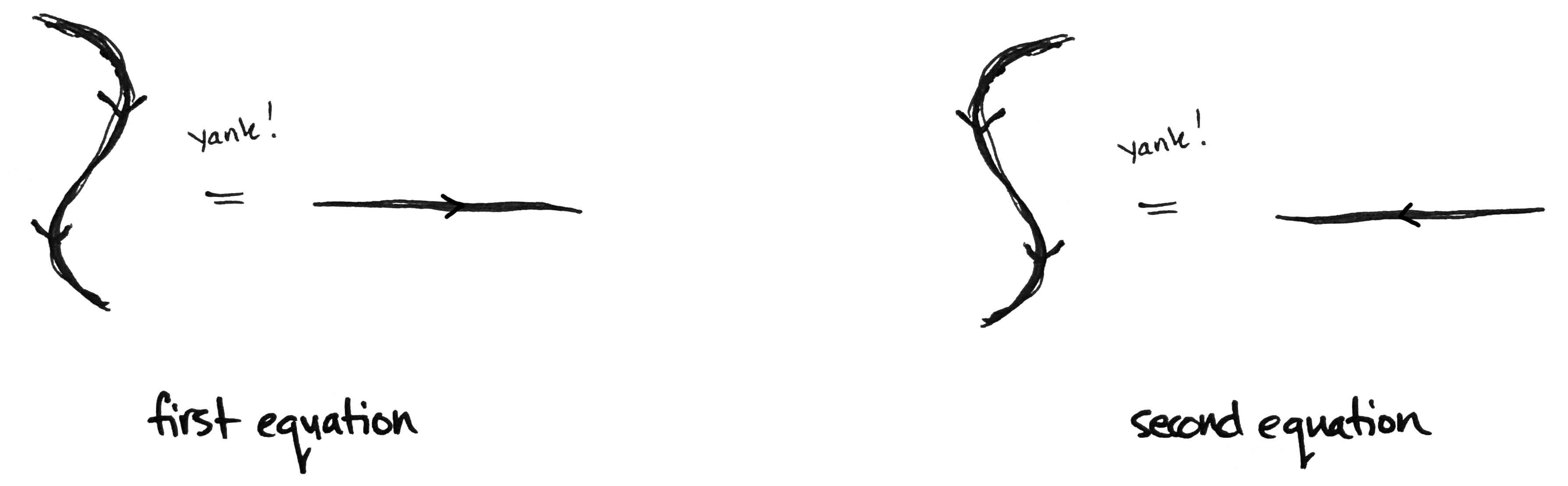}
\end{center}

\newpage
\noindent Since we're in a \textit{symmetric} monoidal category, nothing changes if we reverse the arrows in the pre-yanked strings. If, however, the monoidal product $\otimes$ is not symmetric, then we obtain two more diagrams.
\begin{center}
\includegraphics[width=!,totalheight=!,scale=0.4]{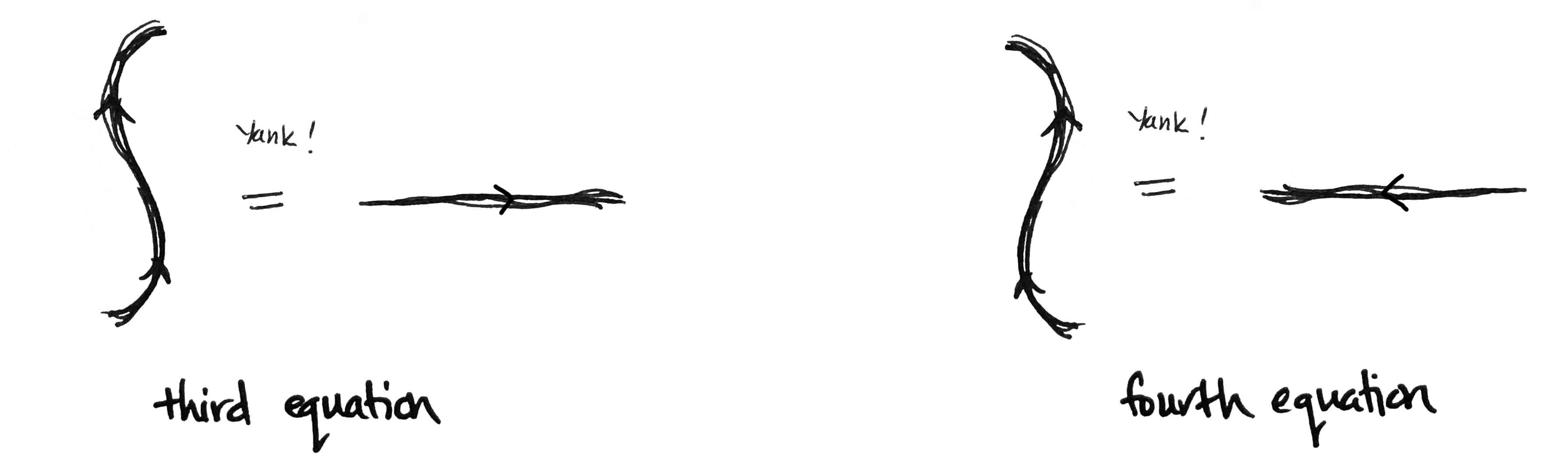}
\end{center}

This gives a grand total of four equations---the four displayed in (\ref{eq:yank}). And that's what gives us a compact closed category. The category $\mathsf{FVect}$ together with the unit and counit defined in (\ref{lis:ucounit}) will make another appearance in Section \ref{sec:Second}. In that same section, we'll also seen an example of a compact closed category that is \textit{not} symmetric.

% \marginnote{\textcolor{Green}{Tai, can you see the adjunction triangle identities in this? It is here: \url{http://i.imgur.com/T1GvC5J.png}}}

\newthought{By the way}, a key feature of $\mathsf{FVect}$ (and more generally, all symmetric compact closed categories) is that processes, i.e. morphisms, $V\to W$ are in bijection with \textit{states} $\mathbb{R}\to W\otimes V^*\cong V^*\otimes W$, which is the special name given to morphisms whose domain is the monoidal unit. This bijection is sometimes called \textbf{process-state duality}, and in the context of $\mathsf{FVect},$ it means we can view linear maps as vectors in a tensor product\footnote{While a linear map $\mathbb{R}\to V^*\otimes W$ is not itself a vector in $V^*\otimes W$, it can be identified with one, namely with the image of 1 in $\mathbb{R}$! More generally, for any finite-dimensional vector space $A$ over $\mathbb{R}$, you can always think of $\hom(\mathbb{R},A)$ as $A$ itself, at least at the set level. That's because the forgetful functor $U\colon \mathsf{FVect}\to \mathsf{Set}$ is representable with representing object $\mathbb{R}$. In other words, linear maps $\mathbb{R}\to A$ are in one-to-one correspondence with the vectors in $A,$ viewed as elements of its underlying set, \[\hom(\mathbb{R},A)\cong UA\] This is completely analogous to how functions $\{\ast\}\to X$ from the one-point set to a set $X$ are in one-to-one correspondence with the elements in $X,$ \[\hom(\{\ast\},X)\cong X\] and is another manifestation of the ``probing'' idea we saw in the margin on page \pageref{fig:probe}.} and vice versa! 
\begin{center}
\includegraphics[width=!,totalheight=!,scale=0.11]{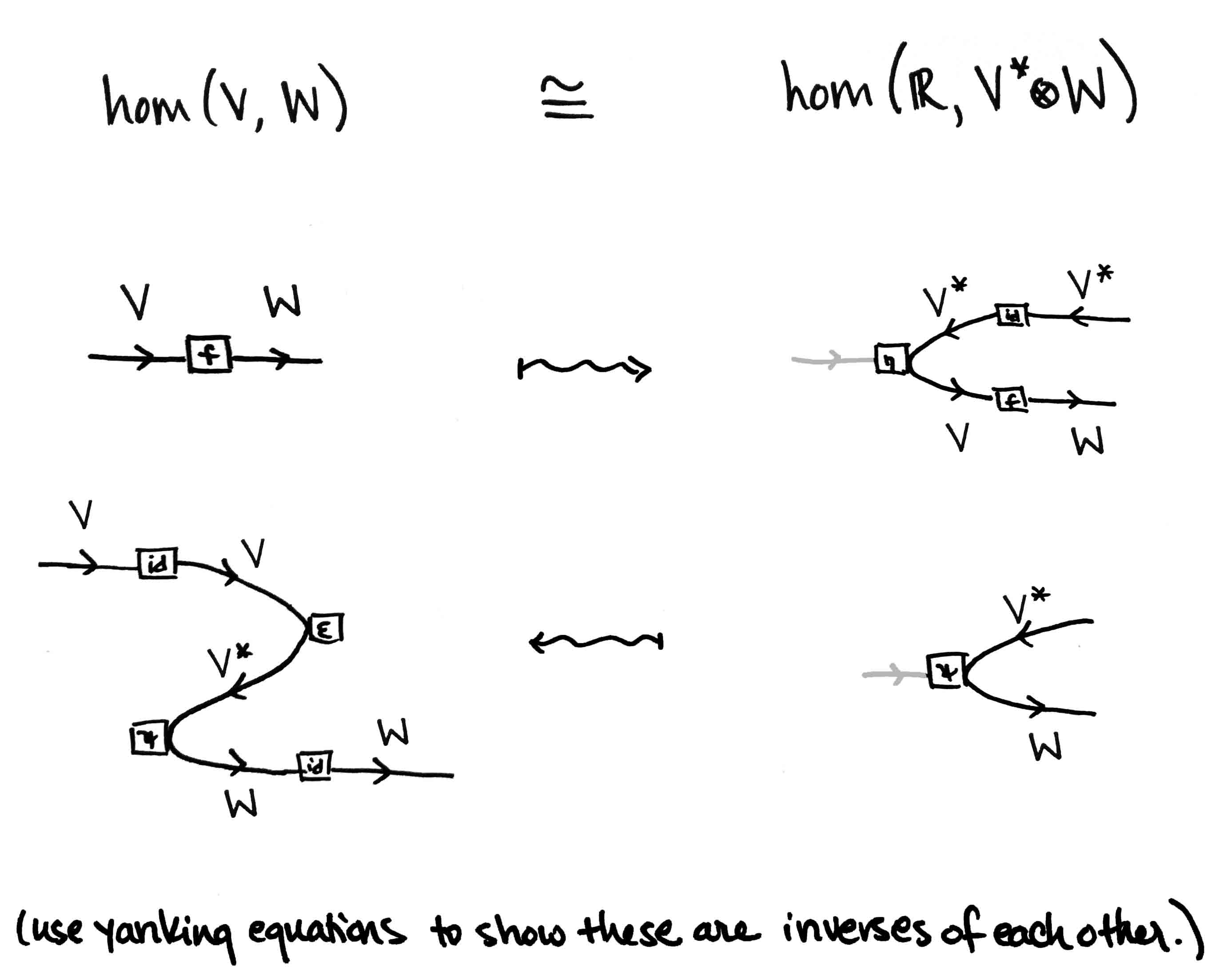}
\end{center}

I like to think of it this way: when $V$ and $W$ are $\mathbb{R}^n$ and $\mathbb{R}^m$ with the standard bases, process-state duality---taken together with the margin note on the previous page---is akin to the observation that matrices can be viewed as vectors and vice versa. That is, a linear map $\mathbb{R}^n\to \mathbb{R}^m$ has an $m\times n$ matrix representation which can be reshaped into an $nm\times 1$ column vector and then identified with a vector in $\mathbb{R}^n\otimes \mathbb{R}^m$. Conversely, there's a way to identify a vector in $\mathbb{R}^n\otimes \mathbb{R}^m$ with an $nm\times 1$ column vector that can be reshaped into an $n\times m$ matrix, which gives rise to a linear map $\mathbb{R}^n\to\mathbb{R}^m$.
		\marginnote{Here we're making the mental identification
	\begin{align*}
	\text{linear map} &\rightsquigarrow \text{process}\\
	\text{vector} &\rightsquigarrow \text{state}
	\end{align*} 
	which is closely related to the process-state duality seen in quantum physics. It's not quite the same, though---we'd need to replace $\mathsf{Mat}(\mathbb{R})$ by the category of completely positive maps, another compact closed category! For more, see Example 2.4 of \href{https://arxiv.org/pdf/1701.04732v4.pdf}{``A categorical semantics for causal structure''} by Aleks Kissinger and Sander Uijlen, as well as section 4.1.2 and chapter 6 of \href{https://www.amazon.com/Picturing-Quantum-Processes-Diagrammatic-Reasoning/dp/110710422X}{\textit{Picturing Quantum Processes}} by Bob Coecke and Aleks Kissinger.}
	\begin{quote}
	\textit{Aside:} There is, I think, a nice categorical way to piece this together. First note that there is a category $\mathsf{Mat}(\mathbb{R})$ whose objects are natural numbers $n,m,\ldots$ and whose morphisms $f\colon n \to m$ are $m\times n$ matrices with real entries. The identity $n\to n$ is the $n\times n$ identity matrix and composition is given by matrix multiplication. This category is actually a compact closed category! The monoidal product on objects is given by multiplication $n\otimes m:=nm$ and on morphisms is given by the Kronecker product of matrices. The monoidal unit is $1\in \mathbb{N}$. For the compact closed structure, each object is self-dual, $n^*:=n$, and for each $n\in\mathbb{N}$ the unit map $\eta_n\colon 1\to n^2$ is the $n^2\times 1$ column vector obtained by stacking the standard bases vectors $\mathbf{e}_1,\ldots,\mathbf{e}_n$ on top of each other. In other words, $\eta$ is given by the Kronecker delta function $\eta_{ij}:=\delta_{ij}$. The counit map $\epsilon_n\colon n^2\to 1$ is the $1\times n^2$ row vector obtained by taking the transpose of $\eta_n$. For example, if $n=3$ then $\eta_3=\begin{bmatrix}1&0&0&0&1&0&0&0&1\end{bmatrix}^\top$
		\marginnote{So the entries of $\eta_3$ are $\begin{bmatrix}\delta_{11}&\delta_{12}&\delta_{13}&\delta_{21}&\delta_{22}&\delta_{23}&\delta_{31}&\delta_{32}&\delta_{33}\end{bmatrix}.$}
	 and $\epsilon_3=\eta_3^\top.$ Since $\mathsf{Mat}(\mathbb{R})$ is compact closed, it exhibits process-state duality, too:
	\[\hom(n,m)\cong \hom(1,nm)\]
	This correspondence is precisely the reshaping of $n\times m$ matrices into $nm\times 1$ column vectors and vice versa, which can be verified by using the unit and counit maps in a way analogous to the string diagrams shown at the bottom of the previous page. To tie this in to the remark about $\mathbb{R}^n$ and $\mathbb{R}^m$ in the previous paragraph, note that there is a functor $\mathsf{FVect}\overset{\cong}{\to}\mathsf{Mat}(\mathbb{R})$ sending a vector space $V$ to its dimension $\text{dim}(V)$ and a linear map $f\colon V\to W$ to its corresponding $\text{dim}(W)\times \text{dim}(V)$ matrix representation, and it defines an equivalence of categories! For details, see the discussion on page 30 as well as Corollary 1.5.11 of Emily Riehl's \href{http://www.math.jhu.edu/~eriehl/context.pdf}{\textit{Category Theory in Context}}. 
	\end{quote}

As we'll see in Section \ref{sec:Second}, process-state duality pairs very nicely with our intuition about language. There we'll discover that a \textit{verb} can either be represented as a vector in a tensor product of vector spaces \textit{or} as a linear map, i.e. a process. In other words, a verb is an action in the eyes of both grammar and mathematics!

\subsection*{Another digression: A conjunction with adjunctions?}
If you're familiar with adjunctions in category theory, then you might wonder about this choice of naming and notation:
\[\eta \rightsquigarrow \textbf{``unit''}\qquad\qquad\qquad \epsilon\rightsquigarrow \textbf{``counit''}\]
Is it a coincidence that these two words are \textit{also} used in the definition of an \textit{adjunction}? NOPE. They are closely related. Specifically, the data $V,V^*,\eta,$ and $\epsilon$ together with the yanking equations \textit{are} an instance of a categorical adjunction! I think this is a neat fact,\footnote{which appears on the first page of \href{https://ac.els-cdn.com/0022404980901012/1-s2.0-0022404980901012-main.pdf?_tid=6c45bfb7-9cd3-42b1-ad62-4726b10f48e8&acdnat=1528395229_68b91be8a4ba1266aee54424771f749e}{``Coherence for Compact Closed Categories''} by Kelley and LaPlaza.} so let's take yet another digression. Happily, it will tie in quite nicely with our discussion on string diagrams. We'll begin by recalling the definition of an adjunction.

\begin{defn}\label{def:adj} An \textbf{adjunction} between categories $\cat{C}$ and $\cat{D}$ is a pair of functors
		\marginnote[-1cm]{Equivalently, $L$ and $R$ form an adjunction if for all objects $c\in \cat{C}$, $d\in\cat{D}$ there is an isomorphism \[\text{hom}_\cat{D}(Lc,d) \overset{\cong}{\longleftrightarrow} \text{hom}_\cat{C}(c,Rd)\] that's natural in both $c$ and $d$.
		}
\[ \adj{ L }{ \cat{C} }{ \cat{D} }{ R } \]
and a pair of natural transformations
\[\eta\colon \id_\cat{C}\Longrightarrow RL\qquad\qquad \epsilon\colon LR\Longrightarrow \id_{\cat{D}}\]
called the \textbf{unit} and \textbf{counit} respectively, such that these two triangles commute:
		\marginnote[-1.5cm]{
		\begin{tikzpicture}[>=stealth,overlay,xshift=-15pt]
		\draw[->] (0,0) -- (-1,0);
		\end{tikzpicture}
		\hspace{-8pt}There, $\id_\cat{C}$ denotes the identity functor on $C.$ It assigns each object and morphism in $\cat{C}$ to itself.\\

		\vspace{0.3cm}
		\begin{tikzpicture}[>=stealth,overlay,xshift=-15pt]
		\draw[->] (0,0) -- (-1,-0.3);
		\end{tikzpicture}
		\hspace{-8pt} Here, $L\circ \eta$ denotes the natural transformation whose components are of the form $L\eta_c\colon c\to LRLc$, while $\epsilon \circ L$ is the natural transformation with components $e_{Lc}\colon LRLc\to Lc.$ A similar story holds for $\eta\circ R$ and $R\circ \epsilon.$ (As per the margin comment on page \pageref{fig:monfun1}, I'd prefer to omit the composition symbol $\circ,$ but I'm writing it now for good reason, as we'll soon see!)}
	\[
	\begin{tikzcd}
	L \arrow[r, "L\circ \eta", Rightarrow] \arrow[rd, "\id_L", Rightarrow, swap] & LRL \arrow[d, "\epsilon \circ L", Rightarrow] \\
	& L
	\end{tikzcd}
	\qquad\qquad\qquad
	\begin{tikzcd}
	R \arrow[r, "\eta \circ R", Rightarrow] \arrow[rd, "\id_R"', Rightarrow] & RLR \arrow[d, "R\circ \epsilon", Rightarrow] \\
	 & R
	\end{tikzcd}
	\]
The adjunction is denoted $L\dashv R,$ and $L$ is said to be \textbf{left adjoint} to $R$ while $R$ is said to be \textbf{right adjoint} to $L$.
\end{defn}

Believe it or not, those commuting triangles---often called the \textit{triangle identities}---are closely related to the yanking equations in (\ref{eq:yankvect})! Indeed, ``these triangles commute''  means that these two equations hold:
\begin{equation}
\begin{aligned}\label{eq:triangle}
(\epsilon\circ L)\circ(L\circ \eta)&=\id_L\\[7pt]
(R\circ \epsilon)\circ (\eta\circ R)&=\id_R
\end{aligned}
\end{equation}
Now lets compare them to (\ref{eq:yankvect}):
\begin{align*}
(\epsilon \otimes \id_V)\circ (\id_V\otimes \eta)&=\id_V\\[7pt]
(\id_{V^*}\otimes\:\epsilon)\circ (\eta\otimes \id_{V^*})&=\id_{V^*}
\end{align*}
Why are (\ref{eq:yankvect}) and (\ref{eq:triangle}) so similar? \textit{What's going on here?} Is there a sense in which a vector space and its dual form an \textit{adjunction}? 

\vspace{0.2cm}
\noindent Is $V\dashv V^*$ a thing?
\vspace{0.2cm}

\vspace{0.2cm}
\noindent \textit{Yes!}
\vspace{0.4cm} 

\noindent But to make sense of $V\dashv V^*$ we'll need to venture into the world of \textit{2-categories.} A 2-category is an appropriate setting in which to talk about adjunctions, among other things. Here's why. As we know from the definition above, an adjunction consists of 
\begin{enumerate}[i.]
\item some \textcolor{Orchid}{objects} (categories) \tikzmark{a}
	\marginnote{Recall: The data of an adjunction are functors between categories \[ \adj{ {\color{CornflowerBlue} L} }{ {\color{Orchid} \cat{C} }}{ {\color{Orchid} \cat{D} }}{ {\color{CornflowerBlue} R} } \]
	and natural transformations
	\[\eta\colon {\color{CornflowerBlue} \id_\cat{C}} \; {\color{Melon}\Longrightarrow} \; {\color{CornflowerBlue} RL}\qquad\qquad \epsilon\colon {\color{CornflowerBlue} LR}  \; {\color{Melon}\Longrightarrow} \; {\color{CornflowerBlue}\id_{\cat{D}}}\]}
\item some \textcolor{CornflowerBlue}{arrows} (functors) \;\: \tikzmark{b}
\item some \textcolor{Melon}{arrows \textit{between} the arrows} (natural transformations)
\end{enumerate} 
	\begin{tikzpicture}[overlay, remember picture,decoration={brace,amplitude=1ex}]
  	\draw[decorate,thick] (a.north east) -- (b.south east) node[midway, right=0.1cm] {} node[midway, right=0.5cm,text=black,text width = 2in] { \textbf{These form a category!}};
	 \end{tikzpicture}
You'll notice that the \textcolor{Orchid}{objects}  and the \textcolor{CornflowerBlue}{arrows} \textit{themselves} form a category, namely $\cat{Cat}$, the category of all categories. The objects of $\cat{Cat}$ are \textcolor{Orchid}{categories} and the morphisms are \textcolor{CornflowerBlue}{functors}.

\vspace{0.4cm}
\noindent \textit{Nice}.
\vspace{0.4cm}

\noindent It'd be even nicer, though, if the \textcolor{Melon}{natural transformations} were also part of the data. That is, it'd be \textit{super nice} if the \textcolor{Orchid}{thr}\textcolor{CornflowerBlue}{ees}\textcolor{Melon}{ome} itself constituted a known categorical construction. But as it stands, it doesn't. There is no room for a notion of ``arrows between arrows'' in the definition of a category. 

\vspace{0.4cm}
\noindent So what do we do? 
\vspace{0.4cm}

\noindent We expand the definition. \textit{Literally.} We add an extra dimension, which results in a 2-category. That is, a 2-category consists of
	\marginnote{What's more, in any 2-category there is a composition rule for 2-morphisms just like there is for 1-morphisms in an ordinary category! In fact, in a 2-category we require that the set $\text{hom}(\bullet,\circ)$ be \textit{more} than a set. We ask that it be a \textit{category} itself! Its objects are 1-cells $\bullet\to\circ$, and its morphisms are  2-cells
 		$\begin{tikzcd}[ampersand replacement =\&]
		\bullet
			\arrow[bend left=30]{r}[name=U]{}
			\arrow[bend right=30]{r}[name=D]{} \&
		\circ
			\arrow[shorten <=2pt,shorten >=1pt,Rightarrow,to path={(U) -- (D)}]{}
		\end{tikzcd}$
		There is also an identity 1-cell $\id_{\bullet}$ for each 0-cell $\bullet$ \textit{and} there is an identity 2-cell $\id_{\id_\bullet}$ for each $\id_\bullet$. Confusingly, both of these identity morphisms are sometimes denoted as $\bullet$. And of course, there are the usual identity and associativity axioms, though I won't write them here.}
	\begin{enumerate}[i.]
	\item objects, now called \textbf{\textcolor{Orchid}{0-cells}} \textcolor{Orchid}{$\bullet, \circ,\ldots$}
	\item morphisms that go between objects, which are the usual arrows, but now we'll call them \textbf{\textcolor{CornflowerBlue}{1-cells}}
	${\color{Orchid} \bullet}\;\; {\color{CornflowerBlue}\longrightarrow} \;\; {\color{Orchid}\circ}$ 
	\item morphisms that go between 1-cells, which are not surprisingly called \textbf{\textcolor{Melon}{2-cells}}
	$\begin{tikzcd}[column sep=huge]
	{\color{Orchid}\bullet}
	  \arrow[bend left=30,color=CornflowerBlue]{r}[name=U]{}
	  \arrow[bend right=30,color=CornflowerBlue]{r}[name=D]{} &
	{\color{Orchid}\circ}
	  \arrow[shorten <=4pt,shorten >=1pt,Rightarrow,to path={(U) -- (D)}, color=Melon]{}
	\end{tikzcd}$
	\end{enumerate}
As you might guess, the quintessential example of a 2-category is $\cat{Cat}$, where the
	\begin{enumerate}[i.]
		\item \textcolor{Orchid}{0-cells} are \textit{categories} ${\color{Orchid}\cat{C},\cat{D},\ldots}$
		\item \textcolor{CornflowerBlue}{1-cells} are \textit{functors} ${\color{Orchid}\cat{C}}\; {\color{CornflowerBlue}\overset{F}{\longrightarrow}}\; {\color{Orchid} \cat{D}}$
		\item \textcolor{Melon}{2-cells} are \textit{natural transformations}  
		$\begin{tikzcd}[column sep=huge]
		{\color{Orchid}\cat{C}}
		  \arrow[bend left=30, color=CornflowerBlue]{r}[name=U,label=above:$F$]{}
		  \arrow[bend right=30, color=CornflowerBlue]{r}[name=D,label=below:$G$]{} &
		{\color{Orchid}\cat{D}}
		  \arrow[shorten <=4pt,shorten >=1pt,Rightarrow,to path={(U) -- (D)}, color=Melon]{}
	\end{tikzcd}$
	\end{enumerate}
So a 2-category is a good generalization of the relationship we see exhibited among \textcolor{Orchid}{categories,} \textcolor{CornflowerBlue}{functors,} and \textcolor{Melon}{natural transformations}. Having generalized this trio, it becomes very easy to talk about ``adjunctions'' in any 2-category. Parallel to Definition \ref{def:adj}, we might lay down the following proposed definition:
	\begin{quote}
	\textbf{Definition (proposed).} An \textbf{adjunction} between 0-cells ${\color{Orchid}\bullet}$ and ${\color{Orchid}\circ}$ is a pair of 1-cells ${\color{CornflowerBlue}l}$ and ${\color{CornflowerBlue}r}$,
	\[
	\begin{tikzcd}[ampersand replacement=\&, column sep=4ex]
			{\color{CornflowerBlue}l} \colon	{\color{Orchid}\bullet} \ar[yshift=+.6ex, color=CornflowerBlue]{r}
			\& {\color{Orchid}\circ} \colon {\color{CornflowerBlue}r} \ar[yshift=-.4ex, color=CornflowerBlue]{l}
	\end{tikzcd}
	\]
	and a pair of 2-cells ${\color{Melon}\eta}$ and ${\color{Melon}\epsilon}$, called the \textbf{unit} and \textbf{counit} respectively, 
	\[ 
	\begin{tikzcd}[column sep=huge]
		{\color{Orchid}\bullet}
		  \arrow[bend left=30, color=CornflowerBlue]{r}[name=U,label=above:$\id_\bullet$]{} 
		  \arrow[bend right=30, color=CornflowerBlue]{r}[name=D,label=below:$r\circ l$]{}
		  &
		{\color{Orchid}\bullet}
		  \arrow[shorten <=4pt,shorten >=1pt,Rightarrow,to path={(U) -- node[label=right:$\eta$] {} (D)}, color=Melon]{}
	\end{tikzcd}
	\qquad\text{and}\qquad
	\begin{tikzcd}[column sep=huge]
		{\color{Orchid}\circ}
		  \arrow[bend left=30, color=CornflowerBlue]{r}[name=U,label=above:$\id_\circ$]{}
		  \arrow[bend right=30, color=CornflowerBlue]{r}[name=D,label=below:$l\circ r$]{} &
		{\color{Orchid}\circ}
		  \arrow[shorten <=4pt,shorten >=2pt,Rightarrow,to path={(D) -- node[label=right:$\epsilon$] {} (U)}, color=Melon]{}
	\end{tikzcd}
	\]
	such that these two triangles commute
	\[
	\begin{tikzcd}
	{\color{CornflowerBlue}l} \arrow[r, "l\circ \eta", Rightarrow, color=Melon] \arrow[rd, "\id_l", Rightarrow, swap, color=Melon] & {\color{CornflowerBlue}l\circ r\circ l } \arrow[d, "\epsilon \circ l", Rightarrow, color=Melon] \\
	& {\color{CornflowerBlue}l}
	\end{tikzcd}
	\qquad\qquad\qquad
	\begin{tikzcd}
	{\color{CornflowerBlue}r} \arrow[r, "\eta \circ r", Rightarrow, color=Melon] \arrow[rd, "\id_r"', Rightarrow, color=Melon] & {\color{CornflowerBlue}r\circ l\circ r} \arrow[d, "r\circ \epsilon", Rightarrow, color=Melon] \\
	 & {\color{CornflowerBlue}r}
	\end{tikzcd}
	\]
	i.e. such that the following equations hold\footnote[][-5cm]{
		If you do a Google search for ``definition of 2-category'' you'll soon find that 2-cells can be composed in two ways: ``vertically'' and ``horizontally.'' I didn't mention this earlier, but now's a good time to do so. Suppose we have three 1-cells from {\color{Orchid}$\bullet$} to {\color{Orchid}$\circ$} and 2-cells ${\color{Melon}\eta}$ and ${\color{Melon}\epsilon}$ as shown below on the left,
		\[ 
		\begin{tikzcd}[column sep=huge,ampersand replacement=\&]
		{\color{Orchid}\bullet}
		  \arrow[bend left=45, color=CornflowerBlue]{r}[name=U]{}
		  \arrow[color=CornflowerBlue]{r}[name=D]{}
		  \arrow[shorten <=1pt,shorten >=-3pt,Rightarrow,to path={(U) -- node[label=right:{\scriptsize$\eta$}] {} (D)}, color=Melon]{}
		  \arrow[bend right=45, color=CornflowerBlue]{r}[name=M]{}
		  \arrow[shorten <=1pt,shorten >=-3pt,Rightarrow,to path={(D) -- node[label=right:{\scriptsize$\epsilon$}] {} (M)}, color=Melon]{}
		   \&
		{\color{Orchid}\circ}
		\end{tikzcd}
		\quad \rightsquigarrow \quad
		\begin{tikzcd}[column sep=huge,ampersand replacement=\&]
		{\color{Orchid}\bullet}
		  \arrow[bend left=30, color=CornflowerBlue]{r}[name=U]{}
		  \arrow[bend right=30, color=CornflowerBlue]{r}[name=D]{}
		  \arrow[shorten <=1pt,shorten >=-3pt,Rightarrow,to path={(U) -- node[label=right:{\scriptsize$\epsilon\diamond\eta$}] {} (D)}, color=Melon]{}
		   \&
		{\color{Orchid}\circ}
		\end{tikzcd}
		\]
		then \textit{vertical} composition {\color{Melon}$\diamond$} gives a 2-cell {\color{Melon}$\epsilon\diamond \eta$} as shown above on the right. This is \textit{composition along a common 1-cell} {\color{CornflowerBlue}$\to$}. On the other hand, given four 1-cells as shown below left, 
		\[
		\begin{tikzcd}[ampersand replacement=\&]
		{\color{Orchid}\bullet}
			\arrow[bend left=40, color=CornflowerBlue]{r}[name=A]{}
			\arrow[bend right=40, color=CornflowerBlue]{r}[name=B]{}
			\arrow[shorten <=1pt,shorten >=-3pt,Rightarrow,to path={(A) -- node[label=right:{\scriptsize$\eta$}] {} (B)}, color=Melon]{} \&
		{\color{Orchid}\star}
			\arrow[bend left=40, color=CornflowerBlue]{r}[name=C]{}
			\arrow[bend right=40, color=CornflowerBlue]{r}[name=D]{}
			\arrow[shorten <=1pt,shorten >=-3pt,Rightarrow,to path={(C) -- node[label=right:{\scriptsize$\epsilon$}] {} (D)}, color=Melon]{} \&
		{\color{Orchid}\circ}
		\end{tikzcd}
		\quad \rightsquigarrow \quad
		\begin{tikzcd}[column sep=huge,ampersand replacement=\&]
		{\color{Orchid}\bullet}
		  \arrow[bend left=30, color=CornflowerBlue]{r}[name=U]{}
		  \arrow[bend right=30, color=CornflowerBlue]{r}[name=D]{}
		  \arrow[shorten <=1pt,shorten >=-3pt,Rightarrow,to path={(U) -- node[label=right:{\scriptsize$\epsilon\circ\eta$}] {} (D)}, color=Melon]{}
		   \&
		{\color{Orchid}\circ}
		\end{tikzcd}
		\]
		\textit{horizontal} composition {\color{Melon}$\circ$} gives a 2-cell {\color{Melon}$\epsilon\circ \eta$} as shown above right. This is \textit{composition along a common 0-cell} {\color{Orchid}$\star$}. Moreover, the triangle identities involve \underline{both} compositions. That is, the \textit{actual} equations are
			{\color{Melon}
			\[(\epsilon\circ l)\diamond (l\circ \eta)=\id_l 
			\quad\text{and}\quad
			(r\circ \epsilon)\diamond (\eta\circ r)=\id_r\]}
		Take note of the diamonds vs. the circles!}

	{
	\color{Melon}
	\[(\epsilon\circ l)\circ (l\circ \eta)=\id_l 
	\qquad\text{and}\qquad
	(r\circ \epsilon)\circ (\eta\circ r)=\id_r\]
	}
	where ${\color{Melon}l\circ \eta}:={\color{Melon}\id_{\color{CornflowerBlue}l}\circ\eta}$, and similarly for ${\color{Melon}\epsilon\circ l}$ and so on. We'll say $l$ is a \textbf{left adjoint} of $r$, and $r$ is a \textbf{right adjoint} of $l$, and we'll denote the adjunction by ${\color{CornflowerBlue}l}\dashv {\color{CornflowerBlue}r}$.
	\end{quote}

\vspace{0.4cm}
\noindent Alright, fine. \textit{But what does this have to do with vector spaces?}
\vspace{0.4cm}

\noindent The answer lies in the following neat fact.

\begin{center}
\noindent \textit{Neat Fact: Every monoidal category $(\cat{C},\otimes,1)$ can be viewed\\as a 2-category!}
\end{center}

\noindent Er, actually, I shouldn't spread rumors. 

\begin{center}
\noindent \sout{\textit{Neat Fact: Every monoidal category $(\cat{C},\otimes,1)$ can be viewed}}\\\sout{\textit{as a 2-category!}}
\end{center}

\noindent Here's the \textit{correct} statement:

\begin{center}
\noindent \textbf{Neat Fact: Every monoidal category $(\cat{C},\otimes,1)$ can be viewed\\as a bicategory!}
\end{center}

\noindent A bicategory is basically a 2-category---the data is completely the same. There are 0-cells, 1-cells, and 2-cells. The only difference is what's in the margin.\footnote[][-3cm]{In a 2-category, the composition of 1-cells is associative, i.e. $f(gh)=(fg)h$ for any composable triple of 1-cells $f,g,h.$ In a \textit{bicategory}, however, we weaken this. Instead of asking for equality, we ask for the existence of an invertible 2-cell  $f(gh)\Longleftrightarrow (fg)h$. As we'll see below, the category $\cat{FVect}$ gives rise to a bicategory rather than a 2-category because the two vector spaces $V\otimes (W\otimes U)$ and $(V\otimes W)\otimes U$ are not \textit{equal,} but there certainly is a linear isomorphism \[V\otimes (W\otimes U)\overset{\cong}{\longleftrightarrow} (V\otimes W)\otimes U!\] }
So any monoidal category $(\cat{C},\otimes,1)$ \textit{gives rise to} a bicategory $\mathcal{C}$ where the
\begin{enumerate}[i.]
\item only \textcolor{Orchid}{0-cell} is the \textit{category} $\cat{C}$
\item \textcolor{CornflowerBlue}{1-cells} are the \textit{objects} of $\cat{C}$; composition is  $\otimes$
\item \textcolor{Melon}{2-cells} are the \textit{morphisms} of $\cat{C}$; composition is composition $\circ$ in $\cat{C}$
\end{enumerate}
Therefore it makes sense to talk about 1-cells in $\mathcal{C}$ (i.e. \textit{objects} in $\cat{C}$) having adjoints! And it makes sense to talk about 2-cells in $\mathcal{C}$ (i.e. \textit{morphisms} in $\cat{C}$) being units and counits of the adjunction, vis-a-vis our Proposed Definition! In particular, this is true of the symmetric monoidal category $(\cat{FVect},\otimes,\mathbb{R})$. It gives rise to a bicategory $\mathcal{F}\mathsf{Vect}$ where the
\begin{enumerate}[i.]
\item only \textcolor{Orchid}{0-cell} is the category $\cat{FVect}$ 
\item \textcolor{CornflowerBlue}{1-cells} are \textit{vector spaces}; composition is the tensor product $\otimes$
\item \textcolor{Melon}{2-cells} are \textit{linear maps}; composition is the usual composition $\circ$
\end{enumerate}

So there is an adjunction of vector spaces $V\dashv W$ whenever the conditions of our Proposed Definition hold. Of course, those conditions hold \textit{precisely} when $W=V^*$ and $\eta$ and $\epsilon$ are defined as in (\ref{lis:ucounit}). Explicitly:
	\begin{quote}
	There is an \textbf{adjunction} ${\color{CornflowerBlue}V}\dashv {\color{CornflowerBlue}V^*}$ in the bicategory $\mathcal{F}\mathsf{Vect}$ since there are linear maps
		\marginnote[-2cm]{You'll notice that the monoidal unit $\mathbb{R}$ is taking the place of $\id_\bullet$ in the Proposed Definition. Indeed, $\id_\bullet$ and $\mathbb{R}$ are comparable since both are 1-cells that act as an identity on other 1-cells: For all  1-cells $\bullet\overset{f}{\longrightarrow}\circ$ in $\mathcal{C}$
		\[f\circ \id_\bullet=f=\id_\circ\circ f\]
		and for all vector spaces $V$ in $\mathcal{F}\mathsf{Vect},$
		\[V\otimes \mathbb{R}\cong V\cong \mathbb{R}\otimes V.\]
		}
	\[ 
	{\color{CornflowerBlue}\mathbb{R}} {\color{Melon}\overset{\eta}{\longrightarrow}} {\color{CornflowerBlue}V^*\otimes V}
	\qquad\text{and}\qquad
	{\color{CornflowerBlue}V\otimes V^*} {\color{Melon}\overset{\epsilon}{\longrightarrow}} {\color{CornflowerBlue}\mathbb{R}}
	\]
	so that the following triangles commute
	\[
	\begin{tikzcd}
	{\color{CornflowerBlue}V\cong V\otimes\mathbb{R}} \arrow[r, "V \otimes \eta", Rightarrow, color=Melon] \arrow[rd, "\id_V", Rightarrow, swap, color=Melon] & {\color{CornflowerBlue}V\otimes V^*\otimes V } \arrow[d, "\epsilon \otimes V", Rightarrow, color=Melon] \\
	& {\color{CornflowerBlue}\mathbb{R}\otimes V\cong V}
	\end{tikzcd}
	\qquad\qquad\qquad
	\begin{tikzcd}
	{\color{CornflowerBlue}V^*\cong\mathbb{R}\otimes V^*} \arrow[r, "\eta \otimes V^*", Rightarrow, color=Melon] \arrow[rd, "\id_{V^*}"', Rightarrow, color=Melon] & {\color{CornflowerBlue}V^*\otimes V\otimes V^*} \arrow[d, "V^*\otimes \: \epsilon", Rightarrow, color=Melon] \\
	 & {\color{CornflowerBlue}V^*\otimes\mathbb{R}\cong V^*}
	\end{tikzcd}
	\]
	\marginnote{
		\begin{tikzpicture}[>=stealth,overlay,xshift=-15pt]
		\draw[->] (0,0) -- (-.7,.2);
		\end{tikzpicture}
		\hspace{-8pt} {\color{black}On the leftmost triangle, the notation $V\otimes \eta$ denotes the linear map
		\[\id_V\otimes \eta\colon V\otimes \mathbb{R}\to V\otimes V^*\otimes V\]
		\begin{tikzpicture}[>=stealth,overlay,xshift=-8pt]
		\draw[->] (0,0) -- (-.9,.8);
		\end{tikzpicture}
		\hspace{-5pt} that appears in the first equation. A similar statement holds for $\epsilon\otimes V$, etc. Also, take note of the different symbols \textcolor{Melon}{$\circ$} and \textcolor{Melon}{$\otimes$} and compare them with the diamond \textcolor{Melon}{$\diamond$} and circle \textcolor{Melon}{$\circ$} in the margin on the previous page.}
		}
	i.e. so that the following equations hold
	{
	\color{Melon}
	\[(\epsilon\otimes\id_V)\circ (\id_V\otimes \eta)=\id_V 
	\qquad\text{and}\qquad
	(\id_{V^*}\otimes\: \epsilon)\circ (\eta\otimes \id_{V^*})=\id_{V^*}\]
	}
	and these are precisely the string diagram equations shown in the chart on page \pageref{fig:yank}.
	\end{quote}

\begin{center}
\noindent Voila!
\end{center}

Finally, notice that the above holds for \textit{every} vector space $V$ in $\mathsf{FVect.}$ On the other hand, there are certainly 2-categories in which not every 1-cell is \textit{dualizable}, i.e. has an adjoint. Take $\cat{Cat}$ for instance! Not every functor is part of an adjunction. There \textit{is}, however, a special name given to those bicategories $\mathcal{C}$ that do arise from a monoidal category $\mathsf{C}$ and in which every 1-cell has an adjoint. 
	\marginnote[-0.5cm]{The punchline for this section is that monoidal categories are an appropriate framework for \textit{stacking things together,} and the calculus of string diagrams allows us to replace complicated, messy equations by simple, neat pictures. In Section \ref{sec:Examples}, we'll see two examples of how this can be put into practice.}

\vspace{0.4cm}
\noindent That name is \textbf{compact closed.}

\newpage
\subsection{Decorated Cospans}\label{sec:cospan}
A second construction that appears in some work within applied category is the \textit{decorated cospan.} In any category, a diagram that looks like
\[A\to C\leftarrow B.\]
is called a \textbf{cospan}. In the next section, we'll only consider the case when $A,B,$ and $C$ are finite sets and the arrows are functions between them. A \textbf{decorated cospan} is a cospan where the middle set $C$ has been endowed with some extra structure. That's the intuitive definition, though I'd like to postpone a more precise definition until the next section.

Now you might think it strange to give a name to a simple diagram like $A\to C\leftarrow B$, but cospans come in handy quite often! For instance, if for some reason you can't possibly hope to find a morphism between objects $A$ and $B$, a common technique\footnote{I learned this from Brendan during the 2018 ACT workshop. Thanks, Brendan!} is to instead look for a ``larger'' object $C$ that ``contains'' both $A$ and $B$. Then although you don't have maps between $A$ and $B$, you \textit{do} have maps $A\to C\leftarrow B$. In that case, your cospan is the \href{https://twitter.com/math3ma/status/989233108137533442}{next best thing}.

\newthought{Admittedly,} this section is \textit{bite-sized} compared to the behemoth on monoidal categories that we just finished, but that's not because decorated cospans are any less important! In fact, Brendan Fong developed the theory of decorated cospans as part of his PhD thesis \href{https://arxiv.org/pdf/1609.05382.pdf}{``The Algebra of Open and Interconnected Systems''}, which has served as the foundation for \textit{much} progress in applied category theory, as I mentioned earlier. But in these notes, we'll only use the cospan construction in our brief discussion on chemical reaction networks in Section \ref{sec:First}. On the other hand, we will need the language of monoidal categories in both Sections \ref{sec:First} and \ref{sec:Second}. In fact, as we'll soon see, decorated cospans \textit{themselves} form a monoidal category!

\newpage
\subsection{Further Reading}
\textbf{For more on monoidal categories and string diagrams:}
\begin{itemize}
	\item Read Chapters 3 and 4 of \href{https://www.amazon.com/Picturing-Quantum-Processes-Diagrammatic-Reasoning/dp/110710422X}{\textit{Picturing Quantum Processes}} by Bob Coecke and Aleks Kissinger. There you'll also find more information on the interpretation of morphisms in a monoidal category as \textit{processes} and objects as \textit{systems}.
	\item Take a look at \href{https://www.youtube.com/watch?v=USYRDDZ9yEc&list=PL50ABC4792BD0A086}{TheCatsters videos} on string diagrams, by Eugenia Cheng and Simon Willerton. On second thought, \href{http://www.simonwillerton.staff.shef.ac.uk/TheCatsters/}{their entire collection} of videos is great. Go watch them all!
	\item If you like $\infty$-categories, you'll be delighted to know that a version of string diagrams (affectionately called ``strictly undulating squiggles'') and the yanking equations (!) make an appearance in chapter 8 of \href{http://www.math.jhu.edu/~eriehl/ICWM.pdf}{\textit{Elements of $\infty$-Category Theory}}, a new book on model-independent $\infty$-category theory by Emily Riehl and Dominic Verity.
\end{itemize}
\textbf{For more on decorated cospans:}
\begin{itemize}
	\item Read \href{https://golem.ph.utexas.edu/category/2016/01/decorated_cospans.html}{``Decorated Cospans''} a blog post by John Baez on the $n$-Category Caf{\'e}.
	\item Read Chapter 6 of \href{https://arxiv.org/pdf/1803.05316.pdf}{\textit{Seven Sketches in Compositionality}} by Brendan Fong and David Spivak. In Section 6.1, the authors give the following bit of motivation:
		\begin{quote}
		...we produce a certain monoidal category---namely that of \textit{cospans in} [a category] $\cat{C}$, denoted \textbf{Cospan}$_\cat{C}$---that can conveniently package $\cat{C}$'s colimits in terms of its own basic operations: composition and monoidal structure. In summary, the first part of this chapter is devoted to the slogan `colimits model connection.' (emphasis theirs)
		\end{quote}
	As we'll see in Section \ref{sec:First}, objects in \textbf{Cospan}$_\cat{C}$ are cospans in $\cat{C}$ and a morphism between two cospans is given by a construction called a \textit{colimit}. Like the composition $\circ$ and product $\otimes$ in a general monoidal category, a colimit is a categorical construction that allows you to connect things together. But for the sake of ``time'' (i.e. so that this document doesn't accidentally turn into a \textit{book}...), I'll assume familiarity with colimits. But if you'd like to an intuitive introduction of colimits, as well as their dual construction, limits, I recommend that you
	\item Take a look at \href{http://www.math3ma.com/mathema/2018/1/2/limits-and-colimits-part-1}{``Limits and Colimits (Part 1)''} a blog post on Math3ma. Also see chapters 3 and 6 of \href{https://arxiv.org/pdf/1803.05316.pdf}{\textit{Seven Sketches}} and chapter 3 of \href{http://www.math.jhu.edu/~eriehl/context.pdf}{\textit{Category Theory in Context}} by Emily Riehl.
\end{itemize}

\newpage
%%%%%%%%%% SECTION 3 %%%%%%%%
\section{Two Examples}\label{sec:Examples}
Having taken a leisurely stroll through two themes (functorial semantics and compositionality) and two constructions (monoidal categories and decorated cospans) within applied category theory, it's time to see them come to life in two examples. As mentioned in the introduction, we'll walk through the first example---chemical reaction networks---relatively quickly. There are several excellent resources available online, including John Baez's expositions on the $n$-Category Caf{\'e} as well as on his personal webpage. (I've included a few links to these in Section \ref{ssec:furtherreading3}.) Afterwards we'll take a longer stroll through the second example---natural language processing---in Section \ref{sec:Second}.

\subsection{Chemical Reaction Networks}\label{sec:First}
The first example comes from a paper by John Baez and Blake Pollard called \href{https://arxiv.org/abs/1704.02051}{``A Compositional Framework for Reaction Networks.''} Specifically, they provide a compositional framework for modeling \textit{chemical} reaction networks. A \textbf{chemical reaction network} is, well, a network of chemical reactions. And a \textbf{chemical reaction} is exactly what you think it is. It's what you learned back in high school: You start with some \textcolor{BlueGreen}{reactants} and some \textcolor{BlueGreen}{products}, and there's a \textcolor{BurntOrange}{chemical} process that takes one to the other.
\begin{center}
\includegraphics[width=!,totalheight=!,scale=0.2]{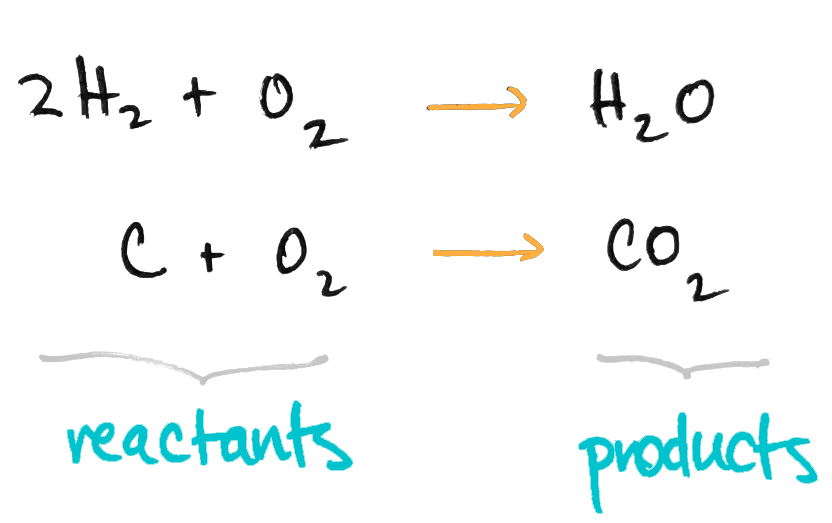}
\end{center}
What's nice is that these reactions can be depicted graphically:
\begin{center}
\includegraphics[width=!,totalheight=!,scale=0.1]{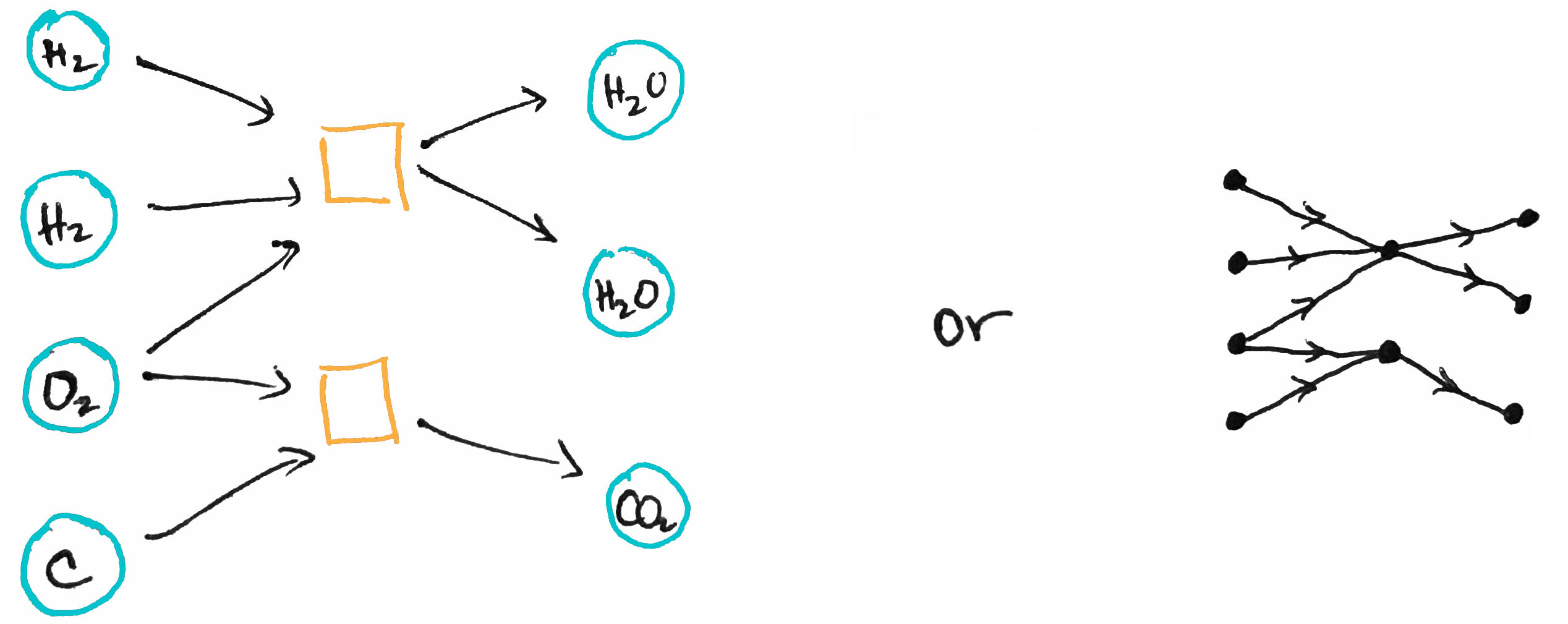}
\end{center}
Of course, you can imagine that there might be \textit{lots} of various \textcolor{BlueGreen}{reactants, products}, and \textcolor{BurntOrange}{chemical processes}. The corresponding \textit{network} would then be a (possibly huge) collection of these graphs stacked side-by-side, perhaps with connecting edges and loops and so on. For instance, this chemical reaction network made a cameo appearance in Baez's 2016 talk \href{https://www.youtube.com/watch?v=IyJP_7ucwWo}{``The Mathematics of Networks''}:
\begin{center}
\includegraphics{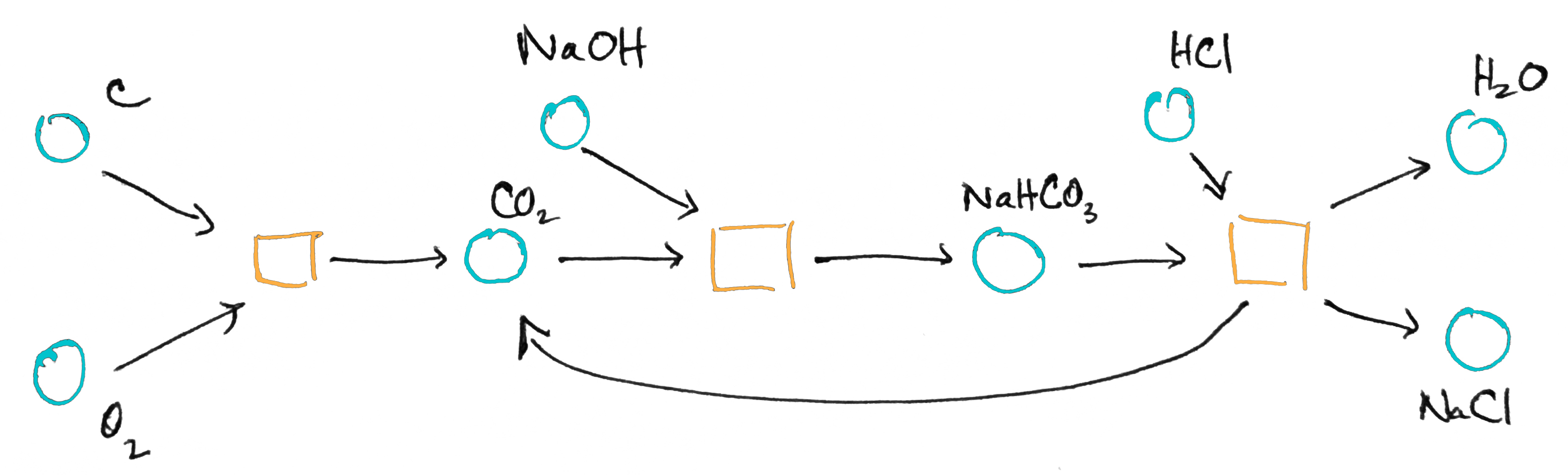}
\end{center}
Graphs such as these are examples of Petri nets. A \textbf{Petri net} is essentially a bipartite directed (multi)graph that allows us to visually represent reactions, though they are used outside of chemistry as well. 

But if we do wish to model chemical reactions, then an important thing we'd like to account for is the \textit{rate} at which one or more chemicals change over to another. A Petri net with rates included is called, appropriately, a  \textbf{Petri net with rates}. More specifically, it's a bipartite directed graph whose two types of vertices are called \textcolor{BlueGreen}{places}, which represent chemical species, and \textcolor{BurntOrange}{transitions}, which represent chemical reactions. Moreover, each transition $\tau_i$ is assigned a \textit{rate} $r_i$, a positive real number that describes how fast or how likely it is for $\tau_i$ to occur. These rates then allow us to write down differential equations that describe the system. A Petri net with rates is thus a pictorial representation of a set of differential equations that describe a system. So, for instance, if you \textit{did} watch Baez's \href{https://www.youtube.com/watch?v=IyJP_7ucwWo}{``The Mathematics of Networks''} talk then this example will look familiar:
\begin{center}
\includegraphics[width=!,totalheight=!,scale=0.13]{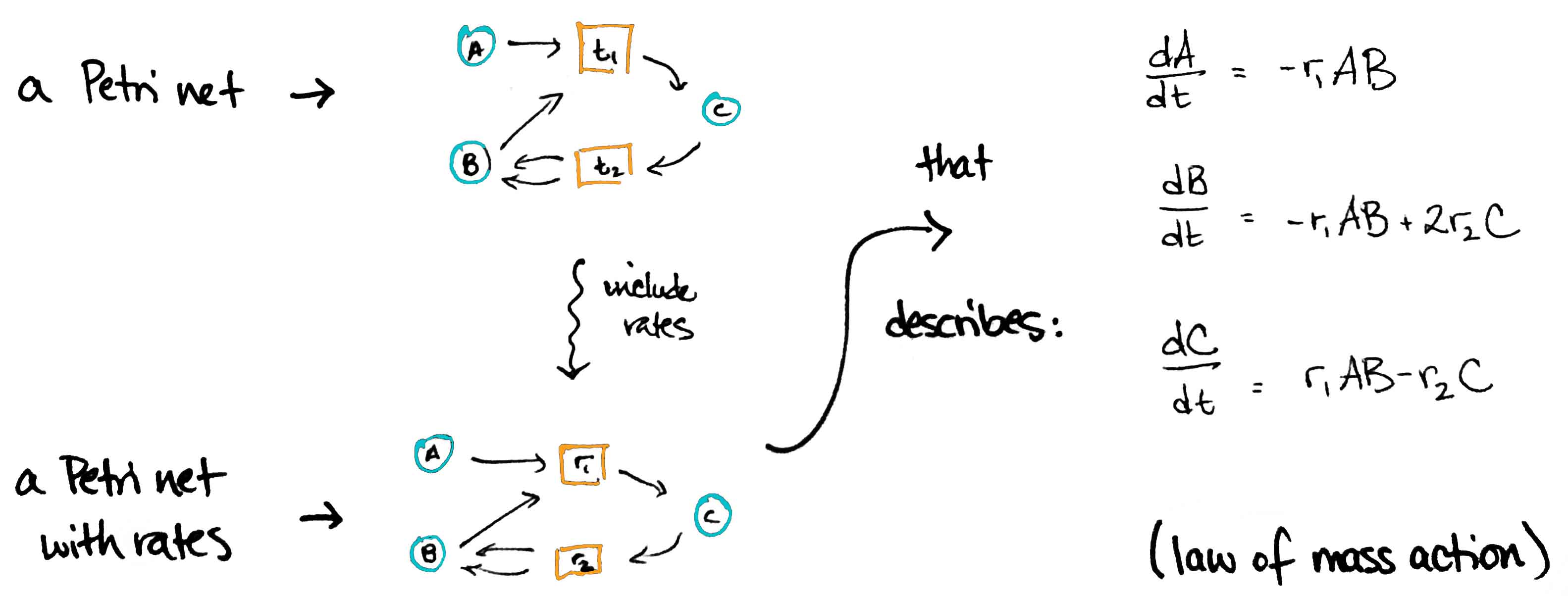}
\end{center}
It tells us that, for example, substances with concentrations $A$ and $B$ combine and produce a substance with concentration $C$ at a rate proportional to $r_1.$ The differential equations you see are due to the \textbf{law of mass action}, which says that the rate with which a chemical reaction will occur is equal to its rate constant $r_i$ multiplied by the product of the concentration of the reactants, i.e. the concentration of the ``inputs'' of the reaction.

By the way, the rates themselves could change with time, which might suggest the presence of a \textbf{dynamical system.}  What's more, a system such as the above could potentially interact with its environment, which is to say there might be some quantities that \textcolor{YellowGreen}{flow in} and some quantities that \textcolor{YellowGreen}{flow out}, resulting in an \textbf{open Petri net with rates}:
\begin{center}
\includegraphics[width=!,totalheight=!,scale=0.1]{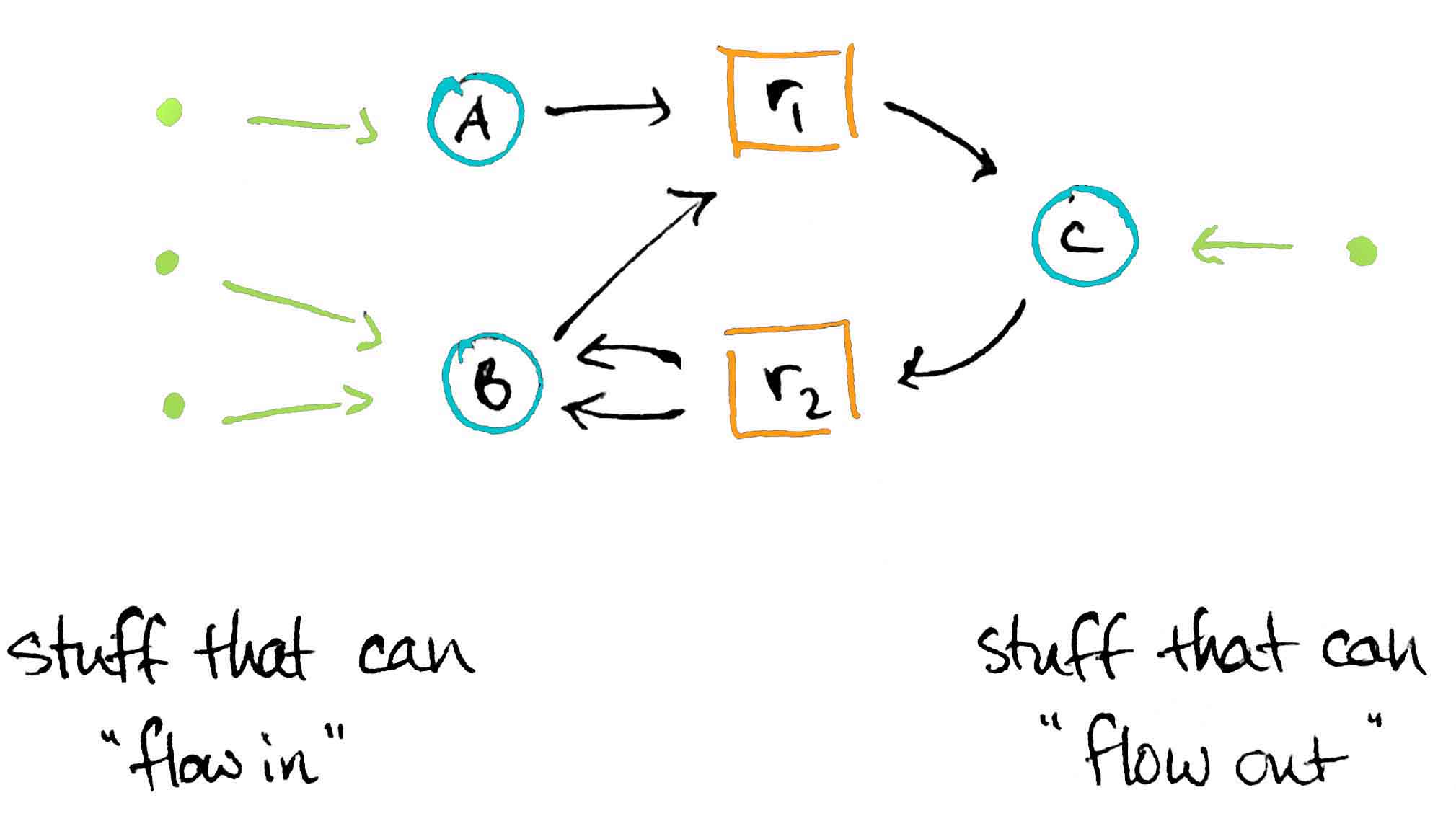}
\end{center}
These quantities can be incorporated into the equations, too, resulting in an \textbf{open dynamical system.} As Baez and Pollard summarize, the goal is to use these observations ``to build up a reaction network from smaller pieces, in such a way that its rate equation can be determined from those of the pieces.'' This is what is meant by a \textit{compositional} framework, and is a prime is example of the principal of compositionality mentioned in Section \ref{sec:Comp}.
What's more, a key step towards achieving this goal is given by functorial semantics! That is, we start by thinking of a Petri net as \textbf{syntax} and a set of differential equations as  \textbf{semantics}. And if we have a \textit{collection} of Petri nets that model a very large network, then---guided by the principal of compositionality---we would like to compose them by gluing graphs together, and we would like to aggregate them by stacking graphs on top of each other. In other words, we hope that Petri nets form a monoidal category! Similarly, one would hope that there is a sense in which dynamical systems form a monoidal category so that differential equations can be ``composed'' and ``aggregated'' as well. One would \textit{also} wish for a monoidal functor $\cat{Petri nets}\to\cat{dynamical systems}.$ 

That's a lot of wishes, but amazingly they all come true, for this is precisely what Baez and Pollard proved in their paper! But how exactly? How was it all formalized? The key is the decorated cospan construction of Brendan Fong that we mentioned in Section \ref{sec:cospan}. (What's amazing is that Fong's construction is general enough to model other \textbf{open reaction networks}\footnote[][-2cm]{This is the catch-all phrase for a network that interacts with its environment so that stuff can either flow in or flow out} as well! But more on that later.)
\begin{center}
\includegraphics{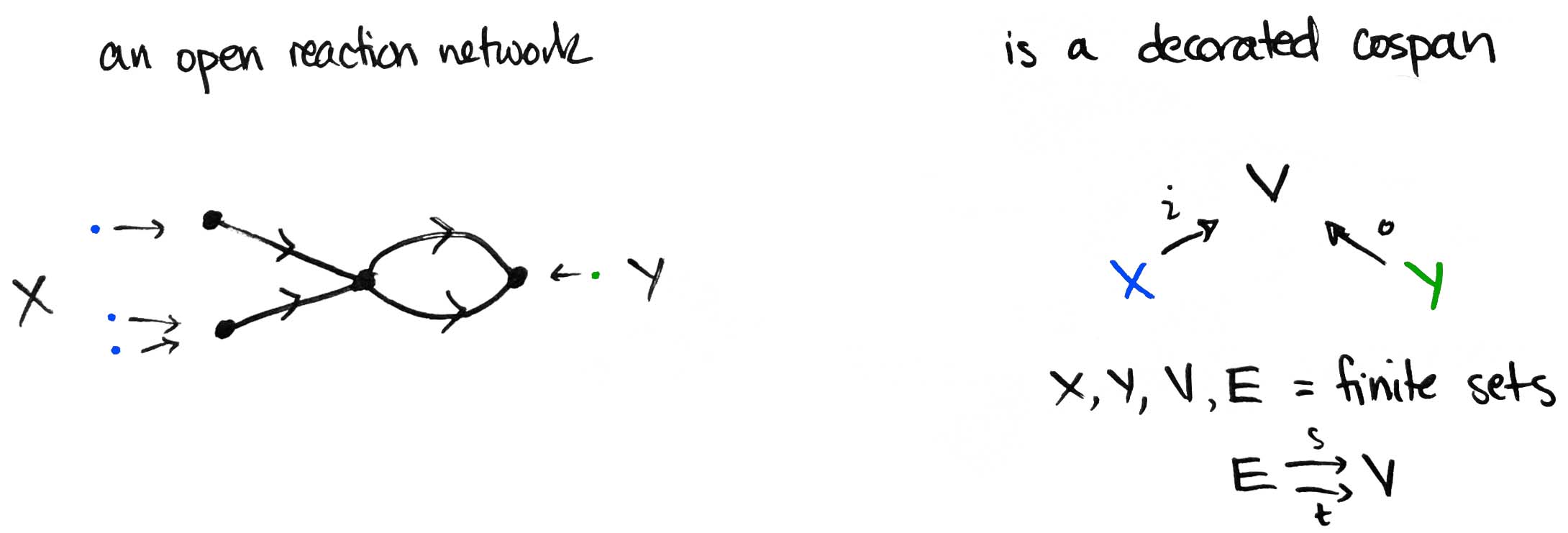}
\end{center}
	\marginnote[-4cm]{The reason the arrows from $X$ and $Y$ point \textit{in} is that $X$ and $Y$ might be thought of as ``leftputs and rightputs'' rather than as inputs and outputs. In other words, you'd like the freedom to think of things as flowing either in \textit{or} out of either end. For example, a physical pipe doesn't know the different between left and right. Water can flow in or out at either direction. As Baez \href{https://golem.ph.utexas.edu/category/2017/07/a_compositional_framework_for_2.html}{notes}, ``The main reason for these designations is to remember that when we screw together two pipes, we attach the output of the first to the input of the second.''} 
A \textbf{cospan} in the category of \textbf{finite sets}, for example, is just a diagram of the form $X\to V\leftarrow Y$, where we're meant to think of $V$ as the set of places (i.e. chemical species) in the Petri net.
% \[
% \begin{tikzcd}
%  & V &  \\
% X \arrow[ru, "i"] &  & Y \arrow[lu, "o"']
% \end{tikzcd}
% \]
To account for the edges in the graph, we ask that $V$ is ``decorated'' with extra structure, namely source and target maps from the set of edges. This results in a \textbf{decorated cospan,}
% $\begin{tikzcd}
% 	  E		\ar[yshift=+.5ex]{r}{s}
% 			\ar[yshift=-.5ex]{r}[swap]{t}
% 	& V
% \end{tikzcd}$
and Fong proved that these constructions form a category! That is, there is a category where objects are finite sets $X,Y,\ldots$, and a morphism $X\to Y$ is a decorated cospan whose feet are $X$ and $Y.$ Composition is given by the \textit{pushout}\footnote[][-1.18cm]{A pushout is a type of colimit, a major construction in category theory that`s a bit like like mathematical glue. Anytime you mush two mathematical objects together---like the graphs in the picture---you've probably got a colimit construction. More intuition behind colimits and their dual construction, limits, can be found \href{https://www.math3ma.com/blog/limits-and-colimits-part-1}{on Math3ma}.}, which amounts to gluing graphs together. (This composition is only associative up to isomorphism, so the morphisms are really \textit{isomorphism classes} of cospans. Also, the identity $X\to X$ is the empty graph.)
\begin{fullwidth}
\vspace{0.1cm}
\begin{center}
\includegraphics{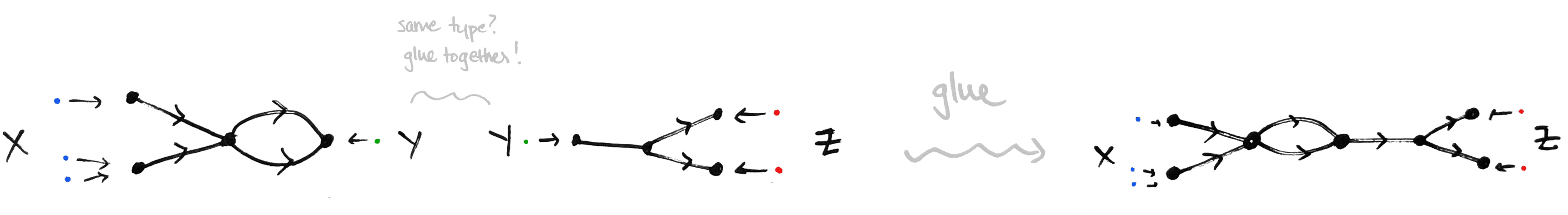}
\end{center}
\vspace{0.1cm}
\end{fullwidth}
What's more, Fong showed that this category has a \textbf{symmetric monoidal} structure by stacking graphs on top of each other, i.e. by taking their disjoint union. In fact, it's \textbf{compact closed} and also a \textbf{hypergraph} category!\marginnote[-1cm]{A hypergraph category is a symmetric monoidal category in which every object has a special commutative Frobenius structure. This allows more freedom (i.e. messiness) when composing morphisms, reflecting the messiness of most network diagrams!}

In summary, Fong's constructions quickly give rise to key results, which I'll summarize here. The first is that the syntax category of Baez and Pollard is indeed a category:

\begin{thm}[Baez, Pollard] There is a symmetric monoidal category $\cat{Petri}$ where
\begin{itemize}
	\item \textbf{objects} are finite sets $X,Y,\ldots$
	\item a \textbf{morphism} $X\to Y$ is a \textbf{open Petri net with rates}, i.e. a cospan\marginnote{\textit{Really,} it's an isomorphism class of cospans. Also, you'll notice that in Theorem 12 of Baez and Pollard's \href{https://arxiv.org/pdf/1704.02051.pdf}{``A Compositional Framework,''} their syntax category is something called $\cat{RxNet}$. That stands for the category of \textit{open reaction networks with rates}. An open reaction network is very nearly the same as an open Petri net, though I'm glossing over this a bit.}
		\[
		\begin{tikzcd}
 		& V &  \\
		X \arrow[ru, "i"] &  & Y \arrow[lu, "o"']
		\end{tikzcd}
		\]
	together with a \textit{Petri net with rates} whose \textit{places} are comprised of $V.$
\end{itemize}
\end{thm}

The next corollary provides the same statement for the semantics category: 
\begin{thm}[Baez, Pollard] There is a symmetric monoidal category $\cat{Dynam}$ where
\begin{itemize}
	\item \textbf{objects} are finite sets $X,Y,\ldots$
	\item a \textbf{morphism} $X\to Y$ is an \textbf{open dynamical system,} i.e. a cospan\marginnote{Again, it's really an isomorphism class of cospans. And again we can think of $V$ as the set of all places in a Petri net where, as before, there may be a real number $r_i$ attached to each vertex than can vary with time. The description of \textit{how} these things vary in time is precisely a vector field on $\mathbb{R}^V$.}
		\[
		\begin{tikzcd}
 		& V &  \\
		X \arrow[ru, "i"] &  & Y \arrow[lu, "o"']
		\end{tikzcd}
		\]
	together with a \textit{smooth vector field} on $\mathbb{R}^V$.
\end{itemize}
\end{thm}

Finally, another result of Baez and Pollard shows the existence of a symmetric monoidal functor $\square\colon \cat{Petri}\to\cat{Dynam}$ from the syntax to the semantics. 
%(Amazingly, the functor comes for ``free'' from the groundwork that Fong laid!) 
	% \begin{quote}{\color{gray}
	% \begin{thm}[Fong] Suppose $\cat{C}$ is a category with finite colimits. Let $F,G\colon (\cat{C},\sqcup)\to (\cat{Set},\times)$ be two lax monoidal functors and let $\theta\colon F\Rightarrow G$ be a monoidal transformation. Then there exists a symmetric monoidal functor $T_\theta\colon F\cat{Cospan}\to G\cat{Cospan}$.
	% \end{thm}
	% So think of $\cat{Petri}$ as $F\cat{Cospan}$ and $\cat{Dynam}$ as $G\cat{Cospan}$ and $T_\theta$ as $\square.$}
	% \end{quote}
\begin{thm}[Baez, Pollard] There is a symmetric monoidal functor $\square\colon \cat{Petri}\to\cat{Dynam}$
	\marginnote{But see Theorem 18 of Baez and Pollard where, since $\cat{RxNet}$ is used in lieu of $\cat{Petri}$,  the symmetric monoidal functor $\cat{RxNet}\to\cat{Dynam}$ is a slightly different \textit{gray boxing functor}.}
sending any \textit{open Petri net with rates} to the corresponding \textit{open dynamical system.} 
\end{thm}
The upshot is that functoriality and monoidality
\begin{align*}
\square(f\circ g)&=\square f\circ \square g\\
\square(f\otimes g)&=\square f\otimes \square g
\end{align*}
tell us that if you want to understand the open dynamical systems of the composite (or tensor product) of two open Petri nets, then you just have to find the open dynamical systems of each one and then compose (tensor). This is \textit{exactly} the \textbf{principle of compositionality}: to determine the behavior of a big complicated thing, you  need only understand the behaviors of its components, and then assemble them together. And by the way, this works for many other kinds of network graphs, not just Petri nets. It's all part of Baez's larger body of work on a general \href{http://math.ucr.edu/home/baez/networks/}{categorical framework for a theory of networks} which encompasses electrical circuits, Markov processes, signal-flow graphs in control theory, and more! 

\newthought{This rapid tour} through chemical reaction networks is only \textit{one} way that compositionality, functorial semantics, and monoidal categories are being used in applications. The next example gives a second way: natural language processing.

\newpage
%%%%%%%%%%%%% SECTION 3 %%%%%%%%%%%%%%
\subsection{Natural Language Processing}\label{sec:Second}
At long last, we've made it to our second application of category theory---natural language processing! It is, simply put, a branch of artificial intelligence that aims to train computers to understand human language. What's nice is that computers can understand meanings of words (through models like \href{https://en.wikipedia.org/wiki/Word2vec}{Word2vec}, for instance\footnote{In the literature, these models are often called \textit{distributional models of meanings.}}) and computers can understand grammar (through \href{https://en.wikipedia.org/wiki/Part-of-speech_tagging}{parts of speech tagging}, for instance\footnote{These are often called \textit{symbolic} or \textit{compositional models of meaning.}}). But what's not-so-nice is that computers aren't too good at understanding meanings of \textit{sentences} and longer bodies of text.

In 2010, Bob Coecke, Mehrnoosh Sadrzadeh, and Stephen Clark sought to address this problem in \href{https://arxiv.org/abs/1003.4394}{``Mathematical Foundations for a Compositional Distributional Model of Meaning.''} In this paper, the authors rely heavily on the principle of compositionality---the idea that the meaning of a sentence can be determined by the meanings of its individual words together with the grammatical rules for combining them. So if a computer can understand meanings of individual words and if it can understand grammatical rules, then then only thing it needs help with is knowing how to combine them to form a meaningful whole. And that's where the category theory comes in! Guided by functorial semantics, Coecke et. al. model natural language as a (monoidal) functor between \textit{compact closed categories}
\[\mathsf{grammar}\to\mathsf{meanings\;of\;words}\]
This functor assigns a grammar type to a word, and the monoidal structures provide a way to combine the meanings of those words (and their grammar types) to form a sentence, whose meaning can be determined via the principal of compositionality.

\newthought{in the remaining pages,} we'll dive into the details by answering the following questions:
	\begin{enumerate}[i.]
	\item \textbf{(the syntax category)} How can we make sense of grammar, mathematically? Specifically, how does $\mathsf{grammar}$ form a compact closed category?
	\item \textbf{(the semantics category)} How can we make sense of meanings of words, mathematically? That is, how do $\mathsf{meanings\;of\;words}$ form a compact closed category?
	\item \textbf{(the functor)} How is the functor $\mathsf{grammar}\to\mathsf{meanings\;of\;words}$ defined and how does it allow one to determine the meaning of a full sentence?
	\end{enumerate}

\noindent Let's start by answering the first question.

\newpage
%%%%%%%%%%%%%%%%%
\subsection*{The Syntax Category: Pregroup}
Following the work of Coecke et. al., we can model grammar algebraically via a \textit{pregroup}, a construction due to
	\marginnote[-1cm]{If you were to develop your own categorical-compositional-distributional model (often called \textbf{DisCoCat} models) of meaning, you might wish to work with a construction \textit{other} than pregroups. And that's fine. But the nice thing about a pregroup is that (as we'll soon see) it has the exact same categorical structure as our semantics category, namely compact closure. That means we consider a functor from a pregroup into the semantics category that \textit{preserves} the compact-closed structure. More generally, then, you might define a DisCoCat-type language model to be any monoidal functor $\cat{C}\to\cat{meanings\;of\;words}$ where $\cat{C}$ is any compact closed category accounting for grammar.}
mathematician Joachim Lambek in the early 1990s. Informally, a pregroup is cooked up from the following recipe: 
\begin{center}
\textbf{poset} + \textbf{monoid} + \textbf{``duals''} = \textbf{pregroup}
\vspace{0.3cm}
\includegraphics[width=!,totalheight=!,scale=0.5]{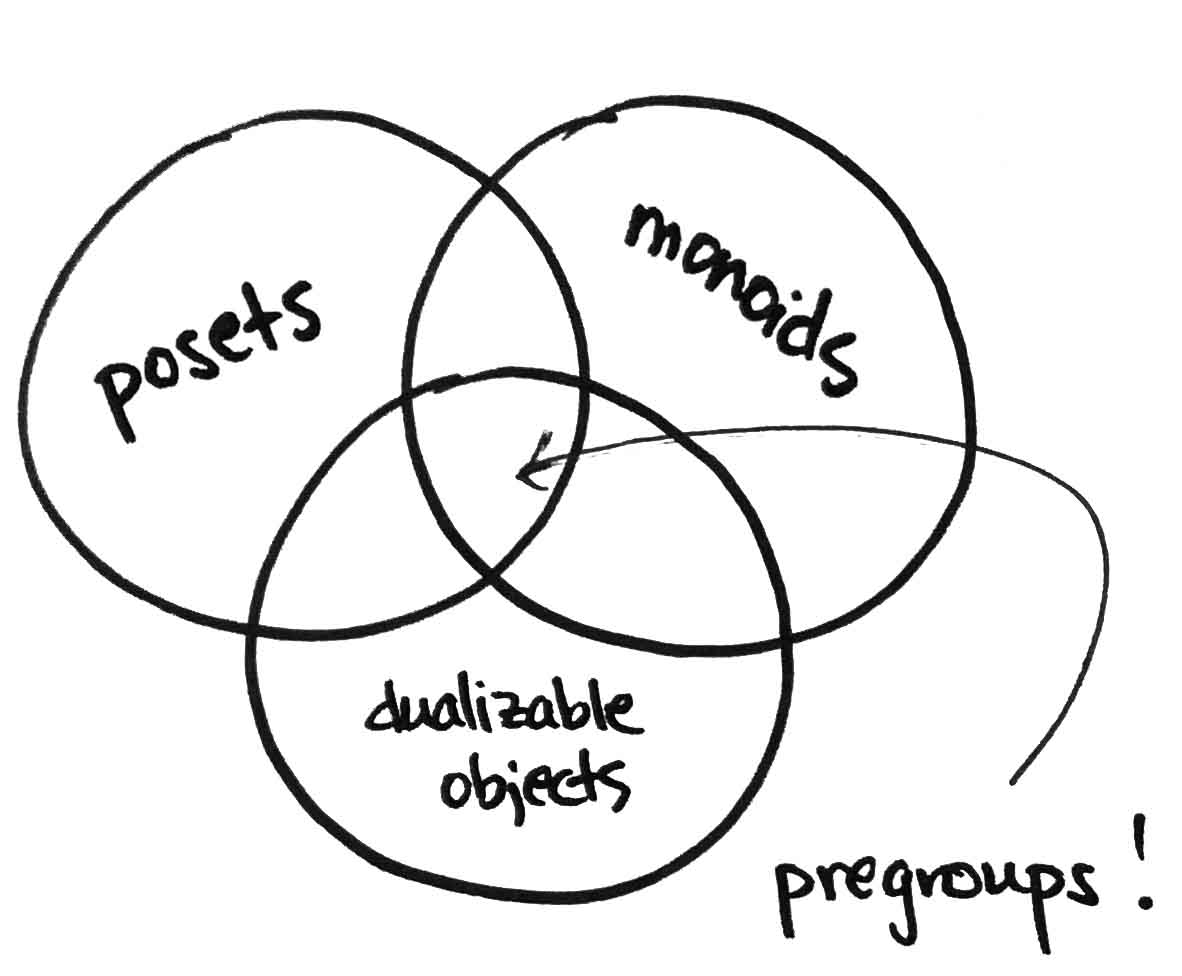}
\end{center}
In other words, a pregroup is 
	\begin{itemize}
	\item a poset $(P,\leq)$ 
	\item that has a multiplication (we'll denote it by juxtaposition) that's compatible with the partial order, i.e. if $p\leq q$ then $ap\leq aq$ and $pa\leq qa$ for all $a\in P$,
	\item together with a unit 1 satisfying $1p=p1=p$ for all $p\in P$,
	\item and moreover each element $p$ has both a left dual $p^l$ and a right dual $p^r$ with maps
	\[p^lp\overset{\epsilon^l}{\leq} 1 \overset{\eta^l}{\leq} pp^l 
	\qquad\text{and}\qquad 
	pp^r \overset{\epsilon^r}{\leq} 1 \overset{\eta^r}{\leq} p^rp.\]
that are required to satisfy the yanking equations in (\ref{eq:yank}).
	\end{itemize}
I've referred to the inequalities as \textit{maps} (and have labeled them as such) because they are actually \textit{morphisms} in a category! Indeed, every poset $(P,\leq)$ can be viewed as a category: an object is an element in $P$ and there is an arrow $p\to q$ if and only if $p\leq q.$ (In particular, there is \textit{at most} one arrow between any two elements in a poset.) Composition is given by transitivity: if $p\leq q$ and $q\leq r$ then $p\leq r$, and associativity is immediate. Also, every element has an identity arrow since the partial order is reflexive: $p\leq p$ for all $p\in P$.
So because a pregroup is a poset, we may also view it as a category. Therefore I'll now draw an arrow $\to$ in lieu of the partial order $\leq$. Moreover, a pregroup is really a poset-with-extra-structure and therefore we may view it as a category-with-extra-structure. Not surprisingly, given the reappearance of the yanking equations, that extra structure is compact closure! In summary, a pregroup is an example of a \textbf{compact closed category}. It is, in particular, a compact closed category that is \textit{not symmetric.} Indeed, those four inequalities above are really the unit and counit maps 
	\[p^lp\overset{\epsilon^l}{\longrightarrow} 1 \overset{\eta^l}{\longrightarrow} pp^l 
	\qquad\text{and}\qquad 
	pp^r \overset{\epsilon^r}{\longrightarrow} 1 \overset{\eta^r}{\longrightarrow} p^rp.\]
discussed in Section \ref{sec:Monoidal}, and the yanking equations (\ref{eq:yank}) amount to the following:
	\begin{align*}
	p&=p\cdot1\overset{1\cdot \eta^r}{\longrightarrow} pp^rp \overset{\epsilon^r\cdot 1}{\longrightarrow} 1\cdot p = p\\[10pt]
	p &= 1\cdot p \overset{\eta^l\cdot 1}{\longrightarrow} pp^lp \overset{1\cdot \eta^l}{\longrightarrow} p \cdot 1=p\\[10pt]
	p &= 1\cdot p^r \overset{\eta^r\cdot 1}{\longrightarrow} p^rpp^r\overset{1\cdot \epsilon^r}{\longrightarrow} p^r\cdot 1=p^r\\[10pt]
	p &= p^l\cdot 1 \overset{1 \cdot \eta^l}{\longrightarrow} p^lpp^l \overset{\epsilon^l\cdot 1}{\longrightarrow} 1\cdot p^l=p^l
	\end{align*}
		\marginnote[-2cm]{In first equality of the third line we're rewriting $p^r$ as $1\cdot p^r$ rather than $p^r\cdot 1$ because neither of the $\eta$s nor $\epsilon$s provide a map $1\to pp^r$. Similarly for the last line, write $p^l=p^l\cdot 1$ rather than $p^l=1\cdot p^l$ since there's no map $1\to p^lp$.}

\newthought{Let's look at two} examples of pregroups. The first is an arithmetic example, which will help to get our feet wet. The second is a grammatical example, which is used in the DisCoCat model of Coecke, Sadrzadeh, and Clark.

\begin{ex} The set \[\{f\colon \mathbb{Z}\to\mathbb{Z}\;\vert\; f\text{ is monotone and unbounded}\}\]
is a pregroup. The partial order is given pointwise: $f\leq g$ if and only if $fn\leq gn$ for all $n.$ The monoid multiplication is given by function composition $f\cdot g:=f\circ g.$ The monoidal unit is $\text{id}_\mathbb{Z}.$ Given such a function $f,$ its left and right duals are given by
\[f^ln:=\min\{m\in\mathbb{Z}\;\vert\; n\leq fm\} \qquad\qquad f^rn=\max\{m\in\mathbb{Z}\;\vert\;fm\leq n\}.\]
For example, if $fm=2m,$ then
	\marginnote[-3cm]{\textbf{Fun fact}: the pair $(f^l,f)$ forms a special kind of categorical adjunction called a \textbf{Galois connection} since it satisfies \[f^ln\leq m \quad \Leftrightarrow \quad n\leq fm.\] Indeed if $n$ is even, then $n/2\leq m \Leftrightarrow n\leq 2m.$ And if $n$ is odd, then $(n+1)/2\leq m$ which means $n+1\leq 2m$ which is true iff $n\leq 2m.$ Similarly, the pair $(f,f^r)$ forms a Galois connection since \[fn\leq m \quad\Leftrightarrow n\leq f^rm.\] Indeed, if $n$ is even then $2n\leq m\Leftrightarrow n\leq m/2.$ And if $n$ is odd, then $2n\leq m$ means $n\leq m/2$ which is true iff $n\leq (m-1)/2.$ 

	For a couple of great introductions to \textit{Galois connections} (They are super cool and appear in lots of places in math!) take a look at \href{https://forum.azimuthproject.org/discussion/1828/lecture-4-chapter-1-galois-connections}{Lecture 4} of John Baez's \href{https://forum.azimuthproject.org/categories/applied-category-theory-course}{online course on applied category theory} as well as Section 1.5 of \href{https://arxiv.org/pdf/1803.05316.pdf}{\textit{Seven Sketches}} by Fong and Spivak.}
\[
f^ln=\begin{cases}
\frac{n}{2} &\text{if $n$ is even},\\
\frac{n+1}{2} &\text{if $n$ is odd}
\end{cases}
\qquad \qquad 
f^r=
\begin{cases}
\frac{n}{2} &\text{if $n$ is even,}\\
\frac{n-1}{2} &\text{if $n$ is odd.}
\end{cases}
\]
In short, $f^ln=\lfloor \frac{n+1}{2} \rfloor$ and $f^rn=\lfloor \frac{n}{2}\rfloor$.

You can find this example in ``Iterated Galois Connections in Arithmetic and Linguistics'' by Lambek, which appears in the Springer book \href{https://www.springer.com/us/book/9781402018978}{\textit{Galois Connections and Applications}}. You'll also find mention of it in the \href{https://arxiv.org/abs/1003.4394}{``Mathematical Foundations''} paper of Coecke et. al.
\end{ex}

While arithmetic is fun, this next example is the one we're most interested in. 

\begin{ex}
Given any finite poset $X$, we can construct the \textit{free pregroup generated by $X,$} denoted $\cat{Preg}X$. For a simple example, suppose $X=\{n,s\}$ whose elements we'll think of as basic grammar types: $n$ is the type of a \textit{noun} and $s$ is the type of a (declarative) \textit{sentence}. Elements of $\cat{Preg}\{n,s\}$ are concatenations of the letters $n$ and $s$ and their left and right duals and iterations of those duals and so on. For example, some grammatical types in $\cat{Preg}\{n,s\}$ are:
\begin{center}
\includegraphics[width=!,totalheight=!,scale=4]{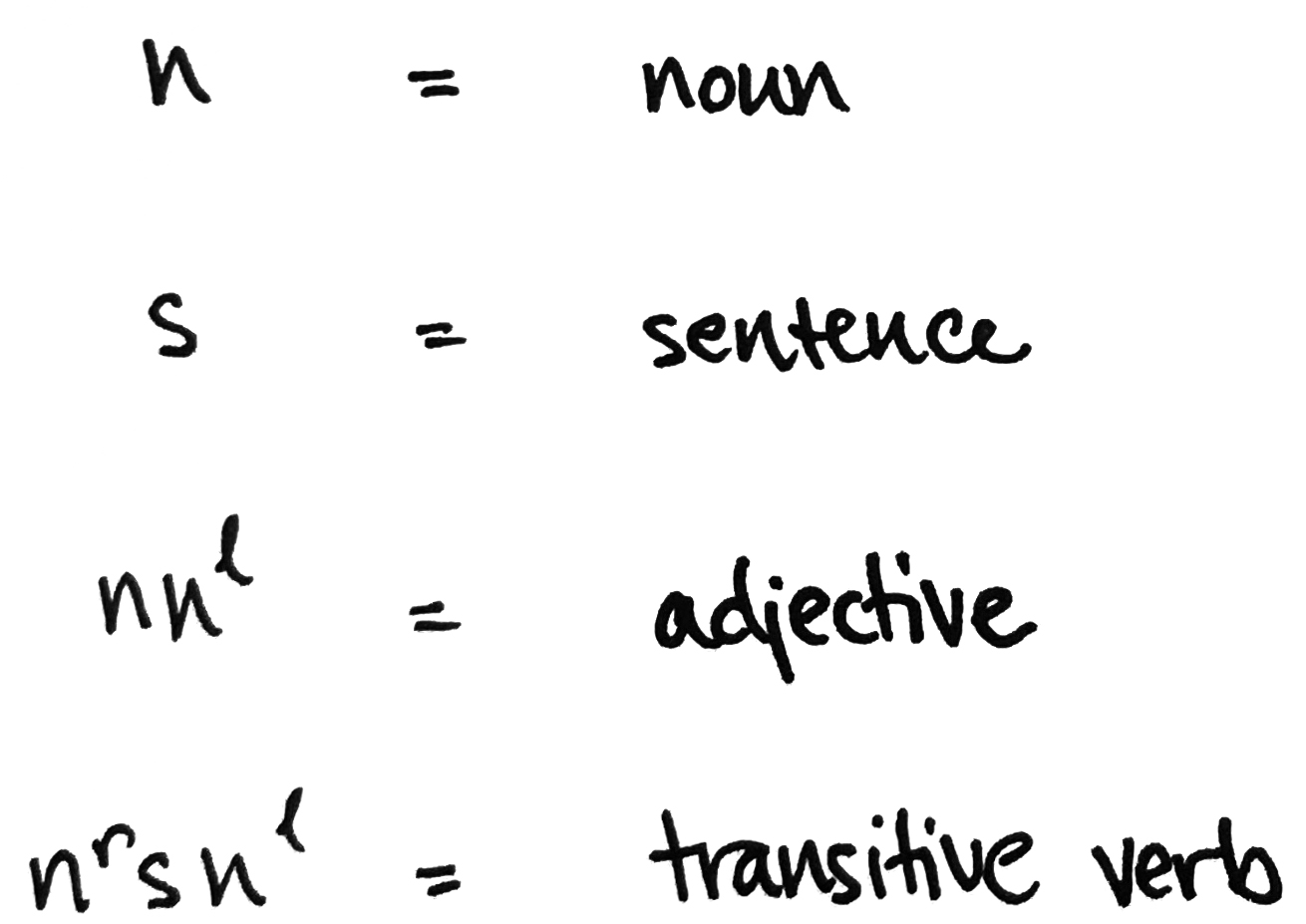}
\end{center}
The strings of letters are called \textit{compound types}, and there is a morphism $a\to b$ between compound types if and only if $a$ can \textit{reduce} to $b$ by application of one or more of the counit maps $\epsilon^r$ and $\epsilon ^l$.

Consider a banana, for example. It has type $n$, of course, while the adjective \textcolor{Goldenrod}{yellow} has type $nn^l$. The reason that adjectives have grammar type $nn^l$ is that an adjective can always be paired on the left with a noun, resulting in a new noun---e.g. \textit{yellow banana.} 
\begin{center}
$\underset{nn^l}{\text{\large yellow}} \quad 
\underset{n^{\phantom{l}}}{\text{\large banana}}$
\end{center}
Indeed, to verify that the grammar type of \textit{yellow banana} is $n$, we start by concatenating the types of the individual words to obtain $nn^ln$. Then we apply the counit map $\epsilon^l\colon n^ln\to~1$ together with the identity map $1_n\colon n\to n$ (this is given to us by the reflexivity axiom of posets: $n\leq n$) to see that $nn^ln$ reduces down to $n$:
	\marginnote[-1.5cm]{The dot $\cdot$ in $1_n\cdot \epsilon^l$ is meant to suggest ``apply $1_n$ to $n$ while simultaneously applying $\epsilon^l$ to $n^ln$.''}
\begin{fullwidth}
\begin{center}
\includegraphics[width=!,totalheight=!,scale=0.5]{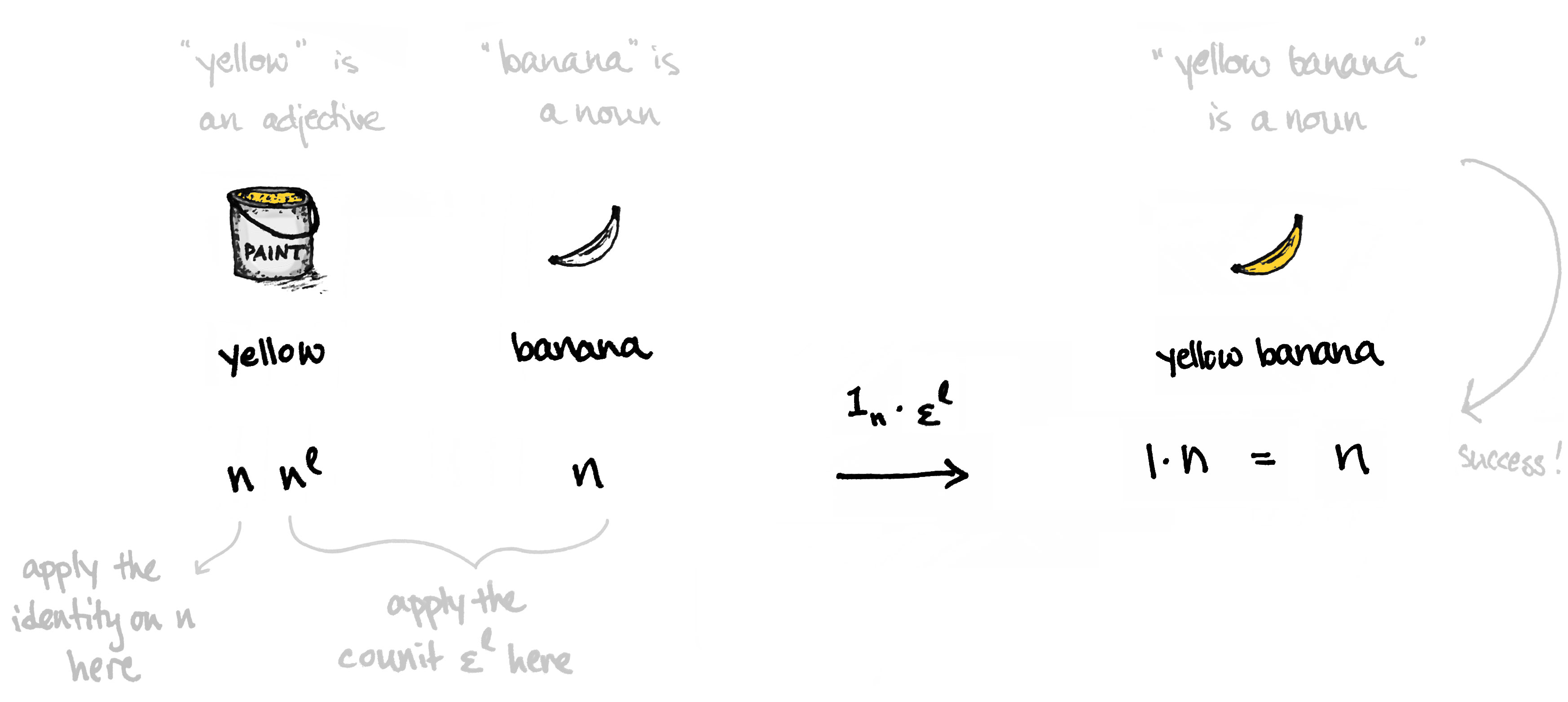}
\end{center}
\end{fullwidth}
This tells us that the phrase \textcolor{Goldenrod}{yellow banana} has grammar type $n$. That's good. A yellow banana is a noun!  

In light of this discussion on yellow bananas, you might enjoy taking a few seconds to think about why $n^rsn^l$ represents the grammar type of \textit{transitive verb}.

\begin{center}
\noindent \texttt{<ponder> ... </ponder>} 
\end{center}

\noindent A transitive verb is a word that accepts a noun on the right and another noun on the left such that the resulting phrase is a full sentence. Since we like bananas, here's another fruit-based example:
% \begin{center}
% bananas $\rightsquigarrow$ $n$\\
% are $\rightsquigarrow$ $n^rsn^l$\\
% fruit $\rightsquigarrow$ $n$ 
% \end{center}
\begin{center}
$\underset{n^{\phantom{l}}}{\text{\large bananas}} \quad 
\underset{n^rsn^l}{\text{\large are}} \quad
\underset{n^{\phantom{l}}}{\text{\large fruit}}$
\end{center}
To determine the grammar type of this phrase, we concatenate the grammatical types of the individual words and then apply the counit maps to reduce, as before:
\begin{fullwidth}
\begin{center}\label{pic:fruit}
\includegraphics{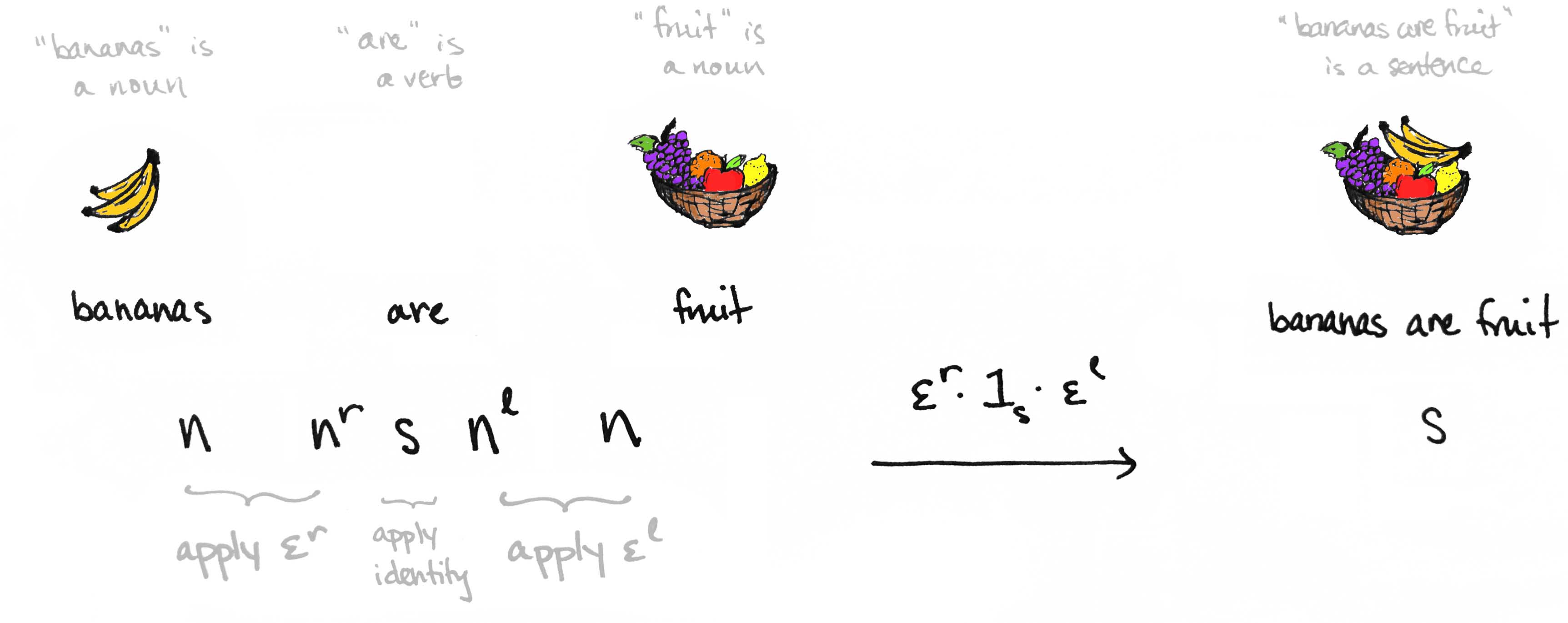}
\end{center}
\end{fullwidth}
Here's that same reduction written out step-by-step. For clarity, I'll indicate the concatenation with a dot: 
	\[
	n n^rsn^ln= nn^r\cdot s\cdot n^ln \overset{\epsilon^r\cdot 1_s\cdot 1_n}{\longrightarrow} 1\cdot s\cdot n^ln = s\cdot n^ln 
	\overset{1_s\cdot \epsilon^l}{\longrightarrow} s\cdot 1
	=s
	\]
In words, we've used the counit maps to reduce the concatenation of the grammar types for \textit{bananas are fruit} to the letter $s$, which confirms that ``bananas are fruit'' is indeed a grammatically correct sentence. More generally, for any $a,b\in \cat{Preg}\{n,s\}$, we draw a morphism $a\to b$ if and only if $a$ can be reduced to $b$ in a similar fashion.
\end{ex}

\newthought{Alright, that's} (a very condensed version of) the pregroup story! To summarize, the language model of Coecke et. al. amounts to a structure-preserving functor
\[\cat{syntax}\to\cat{semantics}\]
In this section, we've just shown that the syntax category is taken to be a pregroup freely generated on a finite set of basic grammar types, i.e. $\cat{syntax} = \cat{Preg}X.$ Let's move on to semantics now.

%%%%%%%%%%%%% SECTION 4 %%%%%%%%%%%%%%
\subsection*{The Semantics Category: Vector Spaces}
As stated in the introduction to Section \ref{sec:Second}, computers are able to understand meanings of individual words pretty well. That's because computers understand numbers! For example, a great way to inform a computer of the meaning of the word \textit{banana}  is to represent \textit{banana} by a number and then give that number to the computer.

\vspace{0.4cm}
\noindent \textit{What number?}
\vspace{0.4cm}

\noindent Well, it's not exactly a number. It's an array of numbers---a vector.

\vspace{0.4cm}
\noindent \textit{Okay, what vector?}
\vspace{0.4cm}

\noindent The answer is simple, though I'd like to motivate it by sharing the following theorem.

\vspace{0.4cm}
\noindent\textbf{Theorem} (The Yoneda Lemma for Linguistics). \textit{You shall know a word by the company it keeps.}
\begin{proof}
John Firth\footnote[][-1cm]{Firth, J. R. A synopsis of linguistic theory, 1930--1955. In \textit{Selected Papers of JR Firth}, 1952--59 (ed. J. Firth and F. Palmer). Indiana University Press.}
\end{proof}
% \begin{thm}[\textbf{Yoneda's Lemma for Linguistics}] \textit{You shall know a word by the company it keeps.} -- John Firth
% \end{thm}

Okay, so it's not a theorem. But it is a great quote! Firth's idea is that \textit{words that appear in similar contexts will have similar meaning.} \marginnote{Are you wondering why I've referred to Firth's idea as the Yoneda Lemma? To find out why, I recommend reading up on the \href{http://www.math3ma.com/mathema/2017/8/30/the-yoneda-perspective}{Yoneda Perspective}.} In the linguistics community, this is referred to as the \href{https://en.wikipedia.org/wiki/Distributional_semantics#Distributional_hypothesis}{distributional hypothesis}. So you might imagine that \textit{apple} is more similar to \textit{banana} than it is to \textit{puppy} since \textit{apples} and \textit{bananas} often occur near words such as sweet, snack, green, eat, etc., whereas \text{puppy} occurs more often near words such as pet, cute, furry, bark, and so on. As another example, you might not know what the word \textit{yegg} means (or perhaps you do), but you can probably infer it from \href{https://www.merriam-webster.com/word-of-the-day/yegg-2018-03-20}{this sentence}:
\begin{center}
\small
The cops grabbed him and another yegg for a Philadelphia store burglary.\footnote{James Lardner and Thomas Reppetto, \textit{NYPD: A City and Its Police,} 2000}
\end{center}

\newthought{So we can represent} the meaning of a word by a vector. This is often called a \textit{distributional model of meaning.} But what, exactly, \textit{is} the assignment \textit{word} $\mapsto$ \textit{vector}? Suppose we have a fixed corpus---your favorite book, say. Start by choosing a set of so-called \textbf{context words} $\{w_1,\ldots, w_n\}$. This can be every word in the corpus or some subset of it. By representing each $w_i$ as the $i$th standard basis vector
\[\mathbf{w}_i=(0,\ldots,\overbrace{1}^{i^\text{th} \text{spot}},\ldots,0)\]
we obtain a basis $\{\mathbf{w_1},\ldots, \mathbf{w_n}\}$ for a vector space $V$. Then any word $w$ in the corpus has a vector representation given by a linear combination of the context words
\[\mathbf{w}=\sum_{i=1}^n c_i\mathbf{w_i}\]
The coefficients $c_i$ are real numbers that indicate the number of times that $w$ occurs near\footnote{You can decide what ``near'' means. That is, the \textit{context of $w$} is the set of words within $k$ words of $w$, where $k=1$ or 2 or 3 or whatever you like.} $w_i$ in the corpus. 

Here's an example. Suppose we're reading a book that contains the words
\[\{ \textit{sweet},\textit{green},  \textit{furry} \}\]
Let's choose them to be our context words and make the assignment
so that \[
\mathbf{sweet} = \begin{bmatrix}1\\0\\0\end{bmatrix} \qquad 
\mathbf{green} = \begin{bmatrix}0\\1\\0\end{bmatrix}\qquad 
\mathbf{furry}=\begin{bmatrix}0\\0\\1\end{bmatrix}\] 
Then if \textit{banana, puppy} and \textit{fruit} are also words in our book, we might have something like
\[\mathbf{banana} = \begin{bmatrix}21\\9\\0\end{bmatrix}\qquad 
\mathbf{puppy} = \begin{bmatrix}8\\1\\32\end{bmatrix} \qquad 
\mathbf{fruit} = \begin{bmatrix}43\\19\\0\end{bmatrix}\]
		\marginnote[-2cm]{\includegraphics[width=!,totalheight=!,scale=0.6]{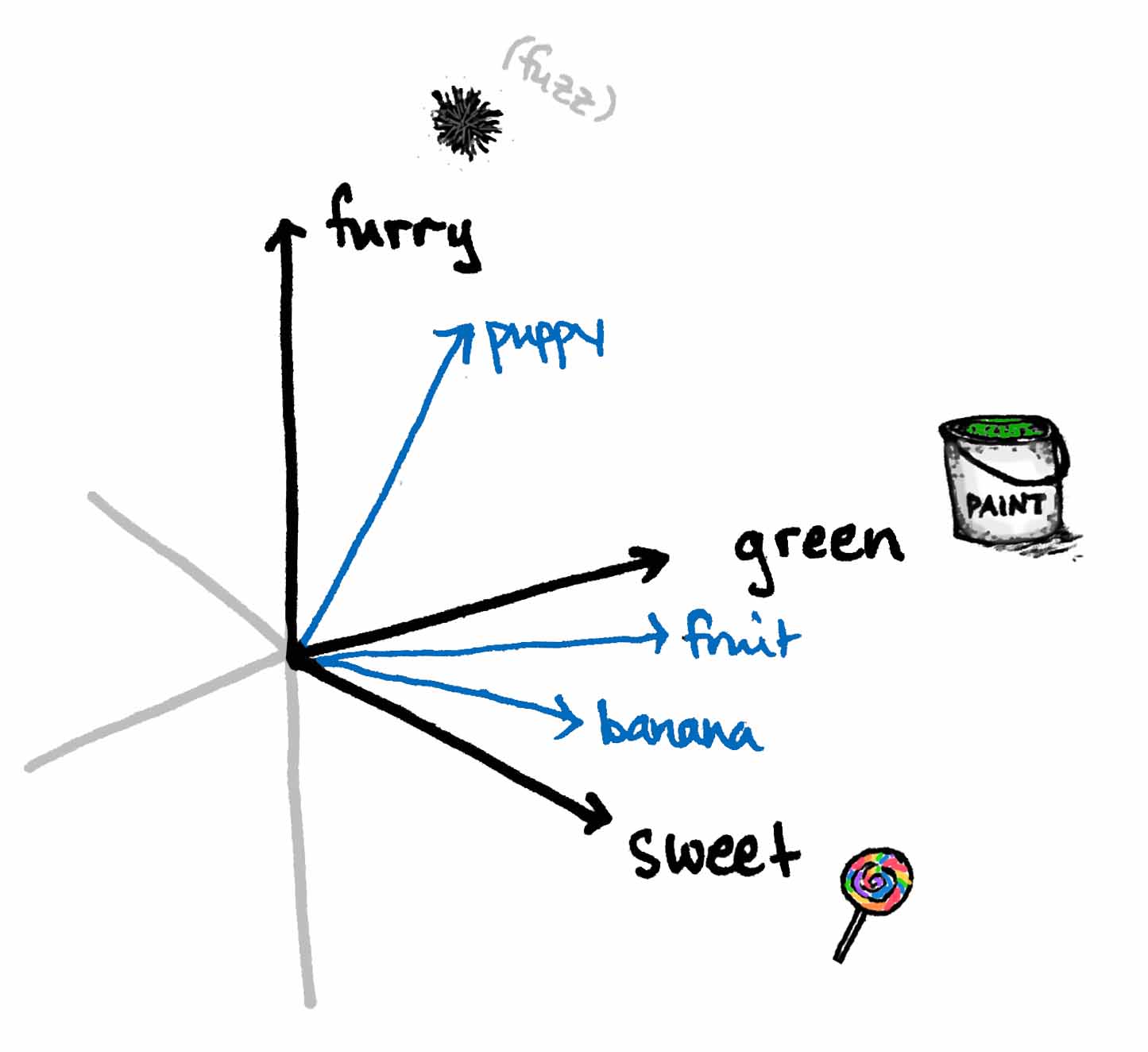}}
In other words, we've used data from the corpus to \textit{embed} these words as vectors inside of a three-dimensional vector space. This prompts us to say that the \textit{meaning} of the word banana \textit{is} the vector $(21,9,0)$, the meaning of puppy is $(8,1,32)$, and the meaning of fruit is $(43,19,0)$.

And this works! That is, you can feed distributional models into your computer, and they'll ace the word-similarity portion of your SAT exam. Or your can compute the dot product between words, and you'll find that vectors are closer together precisely when the words they represent have the same meaning. It's all familiar territory for NLP practitioners. The semantics category for Coecke et. al. is thus the category of finite dimensional vector spaces over $\mathbb{R}$. That is,
\[\cat{meanings\;of\;words} = \cat{FVect}\]

\newthought{Unfortunately, the distributional model} does \textit{not} work for sentences. The same sentence rarely occurs twice in a given document, therefore we can't follow the same procedure above. This is where category theory can help lend a hand. In light of the principle of compositionality, the meaning of a sentence should be able to be computed given the meanings of its individual words and the rules of grammar for combining them. And we can pair meanings-of-words with grammatical types via a map from syntax ($\cat{grammar}$) to semantics ($\cat{meanings\;of\;words}$), i.e. via a functor
\[F\colon \cat{Preg}X\to\cat{FVect}\]
where $X$ is a finite set of basic grammar types. In fact, as we know from Section \ref{sec:Monoidal}, both $\cat{FVect}$ and $\cat{Preg}X$ are compact closed categories, so we'll ask that $F$ be a \textit{strong monoidal functor}---one that preserves the compact closed structure. That's the gist behind the categorical compositional distributional model of Coecke et. al.

\vspace{0.4cm}
\noindent But how is $F$ actually \textit{defined}? Let's talk about that next.

% \newthought{Moreover, as we noted before,} the category $\cat{Vect}$ of finite dimensional vector spaces is compact closed! Here, the dual of a vector space $V$ is its \textit{literal} dual $V^*=\hom(V,\mathbb{R}),$ which in fact is isomorphic to $V$. So instead of writing $V^*$, I'll just write $V.$ The monoidal unit is $\mathbb{R}$, and the unit and counit maps $\eta$ and $\epsilon$ are given \textcolor{red}{in an earlier section.}

%%%%%%%%%%%%% SECTION 6 %%%%%%%%%%%%%%
\subsection*{The Functor: Syntax $\to$ Semantics} 
In this section, we'll give an explicit description of the functor
\[F\colon \cat{syntax}\to\cat{semantics}\]
or more specifically,
\[F\colon \cat{Preg}X\to\cat{FVect}\]
For simplicity, let's take $X=\{n,s\}$ as we did before. Now to define a functor, we need simply to say what it does on objects and morphisms. So let's do that. On objects, 
		\begin{itemize}
			\item $F$ assigns to the noun type $n$ a vector space $N:=Fn$, which we'll call a \textit{noun space}
			\item $F$ assigns to the sentence type $s$ a vector space $S:=Fs$, which we'll call a \textit{sentence space}
		\end{itemize}
and on morphisms
		\begin{itemize}
			\item $F$ assigns to a type reduction $a\overset{r}{\longrightarrow} b$ a linear map $Fa\overset{Fr}{\longrightarrow}Fb$ that sends the vector corresponding to a word or phrase of type $a$ in $Fa$ to the vector corresponding to a word or phrase of type $b$ in $Fb$.
		\end{itemize}
Moreover, asking that $F$ preserve the compact closed structure means that
	\begin{itemize}
		\item units and counits in $\cat{Preg}\{n,s\}$ map to units and counits in $\cat{FVect}$
			
			\begin{quote}
				\underline{e.g.} given $n\in \cat{Preg}\{n,s\}$, 
				\[F(1\overset{\eta_n^r}{\longrightarrow} n^rn)=F(1\overset{\eta_n^l}{\longrightarrow} nn^l)=\eta_N 
				\qquad\text{and}\qquad
				F(nn^r\overset{\epsilon_n^r}{\longrightarrow}1)=F(n^ln\overset{\epsilon_n^l}{\longrightarrow}1)=\epsilon_N 
				\]
				where $\eta_N\colon \mathbb{R}\to N\otimes N$ and $\epsilon_N\colon N\otimes N\to \mathbb{R}$ are the linear maps that we defined on p. \pageref{lis:ucounit}. A similar idea holds if we replace $n$ by any element of $\cat{Preg}\{n,s\}$.
			\end{quote}
			
		\item duals map to duals
			
			\begin{quote} 
				\underline{e.g.} $Fn^r=Fn^l=N^*$. But our vector spaces are finite dimensional and so $N^*\cong N$ and therefore $Fn^r=Fn^l=N$.
			\end{quote}
			
		\item a compound type is assigned to a tensor product of vector spaces. 
		
			\begin{quote}
				\underline{e.g.} $F(n^rsn^l)\cong Fn^r\otimes Fs\otimes Fn^l=N\otimes S\otimes N$
			\end{quote}
			
	\end{itemize}

\vspace{0.4cm}
\noindent \textit{And that's it!}
\vspace{0.4cm}

Except... this might not be very enlightening yet. It'll surely be helpful to look at a toy example. So in the next couple of pages, let's use the DisCoCat model to compute the meaning of the sentence

\begin{center}
\textit{\large bananas are fruit}
\end{center}

\noindent By ``compute the meaning,'' I mean the following: we want to be able to view the sentence \textit{bananas are fruit} as a vector, then feed that vector into the functor $F$ and get an output vector that encodes for the meaning of the sentence. \[F(bananas\;are\;fruit)=\;??\] That output vector will be the ``meaning'' of the sentence. Our goal is to find that meaning.

\vspace{0.4cm}
\begin{center}
\noindent\textbf{\large Goal: Compute the meaning of \textit{bananas are fruit}}.
\end{center}

Let's proceed systematically. \marginnote{\begin{center}\includegraphics{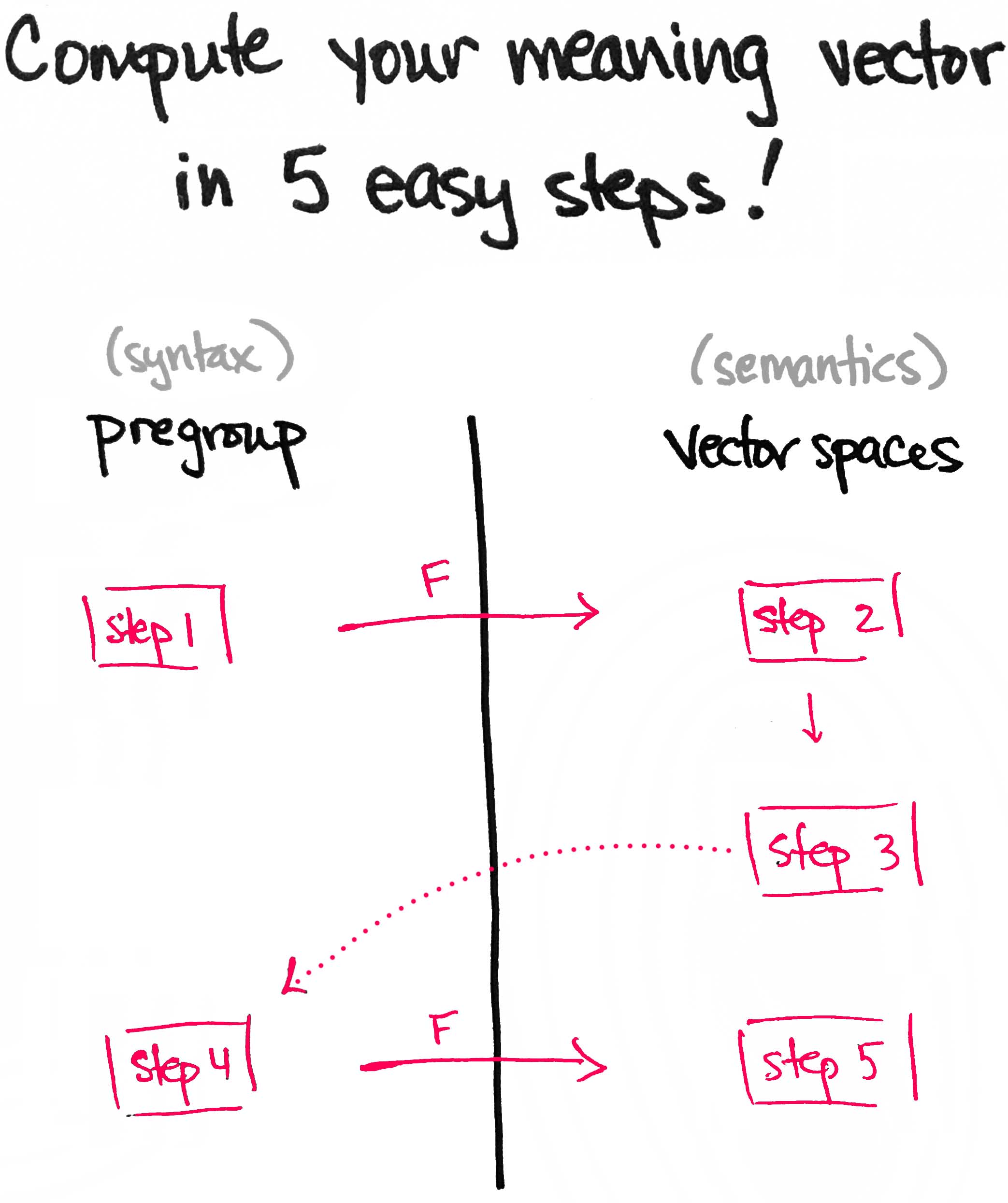}\end{center}} I'll list the computations step-by-step, starting from the beginning.

\vspace{0.4cm}
\noindent\textbf{\large Step 1: Assign each word a grammar type in $\cat{Preg}\{n,s\}$.}\\[10pt]

\noindent That's easy enough:\\[0.3cm]
bananas \textcolor{gray}{$\rightsquigarrow$} $n$\\
are \textcolor{gray}{$\rightsquigarrow$} $n^rsn^l$\\
fruit \textcolor{gray}{$\rightsquigarrow$} $n$

\vspace{0.4cm}
\noindent \textbf{\large Step 2: Fix a noun space $Fn=N$ and a sentence space $Fs=S$.}\\[10pt]

\noindent Let's suppose $N$ is the three-dimensional space spanned by the basis vectors \[\mathbf{sweet}, \mathbf{green}, \mathbf{furry}\] which we can represent as column vectors
\begin{align*}
\mathbf{w_1}=\mathbf{sweet} = \begin{bmatrix}1\\0\\0\end{bmatrix},\qquad
\mathbf{w_2}=\mathbf{green} = \begin{bmatrix}0\\1\\0\end{bmatrix},\qquad
\mathbf{w_3}=\mathbf{furry} = \begin{bmatrix}0\\0\\1\end{bmatrix}
\end{align*}
as before. These basis vectors generate the noun space. But what about the sentence space $S$? For simplicity, let's define $S$ to be a ``true or false'' space so that it's a one-dimensional vector space spanned by a single vector $\vec{1}$. The origin $0\in S$ corresponds to ``false'' while  $\vec{1}$ corresponds to ``true.'' What about scalar multiplies of $\vec{1}$? If you like, you're more than welcome to think of a positive scalar multiple of $\vec{1}$ as the meaning vector for sentence that is \textit{super} true. The larger the scalar, the more true the sentence! 

Finally, note that once we've established $N$ and $S$, the verb space comes for free:
\[F(n^rsn^l)=N\otimes S\otimes N\]
This is a nine-dimensional space spanned by vectors of the form $\mathbf{w_i}\otimes\vec{1}\otimes\mathbf{w_j}$ where $i$ and $j$ range between 1 and 3.

\vspace{0.4cm}
\noindent\textbf{\large Step 3: Determine the vector representations of each word in the sentence}.\\[10pt]

\noindent We'll simply recycle the vectors we used earlier:
\begin{align*}
\mathbf{bananas} = \begin{bmatrix}21\\9\\0\end{bmatrix},\qquad\qquad
\mathbf{fruit} = \begin{bmatrix}43\\19\\0\end{bmatrix}
\end{align*}
Note that both of these vectors live in the noun space $N$ since each word has grammar type $n$. But what about the transitive verb \textit{are}? By Step 1, we know that $are$ has grammar type $n^rsn^l$ and is therefore a vector in the tensor product $N\otimes S\otimes N$. That is, there are coefficients $c_{ij}\in\mathbb{R}$ so that
\[are = c_{11}\mathbf{sweet}\otimes \vec{1}\otimes \mathbf{sweet} \;\;+\;\; 
c_{12}\mathbf{sweet}\otimes \vec{1}\otimes \mathbf{green} 
\;\;+\;\; \cdots \;\;+\;\; c_{33}\mathbf{furry}\otimes \vec{1}\otimes \mathbf{furry}.\]

\vspace{0.4cm}
\noindent Eek. That looks uncomely. 
\vspace{0.4cm}

\noindent Fortunately, we learned in Section \ref{sec:Monoidal} that $\cat{FVect}$ is a compact closed category and therefore it exhibits process-state duality, which is the sophisticated way of saying 
\begin{center}\textit{every vector in a tensor product can be identified with a linear map},
\end{center} 
which is the long way of saying 
\begin{center}\textit{every vector is really a matrix!}
\end{center} 
\noindent And that is excellent news, for if we know ``what is what,'' i.e. if we know that $puppies$ are $furry$ but not $green$ and so on, then we can re-express the vector for $are$ as a $3\times 3$ matrix. The $ij$th entry of this matrix is the coefficient $c_{ij}$ which is
\marginnote[2cm]{It's no surprise that we get the identity matrix. \textit{Being} is all about identity. That is, the verb \textit{are} tells you when something IS something else.}
\[
c_{ij}=
\begin{cases} 1, &\text{if $w_i$ is $w_j$},\\
0, &\text{otherwise}
\end{cases}
\]
The upshot is that the transitive verb $are$ has matrix representation
\begin{align*}\label{eq:are}
\mathbf{are} = 
\begin{bmatrix}
1 & 0 & 0\\
0 & 1 & 0\\
0 & 0 & 1
\end{bmatrix}
\end{align*}

\vspace{0.4cm}

\noindent\textbf{\large Step 4: Choose a type reduction in $\cat{Preg}\{n,s\}$}\\[10pt]

\noindent In this step, which takes place in the grammar category, we simply perform the type reduction already done on page \pageref{pic:fruit}. To recap, we know the grammar types of $bananas$ and $are$ and $fruit,$ and so we concatenate those types to obtain $nn^rsn^ln$. Using the left and right counit maps, this string of letters reduces down to $s$, which confirms that the phrase $bananas\;are\;fruit$ is a tried-and-true sentence. In Step 4, we simply take that reduction morphism
\[nn^rsn^ln\overset{\epsilon_n^r\cdot 1_s\cdot \epsilon_n^l}{\longrightarrow} n\]
and hold on to it. We'll need to use it in Step 5.
\vspace{0.4cm}
	\begin{quote}
	Aisde: You might wonder about the word ``Choose'' in ``Step 4: Choose a type reduction.'' \textit{What's up with that?} Incidentally, no choice was needed in this toy example of ours, so the purpose of this aside might be unclear. Indeed, there's only one way to parse the sentence \textit{bananas are fruit}. But there exist sentences that can be parsed in \textit{more than one} way. Consequently, the grammar type of such sentences may reduce down to type $s$ via \textit{more than one} reduction morphism. In Step 4, we are required to choose one. As an illustration, here is a nice sentence:
	\begin{center}
	\textit{I saw a man with a telescope.}
	\end{center}
	How did you parse it? Perhaps
	\begin{center}
	\textit{I saw (a man with a telescope).}
	\end{center}
	or perhaps
	\begin{center}
	\textit{I saw (a man) with a telescope.}
	\end{center}
	Those are two parsings of the same sentence, each of which corresponds to a different type reduction in the pregroup. In turn, this gives rise to different meaning vectors! And rightly so. Those two sentences have different meanings! Step 4 is simply reminding us of this fact.
	\end{quote} \vspace{0.3cm}
	\begin{center}
	\includegraphics[width=!,totalheight=!,scale=0.15]{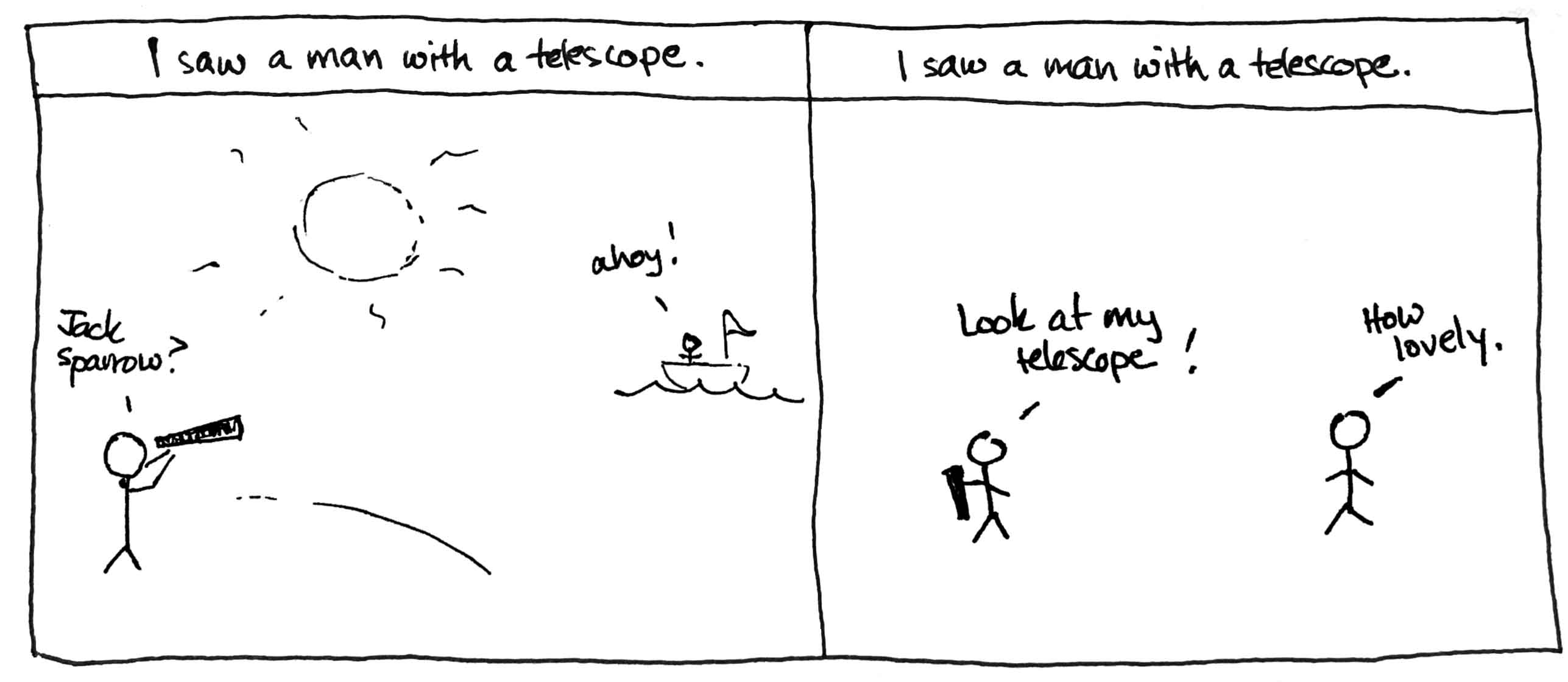}
	\end{center}

\vspace{0.4cm}
\noindent\textbf{\large Step 5: Apply $F$!}\\[10pt]

\noindent This is the fun part! We have a morphism in the pregroup
\[nn^rsn^ln\overset{\epsilon_n^r\cdot 1_s\cdot \epsilon_n^l}{\longrightarrow} n\]
and we can apply $F$ to get a linear map of vector spaces
\[N\otimes N\otimes S\otimes N\otimes N\overset{\epsilon_N\otimes 1_S\otimes \epsilon_N}{\longrightarrow}S\]
where $\epsilon\colon N\otimes N\to \mathbb{R}$ is the linear map given on page \pageref{lis:ucounit} and $1_S\colon S\to S$ denotes the identity map on $S$. Finally, apply this linear map to the vector corresponding to the sentence
\[\epsilon_N\otimes 1_S\otimes \epsilon_N(\mathbf{bananas}\otimes \mathbf{are}\otimes \mathbf{fruit})\]
which amounts to a simple matrix multiplication
\[\epsilon_N\otimes 1_S\otimes \epsilon_N(\mathbf{bananas}\otimes \mathbf{are}\otimes \mathbf{fruit})
\;\;=\;\; \begin{bmatrix}21 & 9 & 0
\end{bmatrix}
\begin{bmatrix}
1 & 0 & 0\\
0 & 1 & 0\\
0 & 0 & 1
\end{bmatrix}
\begin{bmatrix}
43\\ 
19\\
0
\end{bmatrix} \;\;=\;\; 1074\]
Conclusion? The meaning of the sentence $bananas\;are\;fruit$ is 
\[1074\vec{1}\]
which is \textit{super true}. Voila!

\newpage
\subsection*{Some Closing Remarks} 

This functor described above is somewhat reminiscent of a \textit{topological quantum field theory}, which is a functor from the category of cobordisms (another compact closed category) to the category of complex Hilbert spaces. But in 2014 Anne Preller showed that  the only functors from a pregroup $P$ freely generated on a finite set of basic types to $\cat{FVect}$ are those mapping to one-dimensional spaces.	The key to her proof is the fact that $P$ is a \textit{poset} and hence there is at most \textit{one} morphism between any two objects. In particular, any morphism from an object to itself must be the identity. As a consequence, if $a$ is in $P$ then the morphism $(1_a\cdot \eta_a)\circ(\epsilon_a\cdot 1_a)$ from $a a^r a\to a a^r a$ must equal the identity on $a a^r a$. Graphically:
	\[
	\begin{tikzcd}
	a \arrow[r, no head, bend right] & a^r & a \arrow[lldd, no head] &  & a \arrow[dd, no head] & a^r \arrow[dd, no head] & a \arrow[dd, no head] \\
 	&  &  & = &  &  &  \\
	a & a^r \arrow[r, no head, bend left] & a &  & a & a^r & a
	\end{tikzcd}
	\]
Now consider a functor $P\to \mathsf{FVect}$. It assigns $a$ in $P$ to a vector space $A$ in $\mathsf{FVect}$, and it  assigns $(1_a\cdot \eta_a)\circ(\epsilon_a\cdot 1_a)$ to the corresponding linear map, $(1_A\cdot \eta_A)\circ(\epsilon_A\cdot 1_A)\colon A\otimes A^*\otimes A\to A\otimes A^*\otimes A$, which we'll just denote by $f$,
	\[f =
	\begin{tikzcd}
	A \arrow[r, no head, bend right] & A^* & A \arrow[lldd, no head] \\
 	&  &  \\
	A & A^* \arrow[r, no head, bend left] & A
	\end{tikzcd}
	\] 
and which must be an isomorphism. Now if the dimension of $A$ is at least 2, then we can choose orthogonal basis vectors $e_1$ and $e_2$ so that $e_1\otimes e_2^*\otimes e_1\in A\otimes A^*\otimes A$. And since $\epsilon_A$ computes the inner product between $e_1$ and $e_2^*,$ we have $f(e_1\otimes e_2^*\otimes e_1)=0.$ Therefore $f$ is not injective, and so it cannot be an isomorphism.

The intuition is, perhaps, that pregroups have \textit{too few} morphisms to capture the semantics. In particular, pregroups do not allow us to distinguish different parsings of strings of types. One string may reduce in several ways---e.g. (Men and women) whom I like vs. Men and (women whom I like)---and the morphisms in a pregroup do not account for this. So in some sense, there isn't enough ``wiggle room'' for meaning in pregroup syntax, so the output can only be a one-dimensional vector space. But all is not lost! As Preller showed, the problem can be fixed by replacing a free pregroup with a free compact closed category. For more details, see her paper \href{https://arxiv.org/pdf/1412.8527.pdf}{``From Logical to Distributional Models.''}

\newpage
\subsection{Further Reading}\label{ssec:furtherreading3}

\textbf{For more on the work of Baez and Pollard}: 
\begin{itemize}
\item Certainly take a look at their paper, \href{https://arxiv.org/abs/1704.02051}{``A Compositional Framework for Reaction Networks''}. The $n$-Category Caf{\'e} also contains expositions, including \href{https://golem.ph.utexas.edu/category/2017/07/a\_compositional\_framework\_for\_2.html}{a post} by Baez with the same title and \href{https://golem.ph.utexas.edu/category/2018/04/dynamical\_systems\_and\_their\_st.html}{``Dynamical Systems and Their Steady States''} by Maru Sarazola.

\item The examples in Section \ref{sec:First} can be found in \href{https://www.youtube.com/watch?v=IyJP_7ucwWo}{``The Mathematics of Networks''}, a wonderfully accessible talk by Baez on YouTube. 

\item All of the decorated cospan formalism can be found in Brendan Fong's thesis, \href{https://arxiv.org/abs/1609.05382}{``The Algebra of Open and Interconnected Systems.''} As part of the ACT workshop, Jonathan Lorand and Fabrizio Genovese wrote about Fong's thesis in an article titled \href{https://golem.ph.utexas.edu/category/2018/02/hypergraph\_categories\_of\_cospa.html}{``Hypergraph Categories of Cospans''} on the $n$-Category Caf{\'e}.
\item The framework of chemical reaction networks is also being used to model ATP coupling! This project was birthed during the ACT workshop. For more, take a look at  \href{https://johncarlosbaez.wordpress.com/2018/06/27/coupling-through-emergent-conservation-laws-part-1}{``Coupling Through Emergent Conservation Laws ''} on Baez's Azimuth blog.
\end{itemize}

\newpage
\textbf{For more on the work of Coecke, Sadrzadeh, and Clark}:
\begin{itemize}
	\item There is a delightful 5-minute YouTube video (made in the style of \href{https://www.youtube.com/user/minutephysics}{Minute Physics}!) called \href{https://www.youtube.com/watch?v=99keybEZN4g}{``How Quantum Theory Can Help Understanding Natural Language''} that does an excellent job of explaining the ideas behind DisCoCat in a non-technical way.

	\item The \textit{yellow banana} example of the previous section was just a toy example meant to showcase the functor of the DisCoCat model of meaning. But this is a document on \textit{applied} category theory, and so you'd surely like to see some applications! For empirical data arising from actual implementations of the DisCoCat model, take a look at:
		\begin{itemize}
			\item \href{http://www.aclweb.org/anthology/D11-1129}{``Experimental support for a categorical compositional distributional model of meaning''} by Edward Grefenstette and Mehrnoosh Sadrzadeh 
			\item \href{https://www.cs.ox.ac.uk/files/5725/karts_sadr_emnlp.pdf}{``Prior disambiguation of word tensors for constructing sentence vectors''} by Dimitri Kartsaklis and Mehrnoosh Sadrzadeh
			\item \href{https://www.frontiersin.org/articles/10.3389/fphy.2017.00018/full}{``Quantization, Frobenius and bi-algebras from the categorical framework of quantum mechanics to natural language semantics''} by Mehrnoosh Sadrzadeh
		\end{itemize}

	\item The DisCoCat model was also featured on the \href{https://golem.ph.utexas.edu/category/}{$n$-Category Caf{\'e}} as part of the \href{http://www.appliedcategorytheory.org/school/}{Applied Category Theory Workshop.} In fact, it was featured \textit{twice}! The two blog posts are:
		\begin{itemize}
			\item \href{https://golem.ph.utexas.edu/category/2018/02/linguistics_using_category_the.html}{``Linguistics Using Category Theory''} by Corey Griffith and Jade Master. This wonderfully written article gives another recap of the DisCoCat model, as well as another toy example. Check out the comment section, too. There are some nice discussion going on there!
			\item \href{https://golem.ph.utexas.edu/category/2018/03/cognition_convexity_and_catego.html}{``Cognition, Convexity, and Category Theory''} by Brad Theilman and me. This blog post is a summary of \href{https://arxiv.org/pdf/1703.08314.pdf}{``Interacting Conceptual Spaces I,''} a paper in which the authors change the semantics category of the DisCoCat model from $\cat{FVect}$ to something that attempts to model human cognition more closely, namely convex spaces! This is a great example of tweaking the \textit{semantics} part of \textit{functorial \textbf{semantics}}.
		\end{itemize}

	\item The last suggested resource is unrelated to DisCoCat, but it's still in the vein of machine learning + applied category theory and so I thought I'd share: Did you know that \href{https://arxiv.org/pdf/1711.10455.pdf}{backpropagation is a functor}?
\end{itemize}

\newpage
\section{But Wait! There's More...}
There's much more to applied category theory---I've only presented a very tiny subset of hand-selected ideas. But there's so much more to see, learn, and do! So to close out these notes, I'll leave you with a few more links where you can discover other \textit{themes, constructions,} and \textit{examples} of applied category theory. 
\begin{itemize}
	\item To start, there's the \href{http://www.appliedcategorytheory.org/}{main Applied Category Theory webpage}, which has
		\begin{itemize}
			\item a description of the \href{http://www.appliedcategorytheory.org/workshops/}{2018 workshop} that took place at the \href{https://www.lorentzcenter.nl/lc/web/2018/969/info.php3?wsid=969&venue=Oort}{Lorentz Center} in Leiden, Netherlands
			\item a call for the \href{http://www.appliedcategorytheory.org/act-2019/}{2019 workshop}, to be hosted by Bob Coecke in Oxford.
		\end{itemize}

	\item Jelle Herold and the folks at \href{https://statebox.org/}{Statebox} filmed most of the 2018 workshop talks, and you can watch them here: \href{https://statebox.org/events/act-leiden.html}{https://statebox.org/events/act-leiden.html}. Speakers include Samson Abramsky, John Baez, Bob Coecke,  Kathryn Hess, Aleks Kissinger, Tom Leinster, David Spivak, and many more!

	\item Back in March 2018, there was an applied category theory workshop hosted at the National Institute of Standards and Technology. Slides and videos of the talk can be found here:\\
	\href{http://www.appliedcategorytheory.org/nist-workshop-slides/}{http://www.appliedcategorytheory.org/nist-workshop-slides/}

	\item There is also \href{https://arxiv.org/pdf/1803.05316.pdf}{\textit{Seven Sketches in Compositionality}} (subtitle: ``An Invitation to Applied Category Theory'') by Brendan Fong and David Spivak. I've referenced this book several times already, but that's because it's such a gem! (I was sold just after reading the preface.) It's a delightful and insightful introduction to more themes, more constructions, and more examples within applied category theory. Even better, no prior knowledge of category theory is assumed. The book is based on a MIT course the authors taught. You can find videos of their lectures here: \href{http://math.mit.edu/~dspivak/teaching/sp18/}{http://math.mit.edu/~dspivak/teaching/sp18/}
	
	\item And as if all of these great resources weren't enough, John Baez is running a \textit{free online course} on applied category theory. Participants have been working through the \textit{Seven Sketches} book. The lectures and ensuing discussions are a treasure trove of exciting mathematics: \href{https://forum.azimuthproject.org/categories/applied-category-theory-course}{https://forum.azimuthproject.org/categories/applied-category-theory-course}. Also be sure to take a look at the ``applied category'' tag on his blog, \href{https://johncarlosbaez.wordpress.com/?s=applied+category}{Azimuth}.
\end{itemize}

\end{document}